# An Introduction to Wishart Matrix Moments

Adrian N. Bishop[1], Pierre Del Moral[2] and Angèle Niclas[3]

[1] *University of Technology Sydney (UTS) and CSIRO, Australia; adrian.bishop@uts.edu.au*

[2] *INRIA, Bordeaux Research Center, France; pierre.del-moral@inria.fr*

[3] *École Normale Supérieure de Lyon, France*

---

ABSTRACT

These lecture notes provide a comprehensive, self-contained introduction to the analysis of Wishart matrix moments. This study may act as an introduction to some particular aspects of random matrix theory, or as a self-contained exposition of Wishart matrix moments.

Random matrix theory plays a central role in statistical physics, computational mathematics and engineering sciences, including data assimilation, signal processing, combinatorial optimization, compressed sensing, econometrics and mathematical finance, among numerous others. The mathematical foundations of the theory of random matrices lies at the intersection of combinatorics, non-commutative algebra, geometry, multivariate functional and spectral analysis, and of course statistics and probability theory. As a result, most of the classical topics in random matrix theory are technical, and mathematically difficult to penetrate for non-experts and regular users and practitioners.

The technical aim of these notes is to review and extend some important results in random matrix theory in the specific context of real random Wishart matrices. This special class of Gaussian-type sample covariance matrix plays an important role in multivariate analysis and in statistical theory. We derive non-asymptotic formulae for the full matrix moments of real valued Wishart random matrices. As a corollary, we derive and extend a number of spectral and trace-type results for the case of non-isotropic Wishart random matrices. We also derive the full matrix moment analogues of some classic spectral and trace-type moment results. For example, we derive semi-circle and Marchencko–Pastur-type laws in the non-isotropic and full matrix cases. Laplace matrix transforms and matrix moment estimates are also studied, along with new spectral and trace concentration-type inequalities.

---

# 1

# Introduction

Let $X$ be a centered Gaussian random column vector with covariance matrix $P$ on $\mathbb{R}^r$, for some dimension parameter $r \geq 1$. The rescaled sample covariance matrix associated with $(N+1)$ independent copies $X_i$ of $X$ is given by the random matrix

$$P_N = \frac{1}{N} \sum_{1 \leq i \leq N+1} \left(X_i - m^N\right) \left(X_i - m^N\right)'$$

with the sample mean

$$m^N := \frac{1}{N+1} \sum_{1 \leq i \leq N+1} X_i$$

Here, $(.)'$ denotes the transpose operator. The random matrix $P_N$ has a Wishart distribution with $N$ degrees of freedom and covariance matrix $N^{-1}P$ (a.k.a. the scale matrix). When $N \geq r$,the distribution of the Wishart matrix $P_N$ on the cone of symmetric positive definite matrices is defined by

$$\text{Probability}(P_N \in dQ)$$

$$= \frac{\det(Q)^{(N-r-1)/2}}{2^{Nr/2}\Gamma_r(N/2)\det(P/N)^{N/2}} \exp\left[-\frac{1}{2}\text{Tr}\left((P/N)^{-1}Q\right)\right] \gamma(dQ)$$

where $\det(Q)$ denotes the determinant of $Q$ and $\gamma(dQ)$ is the Lebesgue measure on the cone of symmetric positive definite matrices, and $\Gamma_r$ is the multivariate gamma function

$$\Gamma_r(z) = \pi^{r(r-1)/4} \prod_{1\leq k\leq r} \Gamma\left(z - \frac{k-1}{2}\right)$$

We also have the equivalent formulations

$$P_N \stackrel{law}{=} N^{-1} \sum_{1\leq i\leq N} \mathbb{X}_i = \mathcal{X}\mathcal{X}'$$

with the $(r \times N)$-random matrix $\mathcal{X}$ defined by

$$\mathcal{X} = \frac{1}{\sqrt{N}} [X_1, \ldots, X_N]$$

In the above display, $\mathbb{X}_i$ stand for $N$ independent copies of the rank one random matrix $\mathbb{X} = XX'$. The superscript $(.)'$ denotes the transposition operation.

Random matrices, sample covariance matrices, and more specifically Wishart random matrices, play a role in finance and statistics, physics, and engineering sciences. Their interpretation depends on the application model motivating their study.

For example, in Bayesian inference, Wishart matrices often represent the prior precision matrix of multivariate Gaussian data sets. In this context, the posterior distribution of the random covariance given the multivariate-normal vector is again a Wishart distribution with a scale matrix that depends on the measurements. In other words, Wishart distributions are conjugate priors of the inverse covariance-matrix of a multivariate normal random vector [8, 64].

In multivariate analysis and machine learning, the vectors $X_i$ may represent some statistical data such as image, curves and text data. In this case, $P$ may be defined in terms of some covariance function as in Gaussian processes [72]. As its name indicates, the sample covariance matrix $P_N$ attempts to capture the shape of the data; such as the spread around their sample mean as well as the sample correlation between the features dimensions. Principal component analysis and related techniques amount to finding the eigenvalues and the corresponding eigenvectors of sample covariance matrices. The largest eigenvalues represents the dimensions with the strongest correlation in the data set. Expressing the data on the eigenvectors associated with the largest eigenvalues is often used to compress high dimensional data. For a more thorough discussion on this subject we refer to the articles [2, 7, 46, 69], as well as the monographs [8, 72, 63] and the references therein.

In the context of multiple-input multiple-output systems, more general random matrices may be related to the channel gain matrix [77, 99]. Similarly, the covariance matrix in Gaussian process-based inference may be considered a random matrix defined by the particular covariance structure [72]. In data assimilation problems and filtering theory, non-independent sample covariance matrices arise as the control gain in ensemble (Kalman-type) filters; see e.g. [11, 23, 9] and the references therein. Similar (non-independent) sample covariance matrices may be computed with the particles in classical Markov Chain Monte Carlo and sequential Monte Carlo methods [21]; and in this case often represent the uncertainty in an estimation theoretic sense. In finance, sample covariance matrices arise in risk management and asset allocation; e.g. random matrices may represent the correlated fluctuations of assets [57, 13, 27].

Because of their practical importance, we may illustrate the above specific model via the so-called Wishart process. Consider a time-varying linear-Gaussian diffusion of the following form,

$$dX(t) = A(t)\,X(t)\;dt + R(t)^{1/2}\,dB(t) \tag{1.1}$$

where $B(t)$ is an $r$-dimensional Brownian motion, $X_0$ is a $r$-dimensional Gaussian random variable with mean and variance $(\mathbb{E}(X_0), P_0)$, independent of $B(t)$, and $A(t) \in \mathbb{R}^{r\times r}$, and $R(t) > 0$ is a positive definite symmetric matrix. The covariance matrices

$$P(t) = \mathbb{E}\left([X(t) - \mathbb{E}(X(t))]\,[X(t) - \mathbb{E}(X(t))]'\right)$$

satisfy the (linear) matrix-valued differential equation,

$$\partial_t P(t) = A(t)\,P(t) + P(t)\,A(t)' + R(t)$$

The solution of the preceding equation is given easily via the transition/fundamental matrix defined by $A(t)$. More precisely, the solution of the above equation is given by the formula

$$P(t) = e^{\oint_0^t A(s)ds}\ P(0)\ \left[e^{\oint_0^t A(s)ds}\right]' + \int_0^t\ e^{\oint_s^t A(u)du}\ R(s)\ \left[e^{\oint_s^t A(u)du}\right]'\ ds$$

In the above display, $\mathcal{E}_{s,t} := e^{\oint_s^t A(u)du}$ denotes the matrix exponential semigroup, or the transition matrix, defined by

$$\partial_t\,\mathcal{E}_{s,t} = A_t\,\mathcal{E}_{s,t} \quad \text{and} \quad \partial_s\,\mathcal{E}_{s,t} = -\mathcal{E}_{s,t}\,A_s \quad \text{with} \quad \mathcal{E}_{s,s} = I$$

For time homogeneous models $(A(t), R(t)) = (A, R)$ the above formula reduces to

$$P(t) = e^{tA}P(0) + \int_0^t\ e^{(t-s)A}\ R\ e^{(t-s)A'}\ ds$$

The rescaled sample covariance matrices associated with $(N+1)$ independent copies $(X_i(t))_{1\leq i\leq N+1}$ of the process $X(t)$ are defined by

$$P_N(t) := \frac{1}{N}\sum_{1\leq i\leq N+1}\left[X_i(t) - m^N(t)\right]\left[X_i(t) - m^N(t)\right]'$$

with the sample mean

$$m^N(t) := \frac{1}{N+1}\sum_{1\leq i\leq N+1} X_i(t)$$

Up to a change of probability space, the process $P_N(t)$ satisfies the matrix diffusion equation

$$dP_N(t) = A(t)\,P_N(t) + P_N(t)\,A(t)' + R(t) + \frac{1}{\sqrt{N}}\ M_N(t)$$

with the matrix-valued martingale

$$dM_N(t) = P_N(t)^{1/2}\ d\mathcal{W}(t)\ R(t)^{1/2} + R(t)^{1/2}\ d\mathcal{W}(t)\ P_N(t)^{1/2}$$

where $\mathcal{W}_t$ denotes an $(r\times r)$-matrix with independent Brownian entries. The above diffusion coincides with the Wishart process considered in [15]. When $r=1$ this Wishart model coincides with the Cox–Ingersoll–Ross process (a.k.a. squared Bessel process) introduced in [19]. For a more detailed discussion on Wishart processes and related affine diffusions, we refer to the articles [20, 40, 58], and the references therein.

The preceding exposition is by no means exhaustive of applications of random, and more specifically Wishart, matrices and we point to [26, 30, 59, 89, 97, 98] for further applications and motivators. The typical technical questions arising in practice revolve around the calculation of the spectrum distribution, and the corresponding eigenvector distribution of these random matrices.

The analysis of Wishart matrices started in 1928 with the pioneering work of J. Wishart [100]. Since this, the theory of random matrices has been a fruitful contact point between statistics, pure and applied probability, combinatorics, non-commutative algebra, as well as differential geometry and functional analysis.

The joint distribution of the eigenvalues of real valued Wishart matrices is only known for full rank and isotropic models; that is when the sample size is greater than the dimension, and the covariance matrix $P \propto I$ is proportional to the identity matrix $I$; see for instance [34, 62]. In this situation, the matrix of random eigenvectors is uniformly distributed on the manifold of unitary matrices equipped with the Haar measure. In this context, marginal distributions for these uncorrelated models can also be computed in a tractable form only for the smallest and the largest eigenvalues. Sophisticated integral formulae for the marginal distribution of intermediate eigenvalues are provided by Zanella–Chiani–Win [104]. Upper bounds on the marginal distribution of the ordered eigenvalues are given in [67].

The cumulative distribution of the largest eigenvalue of real Wishart random matrices can be expressed explicitly in terms of the hypergeometric function of a matrix argument. These functionals can also be described in terms of zonal polynomials. The smallest and largest eigenvalue distributions can also be expressed in terms of Tricomi functions [26]. For a detailed discussion on these objects we refer the reader to the book of Muirhead [62]. As shown in [45], these hypergeometric functions depends on alternating series involving zonal polynomials which converge very slowly even in low dimensions. Some explicit calculations for $r = 1, 2, 3$ can be found in [88].

Non-necessarily isotropic Wishart models can be considered if we restrict our attention to linear transforms and other trace-type mathematical objects. We refer to the articles of Letac and his co-authors [29, 47, 48, 91] and the tutorial [49]. See also [35, 50] for a description of the joint distribution of traces of Wishart matrices.

To bypass the complexity of finding computable and tractable closed form solutions, one natural and common method for obtaining useful information is to derive limiting distributions as the dimension tends to $\infty$. In this context, one can analyze the convergence of the histogram of the eigenvalues when the dimension tends to $\infty$. This approach is central in random matrix theory. We refer the reader to the pioneering article by E. Wigner [98] published in 1955, the lectures notes of A. Guionnet [30], the research monographs by M.L. Mehta [59] and T. Tao [89], and the references therein. This commonly used limiting theory has some drawbacks. Firstly, as the name suggests, these limiting techniques cannot capture nor control the non-asymptotic fluctuations arising in practical problems. Moreover, the limiting techniques developed in the literature often yield information only on the limiting behaviour of trace or spectral-type properties of random matrix powers. In addition, in the context of Wishart matrices, these limiting spectral-type techniques only apply to asymptotically isotropic-type models. To be more precise, the convergence analysis relies on strong hypotheses on the bias and the variance of the random matrix entries which are satisfied only for Wishart matrices with a covariance matrix close to the identity (up to some ad-hoc scaling factor). When $P \neq I$, the distribution of the eigenvalues and the corresponding eigenvectors is much more involved. The distribution of the sample eigenvalues depends on sophisticated Harish-Chandra integrals [31].

Importantly, all the spectral and trace-type approaches discussed above (whether in the limit or not) give no information on the random matrix moments themselves, but rather on their eigenvalues or trace, etc. Conversely, in many practical situations, such as in data assimilation theory and signal processing (e.g. ensemble Kalman filter theory [11, 23, 9] and particle filtering [21]), we are typically interested in the direct analysis of full matrix moments of interacting-type (non-independent) sample covariance matrices. This study concerns a step in this latter direction. Specifically, we derive formulae for the full matrix moments of real valued Wishart random matrices. As a corollary, we derive and extend a number of spectral and trace-type results for the case of non-isotropic Wishart random matrices. Laplace matrix transforms and matrix moment estimates are also studied, along with new spectral and trace concentration-type inequalities.

## 1.1 Organisation

Section 1.2 is concerned with the description of real random Wishart matrices and their fluctuation analysis. We also review some central result in random matrix theory, such as the semi-circle law and the Marchenko–Pastur law for isotropic Wishart matrices.

Section 1.3 provides a brief description of the main results of these notes. We provide a closed form Taylor-type formula to compute the matrix moments of $P_N$, and its fluctuations defined in (1.2), w.r.t. the precision parameter $1/N$. We also present the full matrix version of the semi-circle law and the Marchenko–Pastur law for non-necessarily isotropic Wishart matrices. Non-asymptotic matrix moments and exponential Laplace transforms are also provided. The last part of the section is concerned with exponential concentration inequalities for operator norms of fluctuation matrices and the eigenvalues of sample covariance matrices.

The rest of these lecture notes is concerned with the precise statement and proof of the theorems in Section 1.3. Some auxiliary outcomes and discussion surrounding these results are also given.

Section 2 reviews some useful mathematical background on Laplace and exponential inequalities, matrix norms, spectral analysis, tensor products, Fréchet derivatives, and fluctuation-type results. This section also contains a brief review of non-crossing partitions, Catalan, Narayana, and Riordan numbers, Bell polynomials, and Murasaki and circular-type representations of non-crossing partitions. The last part of this section discusses Fréchet differentiable functionals on matrix spaces and Taylor-type approximations.

Section 3 concerns closed form polynomial formulae for computing the matrix moments of $P_N$, and its fluctuations in (1.2), in terms of the precision parameter $1/N$ and some partition-type matrix moments. The isotropic semi-circle law (1.8) and the Marchenko–Pastur law (1.10) are simple consequences of these matrix moments expansions. New matrix versions of the semi-circle and the Marchenko–Pastur law for non-necessarily isotropic Wishart matrices are discussed in Section 3.2 and in Section 3.3.

Section 4 is concerned with some matrix moment estimates and Laplace matrix transforms of the fluctuation matrix. Spectral and trace-type concentration-type inequalities are discussed in Section 5. Section 6 is dedicated to the proof of the main theorems.

An appendix is also given containing the proof of a number of technical results required throughout these notes.

## 1.2 Description of the models

We recall the multivariate central limit theorem

$$P_N = P + \frac{1}{\sqrt{N}}\,\mathcal{H}_N$$

with

$$\mathcal{H}_N \overset{law}{:=} \frac{1}{\sqrt{N}} \sum_{1\leq i\leq N} (\mathbb{X}_i - P) \hookrightarrow_{N\to\infty} \mathcal{H} \tag{1.2}$$

where $\mathcal{H}$ is a symmetric $(r\times r)$-matrix with centered Gaussian entries equipped with a symmetric Kronecker covariance structure

$$(\mathcal{H}\otimes\mathcal{H})^\sharp = 2\,(P \mathbin{\widehat{\otimes}} P) = \mathbb{E}\left[(\mathcal{H}_N\otimes\mathcal{H}_N)^\sharp\right] \tag{1.3}$$

where $(A\otimes B)^\sharp$ and $(A \mathbin{\widehat{\otimes}} B)$ are the entry-wise and the symmetric tensor product of matrices $A$ and $B$. These products are defined at the beginning of Section 2.3.

A detailed discussion on the fluctuation result (1.2) can be found in [41]; see also [14] for non-necessarily Gaussian variables. The fluctuation result (1.2) can also be deduced from the Laplace matrix transform estimates stated in theorem 1.5 and corollary 4.5.

Combining a perturbation analysis with the continuous mapping theorem, the central limit result (1.2) can be used to analyze the fluctuation of smooth matrix functionals of the sample covariance matrix.

Roughly speaking, given some smooth Fréchet differentiable mapping $\Upsilon : \mathcal{S}_r \mapsto \mathcal{B}$ from symmetric matrices $\mathcal{S}_r$ to some Banach space $\mathcal{B}$, we have the Taylor expansion

$$\begin{aligned}\mathcal{H}_N^\Upsilon &:= \sqrt{N}\ [\Upsilon(P_N) - \Upsilon(P)]\\ &= \nabla\Upsilon(P)\cdot\mathcal{H}_N + \frac{1}{2\sqrt{N}}\ \nabla^2\Upsilon(P)\cdot(\mathcal{H}_N,\mathcal{H}_N) + \ldots\end{aligned}$$

Using the unbiasedness properties of the sample covariance matrix, the second order term gives the bias of the estimate $\Upsilon(P_N)$; that is we have that

$$\mathbb{E}\left[\mathcal{H}_N^\Upsilon\right] = \frac{1}{2\sqrt{N}}\ \mathbb{E}\left[\nabla^2\Upsilon(P)\cdot(\mathcal{H},\mathcal{H})\right] + \mathrm{O}\left(\frac{1}{N}\right)$$

Equivalently, we have

$$\mathbb{E}\left[\Upsilon(P_N)\right] = \Upsilon(P) + \frac{1}{2N}\ \mathbb{E}\left[\nabla^2\Upsilon(P)\cdot(\mathcal{H},\mathcal{H})\right] + \mathrm{O}\left(\frac{1}{N^{3/2}}\right)$$

For a more precise statement and several illustrations we refer the reader to Section 2.6, theorem 2.2. For instance, for power functions $\Upsilon_n(Q) := Q^n$, for any $1 \leq m \leq n$ we have

$$\begin{aligned}&\mathbb{E}\left[\nabla^m\Upsilon_n(P)\cdot\mathcal{H}_N^{\otimes m}\right]\\ &= \ m!\sum_{0\leq i_1<\ldots<i_m\leq n}\mathbb{E}\left(\prod_{1\leq k\leq m}\left[P^{i_k-i_{k-1}-1}\ \mathcal{H}_N\right]\right)\\ &= (n)_m\ \mathbb{E}(\mathcal{H}_N^m)\quad\text{when}\quad P=I\end{aligned} \tag{1.4}$$

with the Pochhammer symbol $(n)_m := n!/(n-m)!$, and the convention $(i_0, i_{m+1}) = (0, n)$. In this context, the $m$-moments of the fluctuation matrices $\mathcal{H}_N$ represents the mean-error of order $m$. This property also holds for rational powers. For instance, we have the non-asymptotic estimate on the Frobenius norm,

$$\Big\| \mathbb{E}\left(\sqrt{P_N}\right) - \sqrt{P} + \frac{1}{2N}\sqrt{P}\left[\frac{1}{4}I + \int_0^\infty t\,\mathrm{Tr}\left(P\,e^{-t\sqrt{P}}\right)\,e^{-t\sqrt{P}}\,dt\right]\Big\|_F \leq \frac{r}{4}\,\frac{1}{N\sqrt{N}}\,\lambda_{min}(P)^{-5/2}\,\left[\mathrm{Tr}(P^2)+\mathrm{Tr}(P)^2\right] \tag{1.5}$$

The proof of this assertion is provided in Section 2.6.

To summarise the consequences of the preceding discussion, to analyze these approximations at any order it is therefore necessary to be able to compute the $m$-moments of the fluctuation matrices $\mathcal{H}_N$.

Gaussian approximation techniques also require one to estimate the fluctuations of the moment $\mathbb{E}(\mathcal{H}_N^m)$ around those of $\mathbb{E}(\mathcal{H}^m)$ given in terms of the limiting Gaussian matrix $\mathcal{H}$, and with respect to the sample size parameter. Moreover, one often wants to control the behavior of these objects when the dimension parameter tends to $\infty$.

The limiting random matrix model $\mathcal{H}$ discussed above is closely related to Gaussian orthogonal ensembles arising in random matrix theory. To be more precise, we can check that

$$\mathcal{H} \overset{law}{=} P^{1/2}\left(\frac{\mathcal{W}+\mathcal{W}'}{\sqrt{2}}\right)P^{1/2} \tag{1.6}$$

where $\mathcal{W} = (\mathcal{W}_{i,j})_{1\leq i,j\leq r}$ is a matrix of independent, centred Gaussian elements of unit variance.

When $P = I$, the random matrix $\mathcal{H}$ introduced in (1.6) reduces to a Gaussian orthogonal ensemble. In this situation, we have

$$r^{-2}\,\mathbb{E}\left(\mathrm{Tr}\left[\mathcal{H}^2\right]\right) = 1 + r^{-1} \quad \text{and} \quad r^{-3}\,\mathbb{E}\left(\mathrm{Tr}\left[\mathcal{H}^4\right]\right) = 2 + 5\,r^{-2} + 5\,r^{-1} \tag{1.7}$$

The trace of the higher moments $\mathbb{E}\,(\mathcal{H}^n)$ can be estimated using the semi-circle law (1.8) in large dimensions. This celebrated limiting result is proved using the convergence of moment property

$$r^{-1}\,\mathbb{E}\left(\mathrm{Tr}\left(\left[\frac{\mathcal{H}}{\sqrt{r}}\right]^n\right)\right) = 1_{2\mathbb{N}}(n)\;C_{n/2} + O\left(1/r\right) \qquad (1.8)$$

$$= \frac{1}{2\pi}\int_{-2}^{2} x^n\,\sqrt{4-x^2}\;dx + O\left(1/r\right)$$

for any $n \geq 1$, with the Catalan numbers

$$C_n := \frac{1}{n+1}\binom{2n}{n} \qquad (1.9)$$

A proof of the above assertion via Wick's theorem, including detailed reference pointers is given in [30, Section 1.4], see also [60, Chapter 1].

When $P = I$ and $N = r/\rho$ for some parameter $\rho > 0$, another important result is the Marchenko–Pastur law

$$\lim_{r\to\infty} r^{-1}\,\mathbb{E}\left[\mathrm{Tr}\left(P_N^n\right)\right]$$

$$= \sum_{0\leq m<n} \frac{\rho^m}{m+1}\binom{n}{m}\binom{n-1}{m}$$

$$= \int_{a_-(\rho)}^{a_+(\rho)} x^n\left[\left(1-\frac{1}{\rho}\right)_+\delta_0(dx) + \frac{1}{2\pi\rho\,x}\sqrt{\left[a_+(\rho)-x\right]\left[x-a_-(\rho)\right]}\;dx\right]$$

with the parameters

$$a_-(\rho) := (1-\sqrt{\rho})^2 \quad \text{and} \quad a_+(\rho) := (1+\sqrt{\rho})^2$$

The proof of the above integral formula can be found in [101, lemma 5.2], see also [28, 65] and the pioneering article by Vladimir Marchenko and Leonid Pastur [56]. A new proof of this result follows from the full matrix version of the Marchenko–Pastur law given in corollary 3.3 of these notes.

When $P \neq I$, formula (1.6) can be combined with Isserlis' theorem [33] (or Wick's theorem [97]) to compute the matrix moments of the random matrix $\mathcal{H}$. For instance, we have $\mathbb{E}(\mathcal{H}^{2n+1}) = 0$, for any $n \geq 0$. After some lengthy combinatorial computations we also find the

matrix polynomials

$$\begin{aligned}\mathbb{E}\left(\mathcal{H}^2\right) &= P^2 + \mathrm{Tr}(P)\,P \\ \mathbb{E}\left(\mathcal{H}^4\right) &= 5\,P^4 + 3\,\mathrm{Tr}(P)\,P^3 + \left[\mathrm{Tr}(P^2) + \mathrm{Tr}(P)^2\right]\,P^2 \\ &\qquad + \left[\mathrm{Tr}(P^3) + \mathrm{Tr}(P)\,\mathrm{Tr}(P^2)\right]\,P\end{aligned}$$

Although these matrix moments are given by a closed form formula, their complex combinatorial structure cannot be used in simple calculations. For example, the calculation of $\mathbb{E}\left(\mathcal{H}^{2n}\right)$ requires the matrix moments associated with $2^{-n}\,(2n)!/n!$ partitions over $[2n] := \{1,\ldots,2n\}$ with $n$-blocks. For $n = 4$, more than one hundred moments need to be computed. The computational complexity to numerically compute the central moments of the multivariate normal distribution is discussed in [68]; see also [3, p. 49], [38, proposition 1], [62, p. 46], and the matrix derivative formula in [92].

The above formulae also show that we cannot expect to have a semi-circle-type law as in (1.7) for any covariance matrix. Different types of behaviour can be expected depending on the behavior of the eigenvalues of $P$ w.r.t. the dimension parameter $r$. For instance, if the largest eigenvalue is $\lambda_1(P) = r$, we have $1 \leq r^{-1}\,\mathrm{Tr}(P) \leq 2$ but

$$\mathbb{E}\left(\mathrm{Tr}\left[\mathcal{H}^4\right]\right) \geq 5\,r^4 \quad\Longrightarrow\quad r^{-1}\,\mathbb{E}\left(\mathrm{Tr}\left(\left[\frac{\mathcal{H}}{\sqrt{r}}\right]^4\right)\right) \longrightarrow_{r\to\infty}\ \infty$$

## 1.3 Statement of some main results

One of the main objectives of these lecture notes is to analyze the properties of real Wishart matrix moments. Let $\mathcal{P}_n$ be the set of all partitions $\pi$ of $[n] := \{1,\ldots,n\}$, $\mathcal{P}_{n,m} \subset \mathcal{P}_n$ be the subset of all partitions with $m$ blocks $\pi_1 \leq \ldots \leq \pi_m$ ordered in a canonical way w.r.t. their smallest element.

Let $\mathcal{Q}_n \subset \mathcal{P}_n$ and $\mathcal{Q}_{n,m} \subset \mathcal{P}_{n,m}$ be the subset of partitions without the singleton. Also let $\alpha^\pi := \sum_{1\leq i\leq m} i\ 1_{\pi_i}$. In other words, $\alpha^\pi(i)$ is the index of the block of $\pi$ containing index $i$.

The $\pi$-matrix moments $M_\pi^{[Q]}(P)$ and $M_\pi^{\circ,[Q]}(P)$ associated with some collection of $(r \times r)$ matrices $(Q_i)_{i\geq 1}$ are defined by

$$M_\pi^{[Q]}(P) := \mathbb{E}\left([\mathbb{X} - P]_\pi^Q\right) \quad \text{and} \quad M_\pi^{\circ,[Q]}(P) := \mathbb{E}\left(\mathbb{X}_\pi^Q\right) \tag{1.10}$$

with the random matrices

$$\mathbb{X}_\pi^Q := \prod_{1\leq i\leq n}\left[Q_i \mathbb{X}_{\alpha^\pi(i)}\right]$$

and

$$[\mathbb{X} - P]_\pi^Q := \prod_{1\leq i\leq n}\left[Q_i(\mathbb{X}_{\alpha^\pi(i)} - P)\right]$$

We also consider the matrix moments

$$M_{n,m}^{[Q]}(P) := \sum_{\pi\in\mathcal{Q}_{n,m}} M_\pi^{[Q]}(P) \quad \text{and} \quad M_{n,m}^{\circ,[Q]}(P) := \sum_{\pi\in\mathcal{P}_{n,m}} M_\pi^{\circ,[Q]}(P)$$

Our first main result provides polynomial formulae w.r.t. the precision parameter $1/N$.

**Theorem 1.1.** For any collection of matrices $Q_n$, and any $2N \geq n \geq 1$, we have the polynomial formulae

$$\mathbb{E}\left[(Q_1\mathcal{H}_N)\ldots(Q_n\mathcal{H}_N)\right] = \sum_{1\leq m\leq \lfloor n/2\rfloor} \frac{1}{N^{n/2-m}}\, \partial_{n,m}^{[Q]}(P) \tag{1.11}$$

with

$$\partial_{n,m}^{[Q]}(P) := \sum_{m\leq l\leq \lfloor n/2\rfloor} s(l,m)\ M_{n,l}^{[Q]}(P)$$

In addition, we have

$$\mathbb{E}\left[(Q_1P_N)\ldots(Q_nP_N)\right] = \sum_{1\leq m\leq n} \frac{1}{N^{n-m}}\, \partial_{n,m}^{\circ,[Q]}(P) \tag{1.12}$$

with

$$\partial_{n,m}^{\circ,[Q]}(P) = \sum_{m\leq l\leq n} s(l,m)\ M_{n,l}^{\circ,[Q]}(P)$$

In the above displayed formulae, $s(l,m)$ are the Stirling numbers of the first kind.

For the detailed discussion of these matrix moments, including several corollaries and examples, we refer to Section 3.1; see for example theorem 3.1 when $Q_i = I$.

To simplify notation, for homogeneous models $Q_i = I$ and sequences of matrices $P : r \mapsto P(r)$ we suppress the indices $(.)^{[I]}$ and $r$, and write

$$(\partial_{n,m}(P), \partial^{\circ}_{n,m}(P), M_{\pi}(P), M^{\circ}_{\pi}(P), M_{n,l}(P), M^{\circ}_{n,l}(P))$$

instead of

$$(\partial^{[I]}_{n,m}(P(r)), \partial^{\circ,[I]}_{n,m}(P(r)), M^{[I]}_{\pi}(P), M^{\circ,[I]}_{\pi}(P), M^{[I]}_{n,l}(P(r)), M^{\circ,[I]}_{n,l}(P(r)))$$

The polynomial formula (1.11) differs from the invariant moments which can be derived using the algorithm presented in [47]. In the latter, the authors provide an elegant spectral technique to interpret these moments in terms of spherical polynomials and matrix-eigenfunctions of Wishart integral operators; see [47, e.g. proposition 4.3]. A drawback of this spectral method is that it requires one to diagonalize and invert complex combinatorial matrices. It is difficult to use this technique to derive estimates w.r.t. the sample size parameter. Matrix moment formulae can also be derived from [73]. Nevertheless the resulting Isserlis-type decompositions will involve complex series of summations over pair partitions.

Beside the fact that the matrix moments $M^{[Q]}_{n,l}(P)$ can be computed using Isserlis' theorem, to be the best of our knowledge no explicit and closed form polynomial formulae in terms of $P$ are known. In the further development of these notes, we provide estimates of the fluctuation matrix moments w.r.t. the sample size in terms of the dominating term of the sum (1.11). A brief description of these estimates are provided in theorem 1.5 below.

To move one step further in our discussion we assume that $Q_i = I$ and $n = 2m$. In this situation the single dominating term in (1.11) is given by the central matrix moments

$$\partial_{2m,m}(P) = M_{2m,m}(P) \tag{1.13}$$

This implies that

$$\mathbb{E}\left[\mathcal{H}_N^{2m}\right] = M_{2m,m}(P) + \mathrm{O}\left(\frac{1}{N}\right) I$$

We say a partition is a crossing partition whenever we can find $i < j < k < l$ with $i, k$ in a block and $j, l$ in the other block. Let $\mathcal{N}_n \subset \mathcal{P}_n$ and $\mathcal{N}_{n,m} \subset \mathcal{P}_{n,m}$ be the subsets of non-crossing partitions.

We denote by $\Sigma_n(P)$ the matrix polynomial given by

$$\Sigma_n(P) \ := \ \mathrm{Tr}(P) \sum_{\pi\in\mathcal{N}_n} \left[\prod_{i\geq 0} \mathrm{Tr}\left(P^{1+i}\right)^{r_i(\pi)}\right] \frac{P^{|\pi_1|}}{\mathrm{Tr}(P^{|\pi_1|+1})} \tag{1.14}$$

In the above display, $r_0(\pi) := n - \sum_{i\geq 1} r_i(\pi)$ where $r_i(\pi)$ is the number of blocks of size $i \geq 1$ in the partition $\pi$.

Also let $\Sigma^{\circ}_{n,m}(P)$ be the matrix polynomial given by

$$\Sigma^{\circ}_{n,m}(P) \ := \ \sum_{\pi\in\mathcal{N}_{n,m}} \left[\prod_{i\geq 1} \mathrm{Tr}(P^i)^{r_i(\Xi(\pi))}\right] \frac{P^{\iota(\pi)}}{\mathrm{Tr}\left(P^{\iota(\pi)}\right)} \tag{1.15}$$

In the above display, $\iota(\pi)$ denotes the number of blocks visible from above in the Murasaki diagram associated with $\pi$; see Section 2.5 for examples. The partition $\Xi(\pi) \in \mathcal{N}_{n+1-m}$ is defined in terms of a circular representation of $\pi$. That is, firstly, we subdivide the $n$ arcs of $\pi \in \mathcal{N}_{n,m}$ by a new series of $n$ nodes placed clockwise. Then $\Xi(\pi)$ is the coarsest non-crossing partition of these nodes whose chords don't cross those of $\pi$. For a detailed description of the mapping $\Xi$, and examples, we refer Section 2.5; see e.g. (2.20).

Lets further assume that $P : r \mapsto P(r)$ is a collection of possibly random matrices satisfying for any $n \geq 1$ the almost sure convergence of the moments

$$r^{-1}\ \tau_n\left(P(r)\right) := r^{-1}\ \mathrm{Tr}(P(r)^n) \quad \longrightarrow_{r\to\infty} \quad \tau_n(P) \tag{1.16}$$

Also, let $\mathcal{H}_N(r)$ and $\mathcal{H}(r)$ be the random matrix model defined as in (1.2) and (1.6) by replacing $P$ by $P(r)$.

To simplify notation, we write

$$(\mathcal{H}, \mathcal{H}_N, M_{2n,n}(P), M^{\circ}_{n,m}(P), \Sigma_n(P), \Sigma^{\circ}_{n,m}(P))$$

instead of

$$(\mathcal{H}(r), \mathcal{H}_N(r), M_{2n,n}(P(r)), M^{\circ}_{n,m}(P(r)), \Sigma_n(P(r)), \Sigma^{\circ}_{n,m}(P(r)))$$

In this notation, the next theorem relates the matrix moments $(\Sigma_n(P), \Sigma^\circ_{n,m}(P))$ with the matrix moments $(M_{n,m}(P), M^\circ_{n,m}(P))$ and the ones of the Gaussian matrix $\mathcal{H}$.

**Theorem 1.2.** Let $P : r \mapsto P(r)$ be a collection of possibly random matrices satisfying the condition (1.16). In this situation, the central matrix moments $M_{2n,n}(P)$ coincide with the ones of the limiting Gaussian matrix. In addition, for any $n \geq m \geq 1$ we have the matrix moment estimates

$$\begin{aligned} M_{2n,n}(P) = \ & \mathbb{E}\left(\mathcal{H}^{2n}\right) \ = \ \Sigma_n(P) + \mathrm{O}\left(r^{n-1}\right) I \\ M^\circ_{n,m}(P) = \ & \Sigma^\circ_{n,m}(P) + \mathrm{O}\left(r^{n-m-1}\right) I \end{aligned} \tag{1.17}$$

For a proof and a more detailed discussion on these matrix moment relations we refer to Section 3.2 and Section 3.3; see e.g. theorem 3.4, theorem 3.5 and theorem 3.7.

The first line estimate in (1.17) is a consequence of the decomposition (3.1) and theorem 3.4. The second line estimate in (1.17) is a consequence of the estimates (3.18) and theorem 3.7.

Theorem 1.2 together with (1.13) yields the estimates

$$r^{-(n+1)}\ \mathbb{E}\left(\mathcal{H}^{2n}\right) = \overline{\Sigma}_n(P) + \mathrm{O}\left(r^{-1}\right) I_r$$

as well as

$$\mathbb{E}\left[\mathcal{H}_N^{2n}\right] = \mathbb{E}\left(\mathcal{H}^{2n}\right) + \mathrm{O}\left(N^{-1}\right) I$$

We also have

$$r^{-(n-m+1)} M^\circ_{n,m}(P) = \overline{\Sigma}^\circ_{n,m}(P) + \mathrm{O}\left(r^{-1}\right) I_r \quad \text{with} \quad I_r := r^{-1}\, I$$

with the matrix polynomials $(\overline{\Sigma}_n(P), \overline{\Sigma}^\circ_{n,m}(P))$ defined similarly to $\left(\Sigma_n(P), \Sigma^\circ_{n,m}(P)\right)$ but with the trace operator replaced by the normalized traces

$$\overline{T}\mathrm{r}(Q) := r^{-1}\,\mathrm{Tr}(Q)$$

A more refined estimate between $\mathbb{E}\left[\mathcal{H}_N^{2n}\right]$ and $\mathbb{E}\left(\mathcal{H}^{2n}\right)$ can be found later in theorem 1.5.

The first line assertion in (1.17) in theorem 1.2 provides a semicircle-type asymptotic theorem when the dimension parameter tends to $\infty$.

**Corollary 1.3.** Under the assumptions of theorem 1.2, we have the extended semi-circle law

$$r^{-1}\,\mathbb{E}\left(\mathrm{Tr}\left(\left[\frac{\mathcal{H}}{\sqrt{r}}\right]^{2n}\right)\right) = \sigma_n(P) + O\left(r^{-1}\right)$$

with

$$\sigma_n(P) := 2\sum_{\mu}\binom{n}{\mu_1\ \mu_2\ \ldots\ \mu_n}\tau_\mu(P) \quad \text{and} \quad \tau_\mu(P) := \prod_{1\leq i\leq n}\tau_i(P)^{\mu_i} \tag{1.18}$$

In the above display, the summation is taken over all collection of non-negative indices $\mu = (\mu_1, \ldots, \mu_n)$ such that

$$\sum_{1\leq i\leq n}\mu_i = n+1 \quad \text{and} \quad \sum_{1\leq i\leq n} i\mu_i = 2n \tag{1.19}$$

For instance, we have

$$\begin{aligned}
\sigma_1(P) &= \tau_1(P)^2\\
\sigma_2(P) &= 2\ \tau_1(P)^2\ \tau_2(P)\\
\sigma_3(P) &= 2\ \tau_1(P)^3\ \tau_3(P) + 3\ \tau_1(P)^2\ \tau_2(P)^2\\
\sigma_4(P) &= 2\ \tau_1(P)^4\ \tau_4(P) + 8\ \tau_1(P)^3\ \tau_2(P)\ \tau_3(P) + 4\ \tau_1(P)^2\ \tau_2(P)^3\\
\sigma_5(P) &= 2\ \tau_1(P)^5\ \tau_5(P) + 10\ \tau_1(P)^4\ \tau_2(P)\ \tau_4(P) + 5\ \tau_1(P)^4\ \tau_3(P)^2\\
&\quad + 20\ \tau_1(P)^3\ \tau_2(P)^2\ \tau_3(P) + 5\ \tau_1(P)^2\ \tau_2(P)^4
\end{aligned}$$

Observe that $\sigma_n(\alpha I) = C_n\ \alpha^{2n}$, for any $\alpha \geq 0$. These formulae can be checked combining (3.16) and (3.17) with corollary 3.6. Matrix-valued free probability techniques can also be used to recover the above matrix moment formula [66, 86, 87, 90, 95].

To the best of our knowledge the extended and matrix version of the semi-circle law stated in the above theorem have not been considered in the literature. See also Section 3.2.

Also recall that Carleman's condition

$$\sum_{n\geq 1}\sigma_n(P)^{-1/(2n)} = \infty$$

ensures the existence and uniqueness of a random variable with null odd moments and the $(2n)$-moments $\sigma_n(P)$ defined in corollary 1.3

(cf. [1, 16] and p. 296 in [78]). For instance, when $P = I$ we have

$$\sigma_n(I) = C_n \simeq \frac{2^{2n}}{n^{3/2}\sqrt{\pi}} \leq 2^{2n}$$

This implies that

$$\sigma_n(I)^{-1/(2n)} \geq \frac{1}{2} \quad \Longrightarrow \quad \sum_{n\geq 1} \sigma_n(I)^{-1/(2n)} = \infty$$

In this case, the random variable with null odd moments and the $(2n)$-moments $\sigma_n(I)$ is given by the semi-circle law (1.8).

Another direct consequence of theorem 1.2 is the Marchenko–Pastur law for non-isotropic Wishart matrices due to Y.Q. Yin [102, 103]; see also [18] and [80].

**Corollary 1.4 ([102, 103]).** Consider a collection $P : r \mapsto P(r)$ of possibly random matrices satisfying the condition stated in (1.16). Let $N = r/\rho$ be a scaling of the sample size in terms of the dimension associated with some parameter $\rho > 0$. For any $n \geq 1$ we have the Kreweras-type formula

$$\lim_{r\to\infty} r^{-1} \operatorname{Tr}\left(\mathbb{E}\left[P_N^n\right]\right) = \sum_{1\leq m\leq n} \rho^{n-m} \sum_{\mu\vdash[n] \;:\; m+|\mu|=n+1} K(\mu) \; \tau_\mu(P)$$

with the trace parameters $\tau_\mu(P)$ and the Kreweras numbers $K(\mu)$ defined in (1.18) and later in (2.13).

When $P = I$ the above limit result reduces to

$$\lim_{r\to\infty} r^{-1} \; \mathbb{E}\left[\operatorname{Tr}\left(P_N^n\right)\right] = \sum_{1\leq m\leq n} \rho^{n-m} \; N_{n,m}$$

In this situation, we also have the centered version limiting result

$$\lim_{r\to\infty} r^{-1}\mathbb{E}\left(\operatorname{Tr}\left([P_N - I]^n\right)\right) = \sum_{1\leq m\leq \lfloor n/2\rfloor} \rho^{n-m} \; R_{n,m}$$

In the above display, $N_{n,m}$ and $R_{n,m}$ denote the Narayana and the Riordan numbers defined later in (2.13) and (2.15). The matrix version of these isotropic results can be found in corollary 3.3 with the trace-type Marchenko–Pastur law (1.10) a simple corollary. See Section 3.3 for a new matrix version of a non-isotropic Marchenko–Pastur law.

Our third main result concerns moment estimates. We let $\|.\|_{op}$ and $\|.\|_F$ denote the operator norm and the Frobenius norm. In this notation, we have the following theorem.

**Theorem 1.5.** For some sufficiently small time horizon and for any sufficiently large sample size and any $n \geq 1$ we have the estimates

$$N \left\|\mathbb{E}\left[\mathcal{H}_N^{2n}\right] - \mathbb{E}\left[\mathcal{H}^{2n}\right]\right\|_F \leq c_1^n \ (2n)_n \ \mathrm{Tr}(P)^{2n}$$
$$\sqrt{N} \left\|\mathbb{E}\left(\exp\left(t\mathcal{H}_N\right)\right) - \mathbb{E}\left(\exp\left(t\mathcal{H}\right)\right)\right\|_F \leq c_2 \ (t\,\mathrm{Tr}(P))^3$$
$$\mathbb{E}\left[\|\mathcal{H}\|_{op}\right] \wedge \mathbb{E}\left[\|\mathcal{H}_N\|_{op}\right] \leq c_3 \sqrt{r} \ \lambda_1(P)$$

for some finite universal constants $c_1, c_2, c_3 < \infty$ whose values do not depend on the dimension parameter, nor on the parameter $n$.

In the above display, $\lambda_1(P) = \|P\|_{op}$ denotes the maximal eigenvalue of $P$ (cf. 2.2).

A more precise statement with a more detailed description of the constants is provided in Section 4 and Section 5; see for instance theorem 4.1, theorem 4.4, corollary 4.5, and theorem 5.4. The operator norm estimate stated in the above theorem extends the norm estimate for isotropic random vectors presented in [75] in the context of Gaussian random matrices. These norm-type bounds are based on non-commutative versions of Khintchine-type inequalities for Rademacher series presented in [54, 55]. More sophisticated approaches based on Burkholder/Rosenthal martingale-type inequalities are also developed in [36, 37]. Nevertheless these inequalities cannot be used to estimate random operator norms and the constants are often not explicit.

The last part of these lecture notes is concerned with non-asymptotic exponential concentration inequalities for traces and the operator norm of the fluctuation matrix. In this context, our main results can be stated as follows.

**Theorem 1.6.** For any symmetric matrix $A$, any $\delta \geq 0$, and any sufficiently large sample size the probabilities of the following events

$$|\mathrm{Tr}(A\mathcal{H}_N)| \leq c_1 \sqrt{(\delta+1) \ \left[\mathrm{Tr}((AP)^2) + \|AP\|_F^2\right]}$$
$$\|\mathcal{H}_N\|_{op} \leq \ c_2 \ \lambda_1(P) \ \sqrt{\delta + r}$$
$$\sup_{1 \leq k \leq r} |\lambda_k(P_N) - \lambda_k(P)| \leq \ c_2 \ \lambda_1(P) \ \sqrt{(\delta+r)/N}$$

are greater than $1-e^{-\delta}$, where $c_1, c_2$ denote some universal constants.

In the above display, $\lambda_i(P_N)$ and $\lambda_i(P)$ denote the ordered (decreasing in magnitude) eigenvalues of $P_N$, resp. $P$ (cf. 2.2). For a precise statement of this result and a detailed description of the constants $c_1$, $c_2$ we refer to Section 5; see in particular Section 5.2 and theorem 5.2, theorem 5.4 and corollary 5.5.

The proof of the trace-type concentration inequality is based on sub-Gaussian Laplace estimates of well known Wishart trace-type Laplace transforms (5.5). See Section 5.2 for a description of these sub-Gaussian estimates; e.g. see (5.7) and the first assertion in theorem 5.2. The operator norm concentration inequality comes from the variational formulation

$$\|\mathcal{H}_N\|_{op} = \sup_{x,y\in\mathbb{B}} \langle \mathcal{H}_N x, y\rangle = \sup_{A\in\mathbb{A}} \operatorname{Tr}(A\mathcal{H}_N) \tag{1.20}$$

where $\mathbb{B}$ is the unit ball in $\mathbb{R}^r$ equipped with the Euclidian distance and $\mathbb{A}$ is the set of matrices

$$\mathbb{A} := \{A = (xy' + y'x)/2 \ : \ x, y \in \mathbb{B}\} \tag{1.21}$$

The last spectral concentration estimate is a direct consequence of Weyl's inequality (2.5).

We end this section with some comparisons of the above concentration inequalities with existing results in random matrix theory. When $P = I$ the joint density of the random eigenvalues of $P_N$ is explicitly known; see e.g. [4]. Elegant Sanov-type large deviation principles for the spectral empirical measures have been developed by G. Ben Arous and A. Guionnet [5]. The literature also consists of non-asymptotic concentration inequalities for sums of independent random matrices. We refer to the seminal book of J. Tropp [93] for the state of the art on these topics. See also the review [94].

We also emphasize that the Laplace transform-type techniques developed in the present study differ from the ones based on Lieb's inequality (4.13). The latter are often used to control the largest eigenvalue of a random matrix using trace-type estimates; see proposition 4.4 and Section 4.5 in [93].

Other types of models have been considered in the literature leading to different results. For example, Gaussian concentration inequalities have been derived for Rademacher and Gaussian series associated with deterministic self-adjoint matrices; see e.g. theorem 2.1 in [93]. Matrix Hoeffding, Bernstein and Azuma-type inequalities have been derived for almost surely bounded random matrices; see theorem 2.8 and theorem 8.1 in [93]. The concentration results developed in the present notes provide more refined estimates, but of course they are restricted to random Wishart matrix models.

# 2

# Some preliminary results

## 2.1 Laplace transforms and exponential inequalities

We recall that for any non-negative random variable $Z$ such that

$$\mathbb{E}\left(Z^n\right)^{1/n} \leq z\ n \quad \left(\text{with } z = \sqrt{\tfrac{e}{2}}\, y\right)$$

$$\Longleftrightarrow \mathbb{E}\left(Z^n\right) \leq y^n\ n! \quad (\text{with } y = e\, z)$$

for some parameters $y, z \neq 0$ and any $n \geq 1$, the probability of the event

$$Z/z \leq \frac{e^2}{\sqrt{2}}\ \left[\frac{1}{2} + \left(\delta + \sqrt{\delta}\right)\right]. \tag{2.1}$$

is greater than $1 - e^{-\delta}$, for any $\delta \geq 0$.

To check this claim, notice that for any $n \geq 1$ we have

$$\mathbb{E}\left((Z/z)^n\right) \leq n^n \leq \frac{e}{\sqrt{2}}\ \left(\frac{e}{2}\right)^n \mathbb{E}(V^{2n})$$

for some Gaussian and centered random variable $V$ with unit variance. We check this claim using Stirling approximation

$$\begin{aligned}\mathbb{E}(V^{2n}) &= 2^{-n}\frac{(2n)!}{n!}\\ &\geq e^{-1}\ 2^{-n}\frac{\sqrt{4\pi n}\ (2n)^{2n}\ e^{-2n}}{\sqrt{2\pi n}\ n^n\ e^{-n}} = \sqrt{2}e^{-1}\ \left(\frac{2}{e}\right)^n\ \ n^n\end{aligned}$$

The estimate (2.1) is now a direct consequence of [22, Proposition 11.6.6].

Tail probability estimates can be used to obtain moments and Laplace transform estimate. The transfer of exponential tail estimates to these statistical models relies on the integral formulae

$$\mathbb{E}\left(Z^n\right) = n\int_0^\infty\ z^{n-1}\ \mathbb{P}(Z\geq z)\ dz$$

and

$$\mathbb{E}\left(\exp\left[tZ\right]\right) = 1+t\ \int_0^\infty\ \exp\left[tz\right]\ \mathbb{P}(Z\geq z)\ dz$$

which are valid for any $n\geq 1$ and $t\in\mathbb{R}$.

Laplace transform estimates can also be used to estimate moments using the formula

$$\mathbb{E}\left(Z^n\right)\ \leq\ \left(\frac{n}{et}\right)^n\ \mathbb{E}\left(\exp\left[tZ\right]\right)$$

which are valid for any $n\geq 1$ and $t\geq 0$. To check this claim we use the decomposition

$$\begin{aligned}\mathbb{E}\left(Z^n\right) = \mathbb{E}\left(Z^n\ e^{-tZ}\ e^{tZ}\right)\ &\leq\ \sup_{z\geq 0}\left[z^n\ e^{-tz}\right]\ \mathbb{E}\left(\exp\left[tZ\right]\right)\\ &= \exp\left[-\inf_{z\geq 0}\left(tz-n\log z\right)\right]\ \mathbb{E}\left(\exp\left[tZ\right]\right)\end{aligned}$$

and we note that the infimum is attained at $z=n/t$.

We end this section with some basic Cramér–Chernov tools to derive quantitative concentration inequalities. We associate with any non-negative convex function $L$ defined on some domain $\mathrm{Dom(L)}\subset\mathbb{R}_+$, $L(0)=0$, the Legendre–Fenchel transform

$$L^\star(\lambda) := \sup_{t\in\mathrm{Dom(L)}}\ \left(\lambda t - L(t)\right)$$

defined for any $\lambda \geq 0$. Note that $L^{\star}$ is a convex increasing function with $L^{\star}(0) = 0$ and its inverse $(L^{\star})^{-1}$ is a concave increasing function. We let $L_Z(t) := \log \mathbb{E}(\exp(tZ))$ be the log-Laplace transform of a random variable $Z$ defined on some domain $\mathrm{Dom}(L_Z) \subset \mathbb{R}_+$.

Using the Cramér–Chernov–Chebychev inequality, for any $\lambda \geq 0$ and any $\delta \geq 0$ we find that

$$\log \mathbb{P}(Z \geq \lambda) \leq -L_Z^{\star}(\lambda)$$

Equivalently, we have

$$\mathbb{P}\left(Z \geq (L_Z^{\star})^{-1}(\delta)\right) \leq e^{-\delta}$$

## 2.2 Matrix norms and spectral analysis

Firstly, we note some general matrix notation. We let $\mathcal{M}_r$ be the space of $(r \times r)$-matrices with real entries. We denote by $\mathcal{S}_r \subset \mathcal{M}_r$ the subspace of symmetric $(r \times r)$-matrices, $\mathcal{S}_r^0 \subset \mathcal{S}_r$ the subset of semi-definite matrices, and $\mathcal{S}_r^+ \subset \mathcal{S}_r^0$ the subset of positive semi-definite matrices. We let $\lambda_i(A) \in \mathbb{C}$, with $1 \leq i \leq r$, denote the eigenvalues of $A \in \mathcal{M}_r$. When $A \in \mathcal{S}_r$ we assume that

$$\lambda_1(A) \geq \lambda_2(A) \geq \ldots \geq \lambda_r(A) \tag{2.2}$$

For any collection $A = (A_r)_{r\geq 1}$ of $(r \times r)$-matrices we set

$$\lambda^{\star}(A) = \max_{1\leq i\leq r} |\lambda_i(A(r))| \quad \text{and} \quad \lambda_{\star}(A) = \sup_{r\geq 1} \lambda^{\star}(A(r)) \in [0, \infty]$$

We have already defined $\|.\|_{op}$ and $\|.\|_F$ as the operator norm and the Frobenius norm. The $n$-Schatten norm of a matrix $Q \in \mathcal{S}_r$, with $n \geq 1$, is defined by

$$\|Q\|_n := \left[\sum_{1\leq k\leq r} |\lambda_k(Q)|^n\right]^{1/n}$$

This implies that

$$\|Q\|_{2n}^{2n} = \mathrm{Tr}(Q^{2n}) = \|Q^n\|_F^2$$

The 2-Schatten norm reduces to the Frobenius norm.

For any $P, Q \in \mathcal{S}_r^0$ and any $m, n \geq 0$ we have the estimate

$$0 \leq \mathrm{Tr}((PQ)^{n+m}) \leq \mathrm{Tr}((PQ)^m)\ \mathrm{Tr}((PQ)^n) \tag{2.3}$$

We check this claim using the fact that

$$\begin{aligned} 0 \leq \mathrm{Tr}((Q^{1/2}PQ^{1/2})^{n+m}) &= \mathrm{Tr}((PQ)^{n+m}) \\ &\leq \mathrm{Tr}((Q^{1/2}PQ^{1/2})^m)\ \mathrm{Tr}((Q^{1/2}PQ^{1/2})^n) \\ &= \mathrm{Tr}((PQ)^m)\ \mathrm{Tr}((PQ)^n) \end{aligned}$$

For any $(r \times r)$-symmetric matrices $P$ and $H$ we have the $n$-Wielandt–Hoffman inequality

$$\sum_{1 \leq k \leq r} |\lambda_k(P+H) - \lambda_k(P)|^n \leq \|H\|_n^n \tag{2.4}$$

The case $n = 2$ is sometimes called Mirsky's inequality. In the above display, $\lambda_k(Q)$ stands for the eigenvalues of a given symmetric matrix $Q \in \mathcal{S}_r$ ranked in decreasing order. We also have Weyl's inequality

$$\sup_{1 \leq k \leq r} |\lambda_k(P+H) - \lambda_k(P)| \leq \|H\|_{op} \tag{2.5}$$

For a further discussion on the perturbation theory of eigenvalues we refer to [96] and [61].

## 2.3 Tensor products and Fréchet derivatives

The tensor products $(A \otimes B)$ and $(A \overline{\otimes} B)$ and the entry-wise tensor product $(A \otimes B)^\sharp$ of matrices $A$ and $B$ with appropriate dimensions are defined by the formulae

$$(A \otimes B)_{(i,j),(k,l)} = (A \overline{\otimes} B)_{(i,j),(l,k)} = (A \otimes B)^\sharp_{(i,k),(j,l)} = A_{i,k}\ B_{j,l} \tag{2.6}$$

We also consider the symmetric tensor product

$$4(A \widehat{\otimes} B) = (A \otimes B) + (B \otimes A) + (A \overline{\otimes} B) + (B \overline{\otimes} A) \tag{2.7}$$

Let $\mathcal{B}_1, \mathcal{B}_2$ be a couple of Banach spaces equipped with some norm $\|.\|_{\mathcal{B}_i}$, with $i = 1, 2$; also let $\mathcal{O}_1 \subset \mathcal{B}_1$ be an open subset of $\mathcal{B}_1$. We let $\mathcal{L}(\mathcal{B}_1, \mathcal{B}_2)$ be the set of bounded linear functional from $\mathcal{B}_1$ into $\mathcal{B}_2$.

We recall that a mapping $\Upsilon : \mathcal{B}_1 \mapsto \mathcal{B}_2$ is Fréchet differentiable at some $A \in \mathcal{O}_1$ if there exists a continuous linear functional

$$\nabla\Upsilon(A) \in \mathcal{L}(\mathcal{B}_1, \mathcal{B}_2)$$

such that

$$\lim_{\|H\|_{\mathcal{B}_1}\to 0} \|H\|_{\mathcal{B}_1}^{-1} \|\Upsilon(A+H) - \Upsilon(A) - \nabla\Upsilon(A)\cdot H\|_{\mathcal{B}_2} = 0 \tag{2.8}$$

The mapping is said to be twice Fréchet differentiable at $A \in \mathcal{O}_1$ when $\Upsilon$ and the mapping

$$\nabla\Upsilon \;:\; A \in \mathcal{O}_1 \mapsto \nabla\Upsilon(A) \in \mathcal{L}(\mathcal{B}_1, \mathcal{B}_2)$$

is also Fréchet differentiable, and so on. Given some Fréchet differentiable mapping $\Upsilon$ of the third order, at some $A \in \mathcal{O}_1$ for any $H \in \mathcal{B}_1$ we have

$$\Upsilon(A+H) = \Upsilon(A) + \nabla\Upsilon(A)\cdot H + \frac{1}{2}\,\nabla^2\Upsilon(A)\cdot(H,H) + \overline{\nabla}^3\Upsilon\,[A,H]$$

with the remainder functional in the Taylor expansion given

$$\overline{\nabla}^3\Upsilon\,[A,H] := \frac{1}{3!}\int_0^1 (1-\epsilon)^2\;\nabla^3\Upsilon\,(A+\epsilon\,H)\cdot(H,H,H)\;d\epsilon$$

We also consider the multi-linear operator norm

$$|||\nabla^n\Upsilon(P)||| = \sup_{\|H\|_{\mathcal{B}_1}=1} \|\nabla^n\Upsilon(P)\cdot H^{\otimes n}\|_{\mathcal{B}_2}$$

We says $\Upsilon$ satisfies the third-order polynomial growth condition when we have

$$|||\nabla^n\Upsilon(P)||| \le \kappa_n(A)$$

for any $n = 1, 2,$ and

$$\|\overline{\nabla}^3\Upsilon\,[A,H]\,\| \le \kappa_3(A)\;\left(1+\|H\|_{\mathcal{B}_1}^{\alpha}\right)\;\|H\|_{\mathcal{B}_1}^3 \tag{2.9}$$

for some parameter $\alpha > 0$ and some constants $\kappa_n(A)$ whose values only depends on $A$. For instance, the principal square root functional

$$\psi \;:\; Q \in \mathcal{S}_r^+ \mapsto \psi(Q) = Q^{1/2} \in \mathcal{S}_r^+$$

is Fréchet differentiable at any order and satisfies the third-order polynomial growth condition stated above. More precisely, as shown in [11] the first and second order derivative given for any $(A,H)\in(\mathcal{S}_r^+\times\mathcal{S}_r)$ s.t. $A+\epsilon\ H\in\mathcal{S}_r^0$ for any $\epsilon\in[0,1]$ by

$$\begin{aligned}\nabla\psi(A)\cdot H &= \int_0^\infty e^{-t\psi(A)}\ H\ e^{-t\psi(A)}\ dt\\ 2^{-1}\nabla^2\psi(A)\cdot(H,H) &:= -\nabla\psi(A)\cdot(\nabla\psi(A)\cdot H)^2\end{aligned}\tag{2.10}$$

In addition, we have the second order approximation

$$\begin{aligned}&\|\psi(A+H)-\psi(A)-\nabla\psi(A)\cdot H-2^{-1}\nabla^2\psi(A)\cdot(H,H)\|\\ &\leq\frac{c^2}{4}\ \lambda_{min}(A)^{-5/2}\ \|H\|^3\end{aligned}\tag{2.11}$$

where $c=\sqrt{r}$ for the Frobenius norm $\|.\|=\|.\|_F$, and $c=1$ for the $\mathbb{L}_2$-norm $\|.\|=\|.\|_{\mathbb{L}_2}$. The first and second derivatives can be estimated by the operator norm inequalities:

$$|||\nabla\psi(A)|||\leq\ 2^{-1}\ \lambda_{min}(A)^{-1/2}\quad\text{and}\quad\left|\left|\left|\nabla^2\psi(A)\right|\right|\right|\leq\frac{K}{4}\ \lambda_{min}(A)^{-3/2}$$

Let $e_k$ be the $r$-column vectors with null entries but the $k$-th unit one, and let $e_{i,j}=e_ie_j'$ the $(r\times r)$-matrix with null entries but the $(i,j)$-th unit one.

In this notation, we can identify $\nabla\Upsilon(Q)$ with the tensor

$$\nabla\Upsilon(Q)_{(i,j),(k,l)}:=(\nabla\Upsilon(Q)\cdot e_{i,j})_{k,l}=e_k\ \ (\nabla\Upsilon(Q)\cdot e_{i,j})\ \ e_l'$$

In this notation, we have

$$H=\sum_{i,j}\ H_{i,j}\ e_{i,j}$$

$$\Longrightarrow(\nabla\Upsilon(Q)\cdot H)_{(k,l)}=\sum_{i,j}\ H_{(i,j)}\ \nabla\Upsilon(Q)_{(i,j),(k,l)}=(H\ \nabla\Upsilon(Q))_{(k,l)}$$

In the same vein, we can identify $\nabla^2\Upsilon(Q)$ with the tensor

$$\nabla^2\Upsilon(Q)_{((i,j),(k,l)),(m,n)}:=\left(\nabla^2\Upsilon(Q)\cdot(e_{i,j},e_{k,l})\right)_{m,n}\tag{2.12}$$

This implies that

$$\left(\nabla^2\Upsilon(Q)\cdot(H_1,H_2)\right)=(H_1\otimes H_2)^\sharp\ \nabla^2\Upsilon(Q)$$

## 2.4 Partitions and Bell polynomials

Let $\mathcal{P}_{n,m}(\mu) \subset \mathcal{P}_{n,m}$ be the subset of partitions of type $\mu \vdash n$, with $|\mu| = \mu_1 + \ldots + \mu_n = m$ blocks. The cardinality of these set partitions are given by the Bell numbers $B(n)$ and the Stirling numbers $S(n,m)$ of the second kind; that is, we have that

$$B(n) := |\mathcal{P}_n| = \sum_{1 \leq m \leq n} S(n,m) \quad \text{with} \quad S(n,m) := |\mathcal{P}_{n,m}|$$

Let $\mathcal{N}_{n,m} \subset \mathcal{P}_{n,m}$ the subset of non-crossing partitions, and $\mathcal{N}_{n,m}(\mu) \subset \mathcal{P}_{n,m}(\mu)$ the subset of non-crossing partitions of type $\mu$. The number of non-crossing partitions of $n$ are computed using the Narayana and the Kreweras numbers

$$\begin{aligned} N_{n,m} := |\mathcal{N}_{n,m}| &= \frac{1}{n} \binom{n}{m-1} \binom{n}{m} \\ &= \frac{1}{m} \binom{n-1}{m-1} \binom{n}{m-1} = |\mathcal{N}_{n,n+1-m}| \end{aligned} \tag{2.13}$$

and

$$K(\mu) = |\mathcal{N}_{n,m}(\mu)| = \frac{1}{(n+1-m)!} \begin{pmatrix} & n & \\ \mu_1 \; \mu_2 & \ldots & \mu_n \end{pmatrix}$$

for any $\mu \vdash n$, with $|\mu| = m \leq n$. A proof of the above formula can be found in [65, e.g. corollary 9.12], see also [24, 51], the pioneering article by G. Kreweras [44] and the survey of R. Simion on non-crossing partitions [81].

Let $\mathcal{Q}^+_{n,m} \subset \mathcal{Q}_{n,m}$ be the subset of non-crossing partitions, and $\mathcal{Q}^-_{n,m} \subset \mathcal{Q}_{n,m}$ the subset of crossing partitions without singleton.

Notice that $\mathcal{Q}^+_{n,m} = 0$ as soon as $m > \lfloor n/2 \rfloor$. In addition we have

$$\pi \in \mathcal{Q}_{n,m} \quad \text{and} \quad m > \lfloor n/2 \rfloor \implies \exists i \in [m] \; : \; |\pi_i| = 1 \tag{2.14}$$

The number of non-crossing partitions of $[n]$ without singletons is given by the Riordan numbers

$$R_n = \sum_{1 \leq m \leq \lfloor n/2 \lfloor} R_{n,m} \tag{2.15}$$

with

$$R_{n,m} = |\mathcal{Q}^+_{n,m}| = \frac{1}{n+1} \left( \begin{array}{c} n+1 \\ m \end{array} \right) \left( \begin{array}{c} n-1-m \\ m-1 \end{array} \right)$$

A proof of the above assertion can be found in [42].

For a centered Gaussian random variable $V$ with unit variance we have the estimates

$$\begin{aligned} R_{2n,n} &= N_{2n,n}(0,n,0,\ldots,0) = C_n \\ &\leq \ |\mathcal{Q}_{2n,n}| = 2^{-n}\ (2n)_n = \mathbb{E}(V^{2n}) := v_{2n} \end{aligned}$$

with the Catalan numbers $C_n$ defined in (1.9).

We also recall that Stirling numbers of the first kind are denoted by $s(n,m)$ and the Pochhammer symbol is defined by $(n)_m := n!/(n-m)!$. We also denote by $r_i(\pi)$ the number of blocks of size $i \geq 1$ in the partition $\pi$.

For a detailed discussion, and the combinatorial interpretation of non-crossing partitions, Catalan, Riordan and Narayana numbers, we refer to the series of articles [6, 12, 74, 82, 83, 84, 85, 101]. For instance, Narayana numbers $N_{n,m}$ count the number of expressions with $n$ pairs of matched parenthesis and $m$ different nestings. These numbers also count the number of ordered trees with $n$ edges and $m$ leaves; as well as the number of Dicks paths of length $n$ and $m$ distinct picks. Recall that Dicks paths are staircase walks with $2n$ steps in the plane starting at the origin $(0,0)$ and ending in $(n,n)$, and staying above the diagonal (they may touch the diagonal). The Riordan number $R_n$ coincides with the number of plane trees with $n$ edges, in which no vertex has outdegree one. The $n$-th Riordan number also counts Motzkin paths of length $n$ with no horizontal steps of null height. Recall that Motzkin paths are lattice paths in the plane starting from the origin up to $(n,0)$, using up steps $(1,1)$, down steps $(1,-1)$, and horizontal steps $(1,0)$, and never going below the $x$-axis. Identifying $k \in [n]$ with $e^{-2ik\pi/n}$, then non-crossing partitions can also be seen as compact subsets of the unit disk [43].

We also consider the complete Bell polynomials on a commutative algebra $\mathcal{A}$ given by $B_0 = I$ the unit neutral state, and for any $n \geq 1$ and any $x = (x_1, \ldots, x_n) \in \mathcal{A}^n$ by the formula

$$\begin{aligned} B_n(x) &= \sum_{1 \leq k \leq n} B_{k,n}(x) \\ &= n! \sum_{\sum_{1\leq i\leq n} il_i = n} \prod_{1 \leq i \leq n} \left[ \frac{1}{l_i!} \left( \frac{x_i}{i!} \right)^{l_i} \right] \end{aligned} \tag{2.16}$$

with the partial Bell polynomials

$$B_{k,n}(x) = \frac{n!}{k!} \sum_{m_1+\ldots+m_k=n,\, m_i \geq 1} \frac{x_{m_1} \ldots x_{m_k}}{m_1! \ldots m_k!} \tag{2.17}$$

We recall that

$$\exp\left[ \sum_{n\geq 1} \frac{t^n}{n!} x_n \right] = \sum_{n \geq 0} \frac{t^n}{n!} y_n \quad \text{with} \quad y_n = B_n(x_1, \ldots, x_n)$$

When $y_0 = I$ we also have the inverse relations

$$x_n = \sum_{1 \leq k \leq n} (-1)^{k-1} (k-1)!\, B_{n,k}(y_1, \ldots, y_n) \tag{2.18}$$

We have the following technical lemma whose proof is given in the Appendix.

**Lemma 2.1.** Assume that $\forall\ 1 \leq k < p$ we have $x_k = 0$ and $\forall n \geq p$ we have $|x_n| \leq n!\, \rho\, \beta^n$ for some $p \geq 1$ and non-negative parameters $\rho, \beta \geq 0$. Then, for any $1 \leq q < p$ we have

$$\left| B_{pn}(x) - \frac{(np)!}{n!} \left( \frac{x_p}{p!} \right)^n \right| \leq \frac{(np)!}{2} (2\beta)^{pn} \sum_{1 \leq k < n} \frac{1}{k!} \left( \frac{\rho}{2^{p-1}} \right)^k$$

as well as

$$|\, B_{pn+q}(x)| \leq \frac{(np+q)!}{2} (2\beta)^{pn+q} \sum_{1 \leq k \leq n} \frac{1}{k!} \left( \frac{\rho}{2^{p-1}} \right)^k$$

In addition, when $x_n \geq 0$ for any $n \geq 1$ we have

$$B_n(x) \geq 0 \quad \text{and} \quad B_{pn}(x) \geq \frac{(np)!}{n!} \left( \frac{x_p}{p!} \right)^n$$

## 2.5 Some combinatorial constructions

Given some partitions $\pi \in \mathcal{P}_{n,m}$ and $\pi' \in \mathcal{P}_{n',m'}$ we let $\pi \oplus \pi' \in \mathcal{P}_{n+n',m+m'}$ the partition of $[n+n']$ with $m+m'$ blocks given by

$$\forall i \in [m+m'] \qquad (\pi \oplus \pi')_i = \begin{cases} \pi_i & \text{if} \quad i \in [n] \\ n + \pi'_i & \text{if} \quad i \in n + [n'] \end{cases}$$

We also consider the closure operator

$$\mathrm{cl} \; : \; \pi \in \mathcal{N}_{n,m} \mapsto \mathrm{cl}\,(\pi) \in \mathcal{N}_{n+1,m}$$

defined by the blocks

$$\mathrm{cl}\,(\pi)_{i+1} = \begin{cases} 1 + \pi_i & \text{if} \quad n \not\in \pi_i \quad \text{and} \quad 1 \leq i \leq m \\ \{1\} \cup (1 + \pi_j) & \text{if} \quad n \in \pi_j \quad \text{and} \quad i = 0 \end{cases}$$

For example, the partition $\pi = (\pi_1, \pi_2, \pi_3, \pi_4) \in \mathcal{N}_{7,4}$ with blocks,

$$\pi_1 := \{1,4\} \leq \pi_2 := \{2,3\} \leq \pi_3 := \{5\} \leq \pi_4 := \{6,7\} \tag{2.19}$$

is defined by the Murasaki diagram presented in Figure 2.1. We recall that $\iota(\pi)$ is the number of blocks visible from above in the Murasaki diagram associated with $\pi$

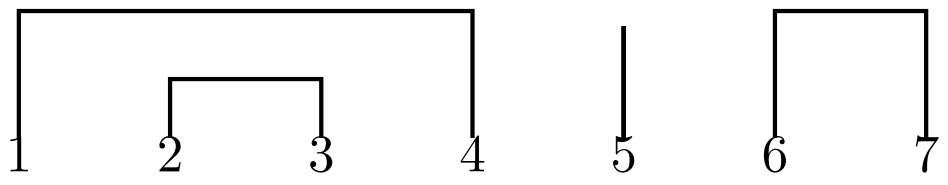

**Figure 2.1:** Murasaki diagram representation of $\pi$

Observe that the partition $\pi$ in (2.19) has $\iota(\pi) = 3$ blocks visible from above. Its closure

$$7 \in \pi_4 \qquad \Longrightarrow \qquad \mathrm{cl}\,(\pi) = (\{1,7,8\}, \{2,5\}, \{3,4\}, \{6\}) \in \mathcal{N}_{8,4}$$

is defined by the Murasaki diagram presented in Figure 2.2.

Writing clockwise the nodes the partition $\pi$ defined in (2.19) also has the circle representation given in Figure 2.3.

Given a circular representation of a partition $\pi \in \mathcal{N}_{n,m}$, we subdivide each of the $n$ arcs by a new blue node. The blue node which subdivides

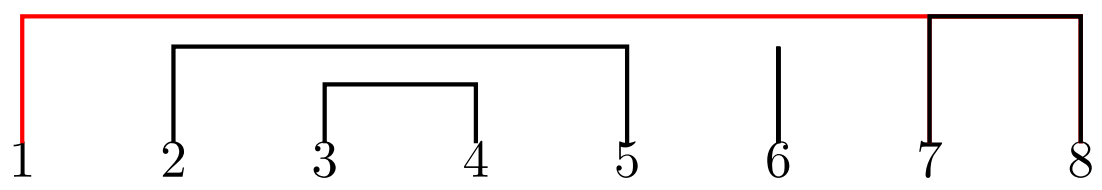

$$\Longrightarrow \qquad \iota\left(\mathrm{cl}\left(\pi\right)\right) = 1.$$

**Figure 2.2:** Murasaki diagram representation of $\pi$ and $\mathrm{cl}\left(\pi\right)$

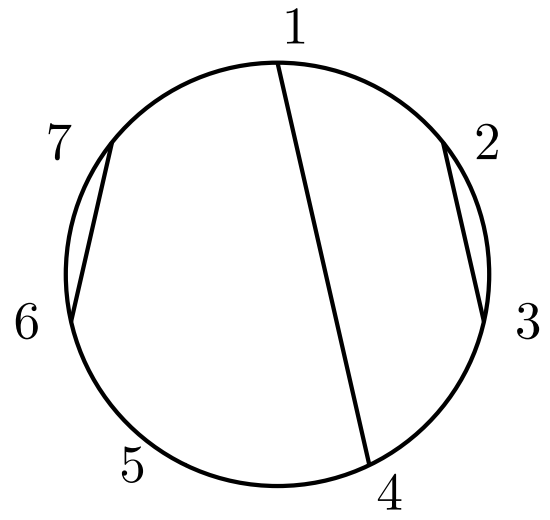

**Figure 2.3:** Circular representation of $\pi$

the arc $\widehat{(1,2)}$ is labelled with 1. The remaining subdivision of nodes $2, \ldots, n$ are placed clockwise. We also consider the mapping

$$\Xi \ : \ \pi \in \mathcal{N}_{n,m} \mapsto \Xi\left(\pi\right) \in \mathcal{N}_{n,n+1-m} \tag{2.20}$$

defined by taking the coarsest non-crossing partition $\Xi\left(\pi\right)$ of the blue nodes whose chords don't cross the chord representation of $\pi$. For instance the blue nodes and the circular representations of $\Xi\left(\pi\right)$ associated with the partition $\pi$ defined in Figure 2.3 are given in Figure 2.4. Up to a change of numbering $\Xi$ coincides with the order-reversing involution introduced in [83].

Observe that for any partition $\pi$ of $[n]$ the block of $\Xi\left(\pi\right)$ containing the node with label $n$ is a block of size $\iota(\pi)$. In other words, the size of the block of $\Xi\left(\pi\right)$ containing the node with label $n$ characterizes the number of visible blocks of $\pi$ from above.

Next we analyze the difference between $\Xi\left(\pi \oplus \pi'\right)$ and $\Xi\left(\pi\right) \oplus \Xi\left(\pi'\right)$. We start with a simple example. Consider the partitions with Murasaki diagrams in Figure 2.5 and Figure 2.6.

The Murasaki diagram of $\pi \oplus \pi' \in \mathcal{N}_{23,12}$ is given as in Figure 2.7.

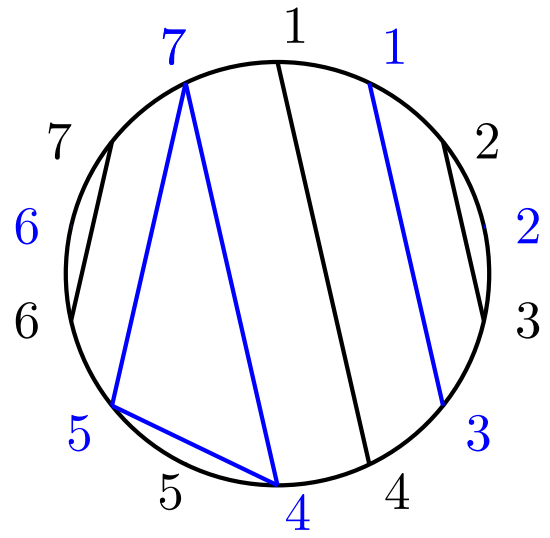

**Figure 2.4:** Circular representations of $\pi$ and $\Xi(\pi)$

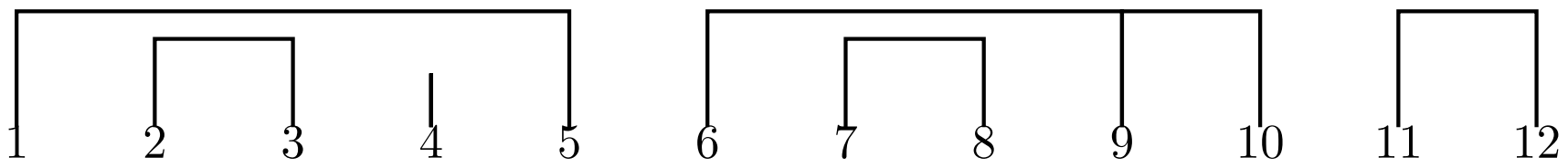

**Figure 2.5:** Murasaki diagram of $\pi \in \mathcal{N}_{12,6}$, with $\iota(\pi) = 3$

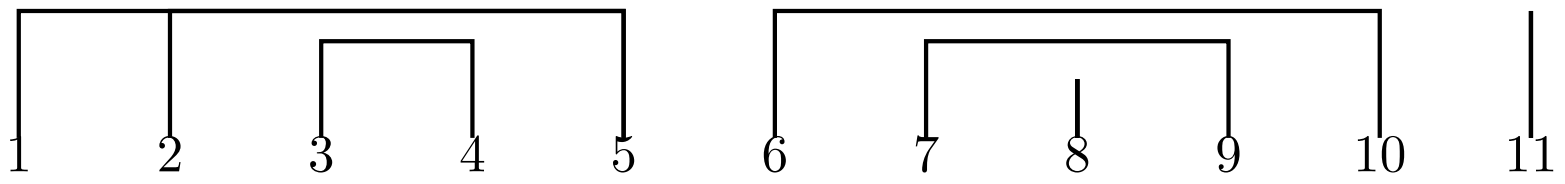

**Figure 2.6:** Murasaki diagram of $\pi' \in \mathcal{N}_{11,6}$, with $\iota(\pi') = 3$

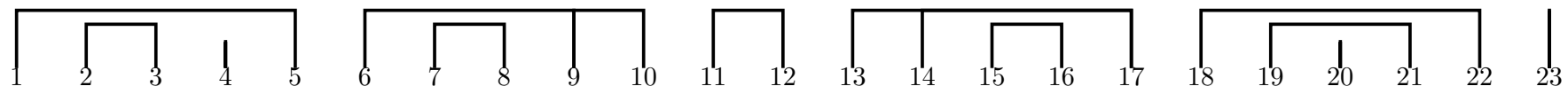

**Figure 2.7:** Murasaki diagram of $\pi \oplus \pi' \in \mathcal{N}_{23,12}$

The Murasaki diagram of $\Xi(\pi) \in \mathcal{N}_{12,6}$ and $\Xi(\pi') \in \mathcal{N}_{11,6}$ are defined in Figure 2.8 and Figure 2.9. Observe that the size of the blue block of $\Xi(\pi)$ containing the node "12" and the one of $\Xi(\pi')$ containing the node "11" coincides with $\iota(\pi)$ and respectively $\iota(\pi')$.

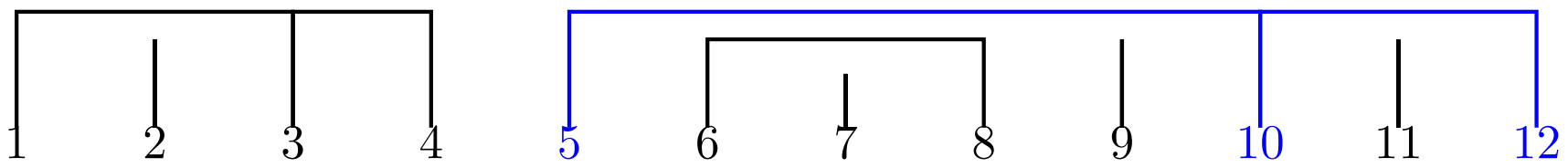

**Figure 2.8:** Murasaki diagram of $\Xi(\pi) \in \mathcal{N}_{12,6}$

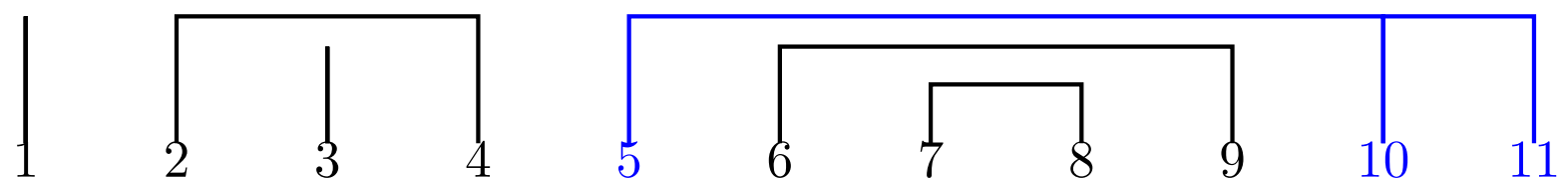

**Figure 2.9:** Murasaki diagram of $\Xi(\pi') \in \mathcal{N}_{11,6}$

The Murasaki diagram of $\Xi(\pi) \oplus \Xi(\pi') \in \mathcal{N}_{23,12}$ is given in Figure 2.10.

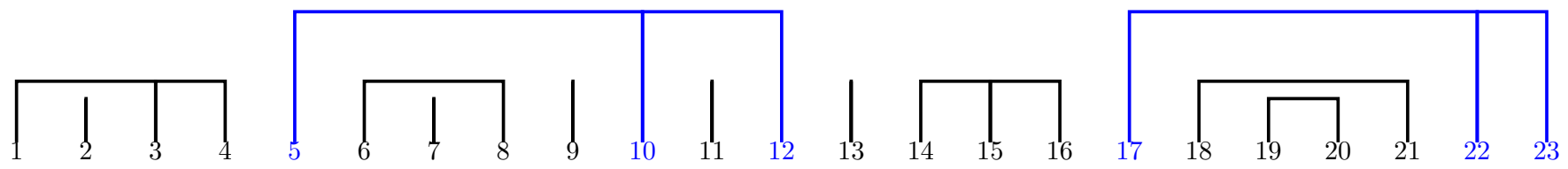

**Figure 2.10:** Murasaki diagram of $\Xi(\pi) \oplus \Xi(\pi') \in \mathcal{N}_{23,12}$

In terms of Murasaki diagrams, $\Xi(\pi \oplus \pi')$ is deduced from $\Xi(\pi) \oplus \Xi(\pi')$ by merging the blue blocks of sizes $\iota(\pi)$ and $\iota(\pi')$. This shows that $\Xi(\pi \oplus \pi')$ and $\Xi(\pi) \oplus \Xi(\pi')$ have the same combinatorial structure, but $\Xi(\pi \oplus \pi')$ looses a block of size $\iota(\pi)$ and a block of size $\iota(\pi')$, and gains a new block of size $\iota(\pi) + \iota(\pi')$ (so that the total number of blocks is kept unchanged). See Figure 2.11 (and compare with Figure 2.7 and Figure 2.10).

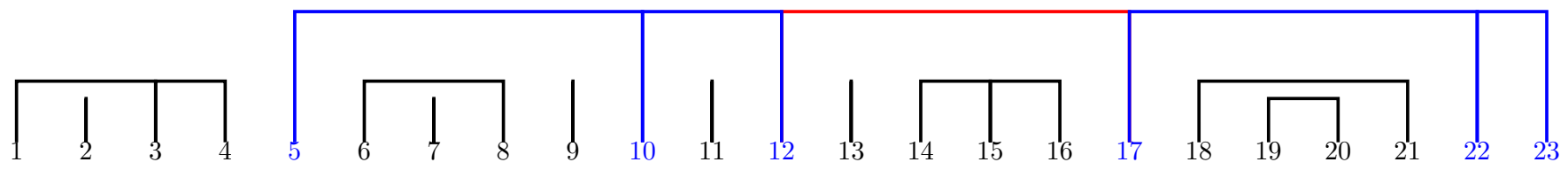

**Figure 2.11:** Murasaki diagram of $\Xi(\pi) \oplus \Xi(\pi') \in \mathcal{N}_{23,12}$ and the one of $\Xi(\pi \oplus \pi') \in \mathcal{N}_{23,12}$

Non-crossing partitions can be generated sequentially using for any $0 \le m \le n$ the formula

$$\mathcal{N}_{n+1,m+1} = [\mathcal{N}_{n,m} \oplus \mathrm{cl}(\mathcal{N}_{0,0})] \cup \bigcup [\mathcal{N}_{n_1,m_1} \oplus \mathrm{cl}(\mathcal{N}_{n_2,m_2})] \tag{2.21}$$

In the above display, the second set-union is taken over all non-negative parameters

$$0 \le m_1 \le n_1 \quad \text{and} \quad 0 \le m_2 \le n_2$$

such that

$$n_1 + n_2 = n \quad \text{and} \quad m_1 + m_2 = m + 1 \tag{2.22}$$

We also use the conventions $\mathcal{N}_{0,0} = \{\pi_0\}$, where $\pi_0 = \{0\}$ stands for the null partition on $[0]$, as well as $\{0\} \oplus \pi = \pi$. In addition for any $n \geq 1$ we set $\mathcal{N}_{n,0} = \emptyset$ so that $\emptyset \oplus \pi = \emptyset$.

For instance, the sets of non-crossing partitions with one or two elements are given by the formulae

$$\mathcal{N}_{1,1} = \mathcal{N}_{0,0} \oplus \mathrm{cl}\,(\mathcal{N}_{0,0}) = \mathrm{cl}\,(\mathcal{N}_{0,0}) = \{\{1\}\} =$$ 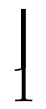

$$\mathcal{N}_{2,1} = \mathcal{N}_{0,0} \oplus \mathrm{cl}\,(\mathcal{N}_{1,1}) = \mathrm{cl}\,(\mathcal{N}_{1,1}) = \{\{1,2\}\} =$$ 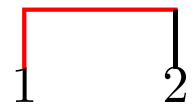

$$\mathcal{N}_{2,2} = \mathcal{N}_{1,1} \oplus \mathrm{cl}\,(\mathcal{N}_{0,0}) = \mathcal{N}_{1,1} \oplus \mathcal{N}_{1,1} = \{(\{1\},\{2\})\} =$$ 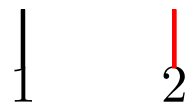

In the same vein, the sets of non-crossing partitions on $[3]$ are given by the decompositions

$$\mathcal{N}_{3,1} = \mathcal{N}_{0,0} \oplus \mathrm{cl}\,(\mathcal{N}_{2,1}) = \mathrm{cl}\,(\mathcal{N}_{2,1}) = \{\{1,2,3\}\} =$$ 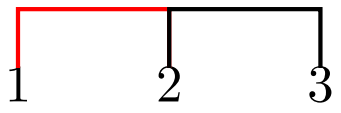

$$\mathcal{N}_{3,2} = [\mathcal{N}_{2,1} \oplus \mathrm{cl}\,(\mathcal{N}_{0,0})] \cup [\mathcal{N}_{1,1} \oplus \mathrm{cl}\,(\mathcal{N}_{1,1})] \cup [\mathcal{N}_{0,0} \oplus \mathrm{cl}\,(\mathcal{N}_{2,2})]$$

$=$ 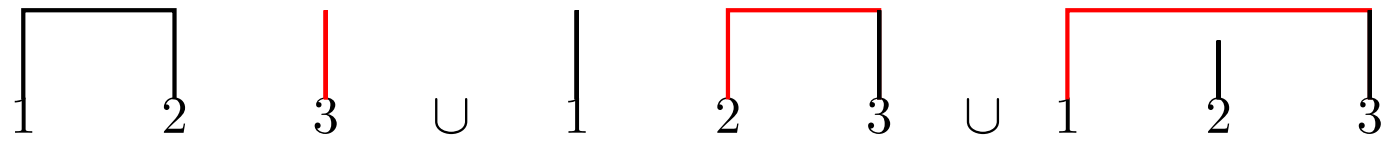

$$\mathcal{N}_{3,3} = \mathcal{N}_{2,2} \oplus \mathrm{cl}(\mathcal{N}_{0,0}) =$$ 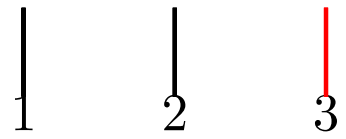

We also have the decompositions

$$\begin{aligned}
\mathcal{N}_{4,1} &= N_{0,0} \oplus \mathrm{cl}(\mathcal{N}_{3,1}) \\
\mathcal{N}_{4,2} &= [N_{3,1} \oplus \mathrm{cl}(\mathcal{N}_{0,0})] \cup [N_{2,1} \oplus \mathrm{cl}(\mathcal{N}_{1,1})] \\
&\qquad \cup [N_{1,1} \oplus \mathrm{cl}(\mathcal{N}_{2,1})] \cup [N_{0,0} \oplus \mathrm{cl}(\mathcal{N}_{3,2})] \\
\mathcal{N}_{4,3} &= [N_{3,2} \oplus \mathrm{cl}(\mathcal{N}_{0,0})] \cup [N_{2,2} \oplus \mathrm{cl}(\mathcal{N}_{1,1})] \\
&\qquad \cup [N_{1,1} \oplus \mathrm{cl}(\mathcal{N}_{2,2})] \cup [N_{0,0} \oplus \mathrm{cl}(\mathcal{N}_{3,3})] \\
\mathcal{N}_{4,4} &= \mathcal{N}_{3,3} \oplus \mathrm{cl}(\mathcal{N}_{0,0})
\end{aligned}$$

as well as the partitioning of non-crossing partitions on the set $[3]=\{1,2,3\}$ given by

$$
\begin{aligned}
\mathcal{N}_{5,1} &= N_{0,0}\oplus \mathrm{cl}(\mathcal{N}_{4,1})\\
\mathcal{N}_{5,2} &= [N_{4,1}\oplus \mathrm{cl}(\mathcal{N}_{0,0})]\cup[N_{3,1}\oplus \mathrm{cl}(\mathcal{N}_{1,1})]\cup[N_{2,1}\oplus \mathrm{cl}(\mathcal{N}_{2,1})]\\
&\qquad \cup[N_{1,1}\oplus \mathrm{cl}(\mathcal{N}_{3,1})]\cup[N_{0,0}\oplus \mathrm{cl}(\mathcal{N}_{4,2})]\\
\mathcal{N}_{5,3} &= [N_{2,2}\oplus \mathrm{cl}(\mathcal{N}_{2,1})]\cup[N_{2,1}\oplus \mathrm{cl}(\mathcal{N}_{2,2})]\\
&\qquad \cup[N_{1,1}\oplus \mathrm{cl}(\mathcal{N}_{3,2})]\cup[N_{0,0}\oplus \mathrm{cl}(\mathcal{N}_{4,3})]\\
\mathcal{N}_{5,4} &= [N_{3,3}\oplus \mathrm{cl}(\mathcal{N}_{1,1})]\cup[N_{2,2}\oplus \mathrm{cl}(\mathcal{N}_{2,2})]\\
&\qquad \cup[N_{1,1}\oplus \mathrm{cl}(\mathcal{N}_{3,3})]\cup[N_{0,0}\oplus \mathrm{cl}(\mathcal{N}_{4,4})]\\
\mathcal{N}_{5,5} &= N_{4,4}\oplus \mathrm{cl}(\mathcal{N}_{0,0})
\end{aligned}
$$

## 2.6 Fluctuation and bias approximations

Let $\mathcal{B}$ be some Banach space equipped with some norm $\|.\|_{\mathcal{B}}$.

**Theorem 2.2.** For any Fréchet differentiable mapping $\Upsilon:\mathcal{S}_r\mapsto\mathcal{B}$ at the third order, and satisfying the third-order polynomial growth condition (2.9), we have

$$\sqrt{N}\;[\Upsilon(P_N)-\Upsilon(P)]\;\longrightarrow_{N\to\infty}\;\nabla\Upsilon(P)\cdot\mathcal{H}$$

In addition, we have the bias estimate

$$\|\mathbb{E}\left(\Upsilon(P_N)\right)-\Upsilon(P)-N^{-1}\;(P\,\widehat{\otimes}\,P)\,\nabla^2\Upsilon(P)\|_{\mathcal{B}}\;\leq\;c\;\kappa_3(P)\;N^{-3/2}$$

with the tensor composition (2.12), and the parameter $\kappa_3(A)$ introduced in (2.9).

*Proof.* The first assertion is a direct consequence of the continuous mapping theorem. In terms of the tensor composition (2.12), we also have the bias estimate

$$N^{3/2}\;\|\mathbb{E}\left(\Upsilon(P_N)\right)-\Upsilon(P)-N^{-1}\;(P\,\widehat{\otimes}\,P)\,\nabla^2\Upsilon(P)\|_{\mathcal{B}}$$

$$\leq c_1\;\kappa_3(A)\;\left(1+N^{-\alpha}\,\mathbb{E}\left[\|\mathcal{H}_N\|_F^{p\alpha}\right]^{1/p}\right)\;\mathbb{E}\left[\|\mathcal{H}_N\|_F^{3q}\right]^{1/q}\leq c_2\;\kappa_3(A)$$

This ends the proof of the theorem. ∎

For instance, when $(\mathcal{B}, \|.\|_{\mathcal{B}}) = (\mathcal{S}_r^+, \|.\|_F)$ and $\Upsilon = \psi$ is the square-root functional, we find that

$$\sqrt{N}\ \left[\sqrt{P_N} - \sqrt{P}\right]\ \hookrightarrow_{N\to\infty}\ \int_0^\infty\ e^{-t\psi(P)}\ \mathcal{H}\ e^{-t\psi(P)}\ dt$$

and the negative bias estimate

$$\begin{aligned}
&N^{3/2}\ \|\mathbb{E}\left(\sqrt{P_N}\right) - \sqrt{P} + N^{-1}\ \nabla\psi(P)\cdot\mathbb{E}\left[(\nabla\psi(P)\cdot(\mathbb{X}-P))^2\right]\|_F \\
&\leq \frac{r}{4}\ \lambda_{min}(P)^{-5/2}\ \mathbb{E}\left(\|\mathcal{H}_N\|_F^3\right) \\
&\leq \frac{r}{4}\ \lambda_{min}(P)^{-5/2}\ \left(1+\mathbb{E}\left(\mathrm{Tr}(\mathcal{H}_N^2)\right)\right)
\end{aligned}$$

We conclude that

$$\begin{aligned}
&\left\|\mathbb{E}\left(\sqrt{P_N}\right) - \sqrt{P} + \frac{1}{N}\ \nabla\psi(P)\cdot\mathbb{E}\left[(\nabla\psi(P)\cdot(\mathbb{X}-P))^2\right]\right\|_F \\
&\qquad\qquad \leq \frac{r}{4}\ \frac{1}{N\sqrt{N}}\ \lambda_{min}(P)^{-5/2}\ \left[\mathrm{Tr}(P^2)+\mathrm{Tr}(P)^2\right]
\end{aligned}$$

We check (1.5) using the matrix moment formula (6.6). Choosing

$$(\mathcal{B}, \|.\|_{\mathcal{B}}) = (\mathbb{R}, |.|) \quad\text{and}\quad \Upsilon(A) = \mathrm{Tr}(AB) \Longrightarrow \nabla\Upsilon(A)\cdot H = \mathrm{Tr}(HB)$$

for some given $(r\times r)$-matrix $B$, we find that

$$\begin{aligned}
&\sqrt{N}\ [\mathrm{Tr}(P_N B) - \mathrm{Tr}(PB)] \\
&\overset{law}{=}\ \frac{1}{\sqrt{N}}\sum_{1\leq i\leq N}\mathrm{Tr}\left[(\mathbb{X}_i - P)B\right] \\
&=\ \frac{1}{\sqrt{N}}\sum_{1\leq i\leq N}\left[X_i' B X_i - \mathbb{E}\left[\mathrm{Tr}[X'BX]\right)\right] \\
&\hookrightarrow_{N\to\infty}\quad \mathrm{Tr}(\mathcal{H}B)\sim\mathcal{N}\left[0, 2\ \mathrm{Tr}((PB)^2)\right]
\end{aligned}$$

# 3

# Matrix moment polynomial formulae

## 3.1 Fluctuation and partition-type matrix moments

Note the collection of matrix moments $M_{n,m}(P)$ and $M^{\pm}_{n,m}(P)$ defined for any $1 \leq m \leq n$ by

$$M_{n,m}(P) := M^{-}_{n,m}(P) + M^{+}_{n,m}(P) \tag{3.1}$$

with

$$M^{\pm}_{n,m}(P) := \sum_{\pi \in \mathcal{Q}^{\pm}_{n,m}} M_{\pi}(P)$$

Notice that $M_{n,m}(P)$ is symmetric but $M_{\pi}(P)$ is not necessary symmetric.

We also consider the "un-centered" matrix moments $M^{\circ}_{n,m}(P)$ and $M^{\circ\pm}_{n,m}(P)$ defined as $M_{n,m}(P)$ and $M^{\pm}_{n,m}(P)$ by replacing $M_{\pi}(P)$ by the $\pi$-matrix "un-centered" moment $M^{\circ}_{\pi}(P)$ defined in (1.10). Observe that in this case, the summation indices $(\mathcal{Q}^{+}_{n,m}, \mathcal{Q}^{-}_{n,m})$ are replaced by $(\mathcal{N}_{n,m}, \mathcal{P}_{n,m} - \mathcal{N}_{n,m})$.

For any $1 \leq m \leq \lfloor n/2 \rfloor$ we also have the rather crude estimates

$$\|M_{2n,n}(P)\|_F \leq \upsilon_{2n} \ \left[\mathrm{Tr}(P^2) + \mathrm{Tr}(P)^2\right]^n \tag{3.2}$$

and

$$\|M_{n,m}(P)\|_F \leq \frac{n!}{m!}\; 2^{3n-(2m+1)}\; \mathrm{Tr}(P)^n$$

A derivation of this estimate is technical, thus it is provided in the Appendix.

When $P = I$, integrating sequentially adjacent non-crossing blocks we find that

$$M^+_{2n,n}(I) := C_n\; (1+r)^n\, I \qquad \Longrightarrow \qquad (1.8) \tag{3.3}$$

This formula can also be proved using the Catalan-type recursive formulae (3.15) provided in theorem 3.4. For any $\pi \in \mathcal{N}_{n,m}(\mu)$ we also have

$$\begin{aligned} M^\circ_\pi(I) &= \prod_{1\leq k\leq n} \mathbb{E}(\mathbb{X}^k)^{\mu_k} \\ &= r^{n-m} \prod_{1\leq k\leq n} \left[\prod_{0\leq l<k}\left(1+\frac{2l}{r}\right)\right]^{\mu_k} I \end{aligned} \tag{3.4}$$

This shows that

$$M^\circ_\pi(I) = r^{n-m}\;\left(1+\mathrm{O}(1/r)\right)\, I$$

The proof of the above assertion follows standard calculations for isotropic matrix moments. For instance we can use the $\chi$-square decomposition (A.5) discussed in the Appendix.

For crossing partitions $\pi \in \mathcal{P}_{n,m} - \mathcal{N}_{n,m}$ we also have

$$M^\circ_\pi(I) = c_\pi\; r^{n-m}\,(1+\mathrm{O}(1/r))\;\; r^{-1} I \tag{3.5}$$

for some finite constant $c_\pi$ whose values doesn't depend on $r$. This implies that

$$r^{-(n-m+1)} \sum_{\pi\in\mathcal{N}_{n,m}(\mu)} \mathrm{Tr}(M^\circ_\pi(I)) = N_{n,m}(\mu) + \mathrm{O}(1/r) \tag{3.6}$$

This implies that

$$r^{-(n-m+1)} \sum_{\pi\in\mathcal{N}_{n,m}} \mathrm{Tr}(M^\circ_\pi(I)) = N_{n,m} + \mathrm{O}(1/r) = r^{-(n-m+1)}\; \mathrm{Tr}(M^\circ_{n,m}(I))$$

The proof of (3.5) is provided in the Appendix.

In the same vein, for any $\pi \in \mathcal{N}_{n,m}(\mu)$ we have

$$M_\pi(I)$$

$$= 1_{\mu_1=0} \prod_{2\le k\le n} \left[\frac{1}{r}\, B_k\left(r-1,\ r\ 2^{2-1}\ (2-1)!, \dots,\ r\ 2^{k-1}\ (k-1)!\right)\right.$$

$$\left. +(-1)^k\ \left(1-\frac{1}{r}\right)\right]^{\mu_k} I = 1_{\mu_1=0}\ r^{n-m}\ (1+\mathrm{O}(1/r))\ I \tag{3.7}$$

with the complete Bell polynomials (cf. (2.16)). In addition, for crossing partitions we have

$$\pi \in \mathcal{Q}^-_{n,m} \Longrightarrow M_\pi(I) = c_\pi\ r^{n-m}\ (1+\mathrm{O}(1/r))\ r^{-1}\ I \tag{3.8}$$

for some finite constant $c_\pi$ whose values doesn't depend on $r$. The proof of (3.8) is provided in the Appendix. This implies that

$$r^{-(n-m+1)} \sum_{\pi\in\mathcal{N}_{n,m}(\mu)} \mathrm{Tr}(M_\pi(I))$$

$$= 1_{\mu_1=0}\ N_{n,m}(\mu) + \mathrm{O}(1/r) \tag{3.9}$$

from which we check that

$$r^{-(n-m+1)} \sum_{\pi\in\mathcal{Q}^+_{n,m}} \mathrm{Tr}(M_\pi(I)) = R_{n,m} + \mathrm{O}(1/r)$$

$$= r^{-(n-m+1)}\ \mathrm{Tr}(M_{n,m}(I))$$

with the Narayana and the Riordan numbers introduced in (2.13) and (2.15). The proof of (3.7) is provided in the Appendix.

The computation of the matrix moments of the fluctuation matrices $\mathcal{H}_N$ requires some more sophisticated combinatorial techniques.

Next theorem provides a centered matrix moment interpretation of $\mathbb{E}\left(\mathcal{H}^{2n}\right)$, a polynomial decomposition of centered matrix moments.

**Theorem 3.1.** For any $n \geq 1$ and $N \geq 1$ we have matrix polynomial formulae

$$\begin{aligned}\mathbb{E}\left([P_N - P]^n\right) &= \sum_{1 \leq m \leq \lfloor n/2 \rfloor \wedge N} \frac{(N)_m}{N^n} \, M_{n,m}(P) \\ &= \sum_{1 \leq m \leq \lfloor n/2 \rfloor \wedge N} \frac{1}{N^{n-m}} \, \partial_{n,m}(P)\end{aligned}$$

with the differential-type matrix moments $\partial_{n,m}(P)$ defined in (1.11) and (1.13). In addition, for any $1 \leq m \leq \lfloor n/2 \rfloor$ we also have the inversion formula

$$M_{n,m}(P) = \sum_{1 \leq k \leq m} (-1)^{m-k} \, \frac{k^{n/2}}{k!(m-k)!} \, \mathbb{E}\left(\mathcal{H}_k^n\right) \tag{3.10}$$

and

$$M_{2n,n}(P) = \mathbb{E}\left(\mathcal{H}^{2n}\right)$$

The proof of the above theorem is provided in Section 6.1. Combining (1.4) with the above theorem we compute the covariance matrix moments.

**Corollary 3.2.** For any $2N \geq n \geq 1$ we have

$$\mathbb{E}\left[P_N^n\right] = \sum_{1 \leq m \leq n} \frac{(N)_m}{N^n} \, M^{\circ}_{n,m}(P) = \sum_{1 \leq m \leq n} \frac{1}{N^{n-m}} \, \partial^{\circ}_{n,m}(P)$$

with the differential-type matrix moments given by

$$\begin{aligned}\partial^{\circ}_{n,m}(P) &:= \sum_{1 \leq k \leq (n-m) \wedge m} (n)_{k+(n-m)} \, \partial_{k+(n-m),k}(P) \\ &= \sum_{m \leq k \leq n} s(k,m) \, M^{\circ}_{n,k}(P)\end{aligned}$$

for any $1 \leq m \leq n$ with the convention $\partial^{\circ}_{n,n}(P) = P^n$ when $m = n$.

The proof of the corollary is provided in the Appendix.

**Corollary 3.3.** Let $P = I$ and let $N = r/\rho$ be a scaling of the sample size in terms of the dimension associated with some parameter $\rho > 0$.

$$\begin{aligned}\mathbb{E}\left[P_N^n\right] &= \sum_{1 \leq m \leq n} \rho^{n-m} \, N_{n,m} \, (1 + \mathrm{O}(1/r)) \, I \\ \mathbb{E}\left([P_N - I]^n\right) &= \sum_{1 \leq m \leq \lfloor n/2 \rfloor} \rho^{n-m} \, R_{n,m} \, (1 + \mathrm{O}(1/r)) \, I\end{aligned}$$

The first assertion in corollary 3.3 is a direct consequence of (3.4) and (3.5). The second assertion is a consequence of (3.7) and (3.8). Using (2.13) and taking the trace in the first estimate we obtain the Marchenko–Pastur law (1.10).

The matrix polynomial formulae (1.11) and corollary 3.2 can be extended without further work to matrix moments associated with some collection of matrices $Q_n$ and given in (1.11).

When $N = 1$, and $P = I$ the matrix moments

$$\mathbb{E}\left(H_1^n\right) = \mathbb{E}((\mathbb{X} - I)^n) = M_{n,1}(I)$$

resume to the formula

$$\begin{aligned}\mathbb{E}\left[(\mathbb{X} - I)^n\right] = &\left[\frac{1}{r}\, B_n\left(r-1,\ r\ 2^{2-1}\ (2-1)!, \ldots,\ r\ 2^{n-1}\ (n-1)!\right)\right.\\ &\left. +(-1)^n\ \left(1-\frac{1}{r}\right)\right]\ I\ \longrightarrow_{r\rightarrow\infty}\ 1_{n\geq 1}\end{aligned} \tag{3.11}$$

where $B_n$ stands for the complete Bell polynomials. For more general covariance matrices we have the matrix polynomial formula

$$\begin{aligned} 0\ \leq\ \mathbb{E}\left((\mathbb{X} - P)^n\right)\ &=\ \sum_{1\leq k\leq n} \rho_{k,n}(P)\ P^k \\ &\leq\ \mathbb{E}\left(\mathbb{X}^n\right) - P^n \end{aligned} \tag{3.12}$$

for some non-negative parameters $\rho_{k,n}(P)$ which can be computed explicitly in terms of the traces of the covariance matrix powers. For instance, we find that

$$\begin{aligned}
\mathbb{E}\left((\mathbb{X} - P)^2\right) &= P^2 + P\ \mathrm{Tr}(P)\\
\mathbb{E}\left((\mathbb{X} - P)^3\right) &= 4\ P^3 + 2\ P^2\ \mathrm{Tr}(P) + P\ \left[\mathrm{Tr}(P)^2 + \mathrm{Tr}(P^2)\right]\\
\mathbb{E}\left((\mathbb{X} - P)^4\right) &= 25\ P^4 + 15\ \mathrm{Tr}(P)\ P^3 + 2\ \left[2\ \mathrm{Tr}(P)^2 + 3\ \mathrm{Tr}(P^2)\right]\ P^2\\
&\quad + \left[\mathrm{Tr}(P)^3 + 4\ \mathrm{Tr}(P)\ \mathrm{Tr}(P^2) + 5\ \mathrm{Tr}(P^3)\right]\ P
\end{aligned}$$

The matrix moments of the rank-one matrices $\mathbb{X}$ can also be computed by the formula

$$\frac{2^{-(n-1)}}{(n-1)!}\ \mathbb{E}(\mathbb{X}^n) = \sum_{0\leq k<n} \frac{2^{-k}}{k!}$$

$$\times B_k(\mathrm{Tr}(P),\ 2^{2-1}(2-1)!\ \mathrm{Tr}(P^2),\ldots,2^{k-1}(k-1)!\ \mathrm{Tr}(P^k))\ \ P^{n-k}$$

The proof of (3.12) and the above formulae, including a more refined analysis of these rank-one matrix moments can be found in Section 2 in [10].

The first matrix moments of the fluctuation matrix $\mathcal{H}_N$ are given by the formulae

$$\begin{aligned}\mathbb{E}\left(\mathcal{H}_N^2\right) &= \mathbb{E}\left(\mathcal{H}^2\right) = \mathbb{E}\left((\mathbb{X}-P)^2\right)\\ \mathbb{E}\left(\mathcal{H}_N^3\right) &= N^{-1/2}\ \mathbb{E}\left((\mathbb{X}-P)^3\right)\end{aligned}$$

and

$$\mathbb{E}\left(\mathcal{H}_N^4\right) = N^{-1}\ \mathbb{E}\left((\mathbb{X}-P)^4\right) + \left(1-N^{-1}\right)\mathbb{E}\left(\mathcal{H}^4\right)$$

This yields the fourth polynomial matrix moment formula

$$\begin{aligned}&N\left[\mathbb{E}\left(\mathcal{H}_N^4\right) - \mathbb{E}\left(\mathcal{H}^4\right)\right]\\ &= 20\ P^4 + 12\ \mathrm{Tr}(P)\ P^3 + \left[3\ \mathrm{Tr}(P)^2 + 5\ \mathrm{Tr}(P^2)\right]\ P^2\\ &\qquad + \left[\mathrm{Tr}(P)^3 + 3\ \mathrm{Tr}(P)\ \mathrm{Tr}(P^2) + 4\ \mathrm{Tr}(P^3)\right]\ P\end{aligned} \tag{3.13}$$

The explicit description of all matrix moments $M_{n,m}(P)$ in terms of traces $\mathrm{Tr}(P^k)$ and matrix powers $P^l$ is not known.

## 3.2 A semi-circle law for non-isotropic Wishart matrices

The semi-circle law (1.8) combined with (3.3) and (3.10) imply that

$$P = I \quad\Longrightarrow\quad \lim_{r\to\infty} \mathrm{Tr}\left(\left[\frac{M_{2n,n}(I) - M^+_{2n,n}(I)}{\sqrt{r}}\right]^n\right) = 0$$

This shows that the dominating terms in the semi-circle law fluctuation are encapsulated in the matrix moments $M^{+}_{2n,n}(I)$.

Our second main result is a recursive formula for computing the centered moments associated with non-crossing partitions for a general covariance matrix. To describe with some precision this result we need to introduce some notation. For any $m \geq 0$ and any $v = (v_1, \ldots, v_m)$ with $v_i \geq 0$ we set

$$\mathrm{Tr}_v\left(P\right) := \prod_{1\leq i\leq m} \mathrm{Tr}\left(P^i\right)^{v_i} \quad \text{with the convention} \quad \mathrm{Tr}_0\left(P\right) = 1$$

We let $\mathcal{T}_{p,q}$, $p \geq q$, be the vector space of matrix polynomials in the variable $P$ with maximal degree $p$ and at most $q$ trace-product; that is the vector space spanned by the monomials

$$\mathrm{Tr}_v\left(P\right)\ P^l \quad \text{with} \quad \textstyle\sum_{1\leq i\leq k}\ v_i \leq q$$

$$\text{and} \quad \textstyle\sum_{1\leq i\leq k} i v_i + l = p \quad \text{for some} \quad k \geq 0 \quad \text{and} \quad 1 \leq l\ \leq p$$

We also let $\mathcal{T}^{+}_{p,q} \subset \mathcal{T}_{p,q}$ be the subset defined as above and considering exactly $q$ trace product terms. For instance we have

$$\mathbb{E}(\mathcal{H}^2) = \mathbb{E}(H_N^2) \in \mathcal{T}_{2,1} \quad \text{and} \quad \mathbb{E}(\mathcal{H}^4) \in \mathcal{T}_{4,2} \quad \text{but} \quad \mathbb{E}(\mathcal{H}_N^2) \in \mathcal{T}_{4,3}$$

We also have the decomposition

$$\mathbb{E}\left(\mathcal{H}^4\right) - \overbrace{\left[\mathrm{Tr}(P)^2\, P^2 + \mathrm{Tr}(P)\,\mathrm{Tr}(P^2)\, P\right]}^{\in \mathcal{T}^{+}_{4,2}}$$

$$= 5\, P^4 + 3\ \mathrm{Tr}(P)\, P^3 + \mathrm{Tr}(P^2)\, P^2 + \mathrm{Tr}(P^3)\, P \in \mathcal{T}_{4,1}$$

Consider a collection $P : r \mapsto P(r)$ of possibly random matrices satisfying the condition (1.16). In this situation we have

$$r^{-3}\ \mathbb{E}\left(\mathrm{Tr}(\mathcal{H}^4)\right) \quad \longrightarrow_{r\to\infty} \quad \sigma_2(P) := C_2\ \tau_1(P)^2\ \tau_2(P)$$

Our next objective is to extend this result to any trace moment. The next theorem provides an explicit description of the dominating terms in the matrix moments $M_{2n,n}$.

**Theorem 3.4.** For any $n \geq 1$, the matrix polynomial

$$M^+_{2n,n}(P) \in \mathcal{T}_{2n,n}$$

are given by the Catalan-type recursive formulae

$$\begin{aligned} M^+_{2n,n}(P) &= \sum_{k+l=n-1} \left[M^+_{2k,k}(P)P + \mathrm{Tr}(M^+_{2k,k}(P)P)\ I\right]\ M^+_{2l,l}(P)P \\ M^+_{0,0}(P) &= I \end{aligned} \tag{3.14}$$

We also have the decompositions

$$M^-_{2n,n}(P) \in \mathcal{T}_{2n,n-2} \qquad \text{and} \qquad M^+_{2n,n}(P) \in \Sigma_n(P) + \mathcal{T}_{2n,n-1}$$

with the matrix polynomial $\Sigma_n(P) \in \mathcal{T}^+_{2n,n}$ defined by the Catalan-type recursive formulae

$$\begin{aligned} \Sigma_{n+1}(P) &= \sum_{k+l=n} \mathrm{Tr}(\Sigma_k(P)P)\ \Sigma_l(P)P \\ \Sigma_0(P) &= I \end{aligned} \tag{3.15}$$

The proof of this theorem is provided in Section 6.2. The above theorem clearly yields the convergence result

$$(1.16) \quad \Longrightarrow \quad \lim_{r\to\infty} r^{-1}\,\mathbb{E}\left(\mathrm{Tr}\left(\left[\frac{\mathcal{H}}{\sqrt{r}}\right]^{2n}\right)\right) = \lim_{r\to\infty} r^{-(n+1)}\ \mathrm{Tr}(\Sigma_n(P))$$

with the matrix polynomials $\Sigma_n(P)$ which can be computed sequentially using formulae (3.15). For instance, we have

$$\begin{aligned} \Sigma_1(P) &= \mathrm{Tr}(P)\ P \\ \Sigma_2(P) &= \mathrm{Tr}(P)^2\ P^2 + \mathrm{Tr}(P^2)\ \mathrm{Tr}(P)\ P \\ \Sigma_3(P) &= \mathrm{Tr}(P)^3\ P^3 + 2\ \mathrm{Tr}(P^2)\ \mathrm{Tr}(P)^2\ P^2 \\ &\qquad + \left[\mathrm{Tr}(P^2)^2\ \mathrm{Tr}(P) + \mathrm{Tr}(P^3)\ \mathrm{Tr}(P)^2\right]\ P \end{aligned}$$

and

$$\begin{aligned} \Sigma_4(P) &= \mathrm{Tr}(P)^4\ P^4 + 3\ \mathrm{Tr}(P^2)\ \mathrm{Tr}(P)^3\ P^3 \\ &+ [2\ \mathrm{Tr}(P^3)\ \mathrm{Tr}(P)^3 + 3\ \mathrm{Tr}(P^2)^2\ \mathrm{Tr}(P)^2]\ P^2 \\ &+ [\mathrm{Tr}(P^4)\ \mathrm{Tr}(P)^3 + 3\ \mathrm{Tr}(P^3)\ \mathrm{Tr}(P^2)\ \mathrm{Tr}(P)^2 + \mathrm{Tr}(P^2)^3\ \mathrm{Tr}(P)]\ P \end{aligned}$$

Working a little harder we check that

$$
\begin{aligned}
&\Sigma_5(P)\\
&= \mathrm{Tr}(P)^5\ P^5 + 4\ \mathrm{Tr}(P^2)\ \mathrm{Tr}(P)^4\ P^4\\
&+ \left[6\ \mathrm{Tr}(P)^3\ \mathrm{Tr}(P^2)^2 + 3\ \mathrm{Tr}(P^3)\ \mathrm{Tr}(P)^4\right]\ P^3\\
&+ \left[2\ \mathrm{Tr}(P^4)\ \mathrm{Tr}(P)^4 + 8\ \mathrm{Tr}(P^3)\ \mathrm{Tr}(P^2)\ \mathrm{Tr}(P)^3\right.\\
&\qquad\qquad\qquad\qquad \left. +4\ \mathrm{Tr}(P^2)^3\ \mathrm{Tr}(P)^2\right]\ P^2\\
&+ \left[\mathrm{Tr}(P^5)\ \mathrm{Tr}(P)^4 + 4\ \mathrm{Tr}(P^4)\ \mathrm{Tr}(P^2)\ \mathrm{Tr}(P)^3 + 2\ \mathrm{Tr}(P^3)^2\ \mathrm{Tr}(P)^3\right.\\
&\qquad\qquad \left. +6\ \mathrm{Tr}(P^3)\ \mathrm{Tr}(P^2)^2\ \mathrm{Tr}(P)^2 + \mathrm{Tr}(P^2)^4\ \mathrm{Tr}(P)\right]\ P
\end{aligned}
\tag{3.16}
$$

Next theorem provides a closed form multinomial formula to compute the matrix polynomial functions $\Sigma_n(P)$.

**Theorem 3.5.** The matrix polynomial $\Sigma_n(P) \in \mathcal{T}^+_{2n,n}$ defined in (3.15) are given by the formula (1.14). In addition, we have the multinomial-trace-type formula

$$
\mathrm{Tr}\left(\Sigma_n(P)\right) = 2 \sum_{\mu} \begin{pmatrix} n \\ \mu_1\ \mu_2\ \ldots\ \mu_n \end{pmatrix} \mathrm{Tr}_\mu(P) \tag{3.17}
$$

with

$$
\mathrm{Tr}_\mu(P) := \prod_{1\leq i\leq n} \mathrm{Tr}\left(P^i\right)^{\mu_i}
$$

In the above display, the summation is taken over all $\mu$ satisfying (1.19).

The proof of this theorem is given in Section 6.3.

**Corollary 3.6.** Let $P : r \mapsto P(r)$ be possibly random matrices satisfying condition (1.16). For any $n \geq 1$, the functional $\sigma_n(P)$ defined by (1.18) in corollary 1.3, obeys,

$$
\sigma_n(P) = \lim_{r\to\infty} \mathrm{Tr}\left(\Sigma_n(P)\right)
$$

## 3.3 A Marchenko–Pastur law for non-isotropic Wishart matrices

The aim of this section is to provide a matrix version of the Marchenko–Pastur law for non-isotropic Wishart matrices. For any non-crossing partition $\pi \in \mathcal{N}_{n,m}$ we have the matrix moment formulae

$$\mathbb{X}_{\pi} \stackrel{law}{=} Q_1\mathbb{X}Q_2\ldots\mathbb{X}Q_{|\pi_m|}\mathbb{X} \Longrightarrow \mathbb{X}_{\mathrm{cl}(\pi)} \stackrel{law}{=} \mathbb{X}Q_1\mathbb{X}Q_2\ldots\mathbb{X}Q_{|\pi_m|}\mathbb{X}$$

This implies that

$$M^{\circ}_{\mathrm{cl}(\pi)}(P) = \mathbb{E}\left(\mathbb{X}Q_1\mathbb{X}Q_2\ldots\mathbb{X}Q_{|\pi_m|}\mathbb{X}\right) = \mathbb{E}\left(\left[\prod_{1\leq i\leq|\pi_m|}\langle X, Q_iX\rangle\right]\ \mathbb{X}\right)$$

In the above display, $(Q_i)_{1\leq i\leq|\pi_m|}$ stands for some (non-necessarily symmetric) random matrices independent of $\mathbb{X}$.

On the other hand, for any symmetric matrices $(S_i)_{1\leq i\leq k}$ we also have the reduction formula

$$\mathbb{E}\left(\mathbb{X}\left[S_1\mathbb{X}S_2\ldots\mathbb{X}S_k\mathbb{X}\right]\right) = \mathrm{Tr}\left(\mathbb{E}\left[S_1\mathbb{X}S_2\ldots\mathbb{X}S_k\mathbb{X}\right]\right)\ P$$
$$+\sum_{1\leq i\leq k}\mathbb{E}\left(\mathbb{X}S_1\ldots\mathbb{X}S_{i-1}\mathbb{X}S_{i+1}\ldots\mathbb{X}S_k\mathbb{X}S_i\right)P$$

A proof of this decomposition can be found in [10]. This implies that

$$M^{\circ}_{\mathrm{cl}(\pi)}(P) = \mathrm{Tr}\left(M^{\circ}_{\pi}(P)\right)P + \sum_{1\leq i\leq|\pi_m|} M^{\circ}_{\varphi_i(\pi)}(P)P$$

for some non-crossing partitions $\varphi_i(\pi) \in \mathcal{N}_{n,m}$. Using (2.21) and arguing as in the proof of theorem 3.4 we have

$$M^{\circ}_{n,m}(P) = \Sigma^{\circ}_{n,m}(P) + \mathrm{O}(r^{n-m-1})\ I \tag{3.18}$$

with the matrices $\Sigma^{\circ}_{n,m}(P)$ defined sequentially by the formula

$$\Sigma^{\circ}_{n+1,m+1}(P) = \Sigma^{\circ}_{n,m}(P)P + \sum_{n_1,n_2,m_1,m_2}\Sigma^{\circ}_{n_1,m_1}(P)\ \mathrm{Tr}\left(\Sigma^{\circ}_{n_2,m_2}(P)\right)P \tag{3.19}$$

The summation in the above display is taken over all non-negative parameters $n_1, n_2, m_1, m_2$ satisfying (2.22).

**Theorem 3.7.** The matrices $\Sigma^{\circ}_{n,m}(P)$ defined in (3.19) are given by formula (1.15).

The proof of the above theorem is provided in Section 6.4. Combining (1.12) with (3.18) we obtain the Marchenko–Pastur law for non-isotropic Wishart matrices stated in corollary 1.4.

## 3.4 Some combinatorial formulae

The theorems stated above can be used to derive several new combinatorial formulae. We end this section around this theme. Firstly using (3.17) we have

$$\begin{aligned} C_n &= \sum_{1\leq m\leq n} N_{n,m} \\ &= \sum_{1\leq m\leq n} \sum_{\nu\vdash n\ |\nu|=m} K(\nu) = 2\sum_{\mu} \binom{n}{\mu_1\ \mu_2\ \ldots\ \mu_n} \end{aligned}$$

with the last sum running over the same set of indices as in (3.17).

We fix the covariance matrix $P$ and we set $u(i) = \mathrm{Tr}\left(P^{1+i}\right)$ for any $i\geq 0$. In this notation, by (1.14) for any $n\geq 0$ we have

$$\Sigma_n(P) = \sum_{0\leq m\leq n} \alpha_n(m)\ P^m$$

with

$$\alpha_n(m) := \frac{u(0)}{u(m)} \sum_{\pi\in\mathcal{N}_n\ :\ |\pi_1|=m}\ \prod_{i\geq 0} u(i)^{r_i(\pi)} \tag{3.20}$$

Observe that $\alpha_n(0) = 1_{n=0}$. In addition, using (3.15) we check that

$$\alpha_{n+1}(1) = \sum_{1\leq m\leq n} \alpha_n(m)\ u(m)$$

and the matrix formula

$$\begin{bmatrix} \alpha_{n+1}(2) \\ \alpha_{n+1}(3) \\ \alpha_{n+1}(4) \\ \vdots \\ \alpha_{n+1}(n+1) \end{bmatrix}$$

$$= \begin{bmatrix} \alpha_n(1) & \alpha_{n-1}(1) & \ldots & \alpha_3(1) & \alpha_2(1) & \alpha_1(1) \\ \alpha_n(2) & \alpha_{n-1}(2) & \ldots & \alpha_3(2) & \alpha_2(2) & 0 \\ \alpha_n(3) & \alpha_{n-1}(3) & \ldots & \alpha_3(3) & 0 & \vdots \\ \vdots & \vdots & \ldots & & \vdots & \vdots \\ \alpha_n(n) & 0 & & & & 0 \end{bmatrix} \begin{bmatrix} \alpha_1(1) \\ \alpha_2(1) \\ \alpha_3(1) \\ \vdots \\ \alpha_n(1) \end{bmatrix}$$

This provides a rather simple way to compute the coefficients $\alpha_n(m)$ in terms of $\alpha_n(1)$; that is we have that

$$\forall 1 \leq m \leq n \qquad \alpha_n(m) = \sum_{k_1+\ldots+k_m=n,\ k_i\geq 1} \ \prod_{1\leq i\leq m} \alpha_{k_i}(1)$$

For $m = 1$ the above formula doesn't give any information. To compute the evolution of the terms $\alpha_n(1)$ we use the "highly nonlinear" Bell-type induction

$$\begin{aligned} \alpha_{n+1}(1) &= \sum_{1\leq m\leq n} u(m) \sum_{k_1+\ldots+k_m=n,\ k_i\geq 1} \ \prod_{1\leq i\leq m} \alpha_{k_i}(1) \\ &= \frac{m!}{n!} \, B_{n,m}\left(1!\,\alpha_1(1),\ 2!\,\alpha_2(1), \ldots, n!\alpha_n(1)\right) \end{aligned}$$

with the partial Bell polynomials defined in (2.17).

By construction, the solution of these equations is given by (3.20), for any given sequence $u = (u(i))_{i\geq 0}$. For instance when $u(i) = 1$ we obtain the formulae

$$C_n = \sum_{1\leq m\leq n} C_{n,m}$$

with the Catalan triangle $C_{n,m}$ defined for any $1 \leq m \leq n$ by

$$\begin{aligned} C_{n,m} &:= \sum_{m+l_1+\ldots+l_m=n} \ \prod_{1\leq i\leq m} C_{l_i} \\ &= \left|\{\pi \in \mathcal{N}_n \ : \ |\pi_1| = m\}\right| = \frac{m}{n} \binom{2n}{n-m} \end{aligned}$$

In the above display, the summation run over non-negative indices $l_i \geq 0$. The Catalan triangle has been introduced by Shapiro in [76], see also [17] for a review of some combinatorial properties.

# 4

# Laplace matrix transforms

## 4.1 Some matrix moment estimates

From the practical viewpoint, the combinatorial complexity of the matrix moments (1.11) requires to find some useful computable estimates. We set

$$\epsilon_n(N) := \frac{1}{N^{n/2-\lfloor n/2\rfloor}} \sum_{1\leq m\leq \lfloor n/2\rfloor\wedge N} \left(\frac{4}{N}\right)^{\lfloor n/2\rfloor-m} \frac{\lfloor n/2\rfloor!}{m!} \tag{4.1}$$

and

$$\upsilon_n := \frac{1}{2^{\lfloor n/2\rfloor}} \frac{n!}{\lfloor n/2\rfloor!}$$

When $N \geq n$ we find that

$$\epsilon_n(N) = \begin{cases} 1+2n\ N^{-1}+\mathrm{O}\left(N^{-2}\right) & \text{if} \quad n \text{ is even} \\ \\ N^{-1/2}+\mathrm{O}\left(N^{-3/2}\right) & \text{if} \quad n \text{ is odd} \end{cases}$$

In this notation, we have the following estimates.

**Theorem 4.1.** For any $n \geq 1$ and any $N \geq 1$ we have the estimates

$$\|\mathbb{E}\left(\mathcal{H}_N^n\right)\|_F \leq \frac{\upsilon_n}{2^{1+\lfloor n/2\rfloor-3n}}\ \mathrm{Tr}(P)^n\ \epsilon_n(N) \tag{4.2}$$

In addition, for any $N \geq n$ we have

$$\begin{aligned} & \mathrm{Tr}\left(\mathbb{E}\left[\mathcal{H}_N^{2n}\right]\right) \\ & \leq \mathrm{Tr}\left(\mathbb{E}\left[\mathcal{H}^{2n}\right]\right) + \sqrt{r}\ 2^{12n-1}\ \upsilon_{2n}\ \mathrm{Tr}(P)^{2n}\quad [\epsilon_{2n}(N) - 1] \end{aligned} \tag{4.3}$$

and

$$\begin{aligned} & N\,\|\mathbb{E}\left[\mathcal{H}_N^{2n}\right] - \mathbb{E}\left[\mathcal{H}^{2n}\right]\|_F \\ & \leq (n-1)^2\ \upsilon_{2n}\ \mathbb{E}\left[\|\mathcal{H}\|_F^2\right]^n + 2^{6n-1}\ \upsilon_{2n}\ \mathrm{Tr}(P)^{2n}\epsilon_{2n}(N) \end{aligned} \tag{4.4}$$

The proof of this theorem is provided in Section 6.5.

## 4.2 Rank one matrices

Recalling that $\mathbb{X}$ has a Wishart distribution with one degree of freedom and covariance matrix $P$ we check the integral matrix formula

$$\begin{aligned} & \mathbb{E}\left(\exp\left[t\,\mathbb{X}\right]\right) \\ & = I + \int_0^t \exp\left[\frac{1}{2}\sum_{n\geq 1}\ \frac{(2s)^n}{n}\ \mathrm{Tr}(P^n)\right]\ [I - 2sP]^{-1}P\ ds \end{aligned} \tag{4.5}$$

for any $t \in \mathbb{R}$ s.t. $I - 2tP$ is invertible. The proof of (4.5) is provided in the Appendix. When $P = I$ the above formula reduces to

$$\mathbb{E}\left(\exp\left[t\,\mathbb{X}\right]\right) = \left(1 + \frac{1}{r\ (1-2t)^{r/2}}\right)\ I$$

Our next objective is to compute the log-Laplace matrix transforms

$$\mathbb{L}_P(t) := \log E_P(t) \quad \text{with} \quad E_Q(t) := \mathbb{E}\left(\exp\left[t\left(\mathbb{X} - Q\right)\right]\right)$$

We consider the positive mappings

$$\varpi_1(Q) := Q\ \left[\frac{1}{5}\ I\ + \frac{1}{3}\ Q + \ Q^2\right] \geq \varpi(Q) := Q\ \left[\frac{1}{5}\ I\ + \frac{5}{24}\ Q + \ Q^2\right]$$

We also consider the convex function $L\ :\ t \in [0,1[ \mapsto \mathbb{R}_+$ defined by

$$L(t) := \frac{t^2}{1-t}$$

Observe that

$$\begin{aligned} u\in\mathbb{R}_+\mapsto L^\star(u) &:= \sup_{t\in[0,1[}\,(ut-L(t)) = \left(\sqrt{u+1}-1\right)^2 \\ \Longrightarrow u\in\mathbb{R}_+\mapsto \left(L^\star\right)^{-1}(u) &= u+2\sqrt{u} \end{aligned} \tag{4.6}$$

**Theorem 4.2.** For any $0\leq t<2\mathrm{Tr}(P)$ we have

$$I\leq E_0(t)\;\exp\left[-t\,P\right]\leq\;\exp\left[L\left(2t\,\mathrm{Tr}(P)\right)\;\varpi_1\left(P/\mathrm{Tr}(P)\right)\right] \tag{4.7}$$

In addition we have the centered Laplace estimates

$$I\leq E_p(t)\leq\;\exp\left[L\left(2t\,\mathrm{Tr}(P)\right)\;\varpi\left(P/\mathrm{Tr}(P)\right)\right] \tag{4.8}$$

The proof of this theorem is provided in Section 6.6. Recalling that

$$0<A\leq B\Longrightarrow\log\left(A\right)\leq\log\left(B\right)\quad\text{and}\quad A=\log\exp\left(A\right)$$

for any symmetric matrix $A$, we prove the following corollary.

**Corollary 4.3.** For any $0\leq t<2\mathrm{Tr}(P)$ we have the Laplace matrix inequalities

$$0\leq\mathbb{L}_P(t)-t\,P\leq L\left(2t\mathrm{Tr}(P)\right)\;\varpi_1\left(P/\mathrm{Tr}(P)\right) \tag{4.9}$$

$$0\leq\mathbb{L}_P(t)\leq L\left(2t\mathrm{Tr}(P)\right)\;\varpi\left(P/\mathrm{Tr}(P)\right) \tag{4.10}$$

## 4.3 Laplace matrix estimates

The objective of this section is to estimate the Laplace matrix transforms defined by

$$\mathcal{E}_N(t):=\mathbb{E}\left(\exp\left(t\mathcal{H}_N\right)\right)\quad\text{and}\quad\mathcal{E}(t):=\mathbb{E}\left(\exp\left(t\mathcal{H}\right)\right)$$

The existence of $\mathcal{E}(t)$ is ensured by the estimate

$$\log\|\mathcal{E}(t)\|_F\leq\frac{t^2}{2}\;\left[\mathrm{Tr}(P^2)+\mathrm{Tr}(P)^2\right]$$

which is valid for any $t\in\mathbb{R}$. This assertion is a consequence of the technical lemma 6.3.

**Theorem 4.4.** For any $0 \leq 4\sqrt{2}\, t\, \mathrm{Tr}(P) < 1$ and $N \geq 8$ we have the estimates

$$\begin{aligned} \|\mathbb{E}\left[\cosh\left(t\,\mathcal{H}_N\right)\right] - \mathcal{E}(t)\|_F &\leq\ (8t)^4\ \left[\mathrm{Tr}(P^2) + \mathrm{Tr}(P)^2\right]^2 / N \\ \|\mathbb{E}\left[\sinh\left(t\,\mathcal{H}_N\right)\right]\|_F &\leq 4\ \left(t\,\mathrm{Tr}(P)\right)^3\ /\ \sqrt{N} \end{aligned} \tag{4.11}$$

The proof of this theorem is provided in Section 6.7. Recalling that

$$\mathbb{E}\left[\exp\left(t\,\mathcal{H}_N\right)\right]\ = \mathbb{E}\left[\cosh\left(t\,\mathcal{H}_N\right)\right] + \mathbb{E}\left[\sinh\left(t\,\mathcal{H}_N\right)\right]$$

the estimates (4.11) readily yields the following corollary.

**Corollary 4.5.** For any $0 \leq 4\sqrt{2}\, t\, \mathrm{Tr}(P) < 1$ and $N \geq 8$ we have

$$\|\mathcal{E}_N(t) - \mathcal{E}(t)\|_F \leq \frac{4^6}{N}\, t^4\ \left[\mathrm{Tr}(P^2) + \mathrm{Tr}(P)^2\right]^2 + \frac{4}{\sqrt{N}}\ \left(t\,\mathrm{Tr}(P)\right)^3 \tag{4.12}$$

The Laplace matrix estimate (4.12) differs from trace-type estimates which can be derived using Lieb's inequality

$$\mathbb{E}\left(\mathrm{Tr}\left(\exp\left[A+B\right]\right)\ \mid A\right) \leq \mathrm{Tr}\left(\exp\left[A + \log\mathbb{E}\left(e^B \mid A\right)\right]\right) \tag{4.13}$$

which is valid for any possibly random symmetric matrices $A, B$, see for instance [32, 52, 93]. Applying sequentially the estimate (4.13) and using the spectral mapping theorem we obtain the following spectral-sub-additivity estimate.

**Lemma 4.6** (Tropp [93])**.** For any sequence of independent $(r \times r)$ symmetric matrices $A_i$ indexed by some finite set $i \in I$ we have

$$\mathrm{Tr}\left[\mathbb{E}\left(\exp\sum_{i\in I} A_i\right)\right]\ \leq \mathrm{Tr}\left(\exp\left[\sum_{i\in I}\log\mathbb{E}\left(\exp\left(A_i\right)\right)\right]\right) \tag{4.14}$$

Observe that

$$\varpi(Q) \leq \varpi_+(Q)\ I \quad \text{with} \quad \varpi_+(Q) := \lambda_1(\varpi(Q))$$

Using (4.10) and applying this lemma to the collection of matrices

$$A_i = \frac{t}{\sqrt{N}}\ \left(\mathbb{X}_i - P\right) \Longrightarrow \sum_{1\leq i\leq N} A_i = t\ \mathcal{H}_N$$

we check that

$$\begin{aligned}\operatorname{Tr}\left(\mathcal{E}_N(t)\right) &\leq \operatorname{Tr}\left(\exp\left(N\,\mathbb{L}_P(t/\sqrt{N})\right)\right)\\ &\leq r\ \exp\left(NL\left[2t\operatorname{Tr}(P)/\sqrt{N}\right]\ \varpi_+\left[P/\operatorname{Tr}(P)\right]\right)\end{aligned}$$

These trace estimates are often used to derive exponential concentration estimates.

For instance, using conventional large deviation techniques we check that for any $N \geq 1$ and $\delta > 0$ the probability of the event

$$\begin{aligned}&\lambda_1(H^N)\\ &\leq 2\operatorname{Tr}(P)\left[\frac{\delta+\log(r)}{\sqrt{N}} + 2\sqrt{(\delta+\log(r))\ \varpi_+\left(\frac{P}{\operatorname{Tr}(P)}\right)}\right]\end{aligned} \tag{4.15}$$

is greater than $1 - e^{-\delta}$. For the convenience of the reader a proof of this assertion is provided in the Appendix.

# 5

# Concentration and trace-type estimates

## 5.1 A fluctuation theorem

For any symmetric matrix $A$ and for any $t \in \mathbb{R}$ we have

$$\log \mathbb{E}\left(\exp\left[t\,\mathrm{Tr}(A\mathcal{H})\right]\right) = t^2\ \mathrm{Tr}((AP)^2) \tag{5.1}$$

In other words, $\mathrm{Tr}(A\mathcal{H})$ is a centered Gaussian random variable with variance $2\mathrm{Tr}((AP)^2)$. We also have the $\chi$-square Laplace formula

$$\log \mathbb{E}\left(\exp\left(t\|\mathcal{H}\|_F^2\right)\right) = \frac{1}{4}\sum_{n\geq 1} \frac{(4t)^n}{n}\ \left[\mathrm{Tr}(P^n)^2 + \mathrm{Tr}(P^{2n})\right] \tag{5.2}$$

In addition, there exists some $t > 0$ and some finite constant $c < \infty$ such that

$$\sup_{N\geq 1} \mathbb{E}\left[\exp\left(t\,\|\mathcal{H}_N\|_F\right)\right] < \infty \quad \Longleftrightarrow \quad \forall n \geq 1 \quad \sup_{N\geq 1} \mathbb{E}\left[\|\mathcal{H}_N\|_F^n\right] \leq c^n\ n! \tag{5.3}$$

The first two assertions comes from the fact that $\mathcal{H}$ is a matrix with Gaussian entries.

The proofs of (5.3) and the Laplace formulae (5.1) and the Hardy's exponential property (5.2) are provided in the Appendix.

The Laplace estimate (5.1) readily implies that any $\delta \geq 0$ the probability of the event

$$\mathrm{Tr}(A\mathcal{H}) \leq 2\sqrt{\delta \ \mathrm{Tr}((AP)^2)} \quad \text{is greater than } 1 - e^{-\delta}. \tag{5.4}$$

The proof of this assertion is based on conventional Cramér–Chernov arguments. A brief summary of large deviation techniques is provided in Section 2.1.

Recalling that $P_N$ has a Wishart distribution with $N$ degree of freedom and covariance matrix $N^{-1}P$, for any $2|t|\|AP\|_F \leq \sqrt{N}$ we have

$$\begin{aligned} &\log \mathbb{E}\left(\exp\left(t \ \mathrm{Tr}\left(A\mathcal{H}_N\right)\right)\right) \\ &= t^2 \ \mathrm{Tr}((AP)^2) + t^2 \sum_{n\geq 1} \frac{t^n}{n+2} \frac{2^{n+1}}{N^{n/2}} \ \mathrm{Tr}((AP)^{2+n}) \end{aligned} \tag{5.5}$$

Explicit formulae for $\mathrm{Tr}((AP)^n)$ are rarely available. For rank one matrices $A = (xy' + yx')/2$ we have the formula

$$\begin{aligned} &2^{n+1} \ \mathrm{Tr}\left[(AP)^{n+2}\right] \\ &= \sum_{0\leq l\leq \lfloor (n+1)/2\rfloor} \binom{n+1}{2l} \langle x, Py\rangle^{n+2-2l} \ \left(\langle x, Px\rangle \langle y, Py\rangle\right)^l \\ &+ \sum_{1\leq l\leq \lfloor n/2\rfloor -1} \binom{n+1}{2l-1} \langle x, Py\rangle^{n+2-2l} \ \left(\langle x, Px\rangle \langle y, Py\rangle\right)^l \\ &\leq 2^{n+2}\lambda_1(P)^n \end{aligned} \tag{5.6}$$

as soon as $x, y \in \mathbb{B}$, where $\mathbb{B}$ is the unit ball in $\mathbb{R}^r$. A proof of (5.6) is provided in the Appendix. In this situation, we obtain the Laplace estimate

$$\begin{aligned} &\mathbb{E}\left(\exp\left(t \ \mathrm{Tr}\left(A\mathcal{H}_N\right) - t^2 \mathrm{Tr}((AP)^2)\right)\right) \\ &\leq \exp\left(-(2t\lambda_1(P)))^2 \log\left[1 - \frac{1}{N^{1/2}} \ (2t\lambda_1(P))\right]\right) \end{aligned}$$

Whenever $x = u_i$ and $y = u_j$ are two orthonormal eigenvectors associated with some different eigenvalues $\lambda_i(P)$ and $\lambda_j(P)$ with $i \neq j$ we check that

$$2^{n-1}\, \mathrm{Tr}\left[(AP)^n\right] = (\lambda_i(P)\lambda_j(P))^{n/2}$$

To simplify notation, given some collections $P : r \mapsto P(r)$ and $A \ : \ r \mapsto A(r)$ of symmetric matrices, we write

$$(\lambda_1(AP), \mathcal{H}_N, \mathrm{Tr}((AP)^n))$$

instead of the sequence

$$(\lambda_1(A(r)P(r)), \mathcal{H}_N(r), \mathrm{Tr}((A(r)P(r))^n))$$

In this notation, we have the following theorem.

**Theorem 5.1.** For any function $N = N(r)$ s.t. $rN(r) \rightarrow_{r\rightarrow\infty} \infty$ and any collection of symmetric matrices $A \ : \ r \mapsto A(r)$ we have

$$\begin{aligned}
&\sup_{r\geq 1} \lambda_1(AP) < \infty \quad \text{and} \quad \sigma^2_{A,r}(P) := r^{-1}\mathrm{Tr}((AP)^2) \rightarrow_{r\rightarrow\infty} \ \sigma^2_A(P) \\
&\Longrightarrow \quad \frac{1}{\sqrt{2r}}\ \mathrm{Tr}\left(A\mathcal{H}_{N(r)}\right) \hookrightarrow_{r\rightarrow\infty} \mathcal{N}\left(0, \sigma^2_A(P)\right)
\end{aligned}$$

Assume that

$$0 \leq r^{-1}\ |\mathrm{Tr}((AP)^n)| \leq \alpha_A(P)\ \beta_A(P)^n$$

for some parameters $\alpha_A(P) \geq 1$ and $\beta_A(P) \geq 0$ whose values doesn't depend on $n$. In this situation, we have the estimates

$$\begin{aligned}
&\frac{1}{(2n)!}\ |\mathbb{E}\left(\left(\frac{1}{\sqrt{2r}}\ \mathrm{Tr}\left(A\mathcal{H}_N\right)\right)^{2n}\right) - \mathbb{E}\left(\left(\frac{1}{\sqrt{2r}}\ \mathrm{Tr}\left(A\mathcal{H}\right)\right)^{2n}\right)| \\
&\qquad\qquad \leq \ \frac{4(e-1)}{rN}\ \alpha_A(P)^{n-1}\ \beta_A(P)^{2n}\ \left(1 \vee \frac{8}{rN}\right)^{n-1}
\end{aligned}$$

as well as

$$\begin{aligned}
&\frac{1}{(2n+1)!}\ |\mathbb{E}\left(\left(\frac{1}{\sqrt{2r}}\ \mathrm{Tr}\left(A\mathcal{H}_N\right)\right)^{2n+1}\right)| \\
&\leq \frac{(e-1)\sqrt{2}}{\sqrt{rN}}\ \alpha_A(P)^n\ \beta_A(P)^{2n+1}\ \left(1 \vee \frac{8}{rN}\right)^n
\end{aligned}$$

In addition if $A \geq 0$ then the terms in absolute value $|.|$ in the above display are positive.

Under the assumptions of the above theorem we have

$$\frac{1}{\sqrt{2r}}\ \mathrm{Tr}\left(A\mathcal{H}_N\right) \hookrightarrow_{N\to\infty} \frac{1}{\sqrt{2r}}\ \mathrm{Tr}\left(A\mathcal{H}\right) \ \ \hookrightarrow_{r\to\infty} \mathcal{N}\left(0, \sigma_A^2(P)\right)$$

The proof of theorem 5.1 is provided in Section 6.8.

## 5.2 Spectral versus trace-type concentration inequalities

For any symmetric matrix $A$ we have the sub-Gaussian Laplace estimate

$$\begin{aligned} &4|t|\ \|AP\|_F \leq \sqrt{N} \\ &\Longrightarrow \log \mathbb{E}\left(\exp\left(t\ \mathrm{Tr}\left(A\mathcal{H}_N\right)\right)\right) \leq t^2\ \ \left[\mathrm{Tr}((AP)^2) + 2\|AP\|_F^2\right] \end{aligned}$$

as well as

$$\begin{aligned} &4|t|\,\mathrm{Tr}(AP) < \sqrt{N},\ \ A \geq 0 \\ &\Longrightarrow \log \mathbb{E}\left(\exp\left(t\ \mathrm{Tr}\left(A\mathcal{H}_N\right)\right)\right) \leq 3t^2\ \mathrm{Tr}((AP)^2) \end{aligned} \tag{5.7}$$

The proof of these estimates are essentially based on the Laplace formula (5.5). The details of the proof including more refined logarithmic-type estimates are provided in the Appendix.

**Theorem 5.2.** For any symmetric matrix $A$ and any $0 \leq 8\delta \leq N$ the probability of the event

$$|\mathrm{Tr}(A\mathcal{H}_N)| \leq 2\sqrt{(\delta+1)\ \ [\mathrm{Tr}((AP)^2) + 2\|AP\|_F^2]} \tag{5.8}$$

is greater than $1 - e^{-\delta}$. In addition, if $A \geq 0$ then for $0 \leq 8\delta \leq N$ the probability of the event

$$\mathrm{Tr}(A\mathcal{H}_N) \leq 2\sqrt{3}\ \sqrt{\delta\ \mathrm{Tr}((AP)^2)} \quad \text{is greater than } 1 - e^{-\delta}$$

Furthermore for any $0 \leq 4\delta \leq N\ \mathrm{Tr}((AP)^2)/\mathrm{Tr}(AP)^2$ the probability of the event

$$\mathrm{Tr}(A\mathcal{H}_N) \geq -2\ \sqrt{\delta\ \mathrm{Tr}((AP)^2)} \quad \text{is greater than } 1 - e^{-\delta}$$

The proof of the above theorem is provided in Section 6.9. We also have the immediate corollary.

**Corollary 5.3.** For any $0 \leq 4\delta \leq N\left[1/2 \wedge \mathrm{Tr}((AP)^2)/\mathrm{Tr}(AP)^2\right]$ and any matrix $A \geq 0$ the probability of the event

$$|\mathrm{Tr}(AP_N) - \mathrm{Tr}(AP)| \leq 2\sqrt{\frac{\delta+1}{N}\left[2\mathrm{Tr}(AP)^2 + \mathrm{Tr}((AP)^2)\right]}$$

is greater than $1 - e^{-\delta}$.

We end this section with some operator norm estimates which can be deduced almost directly combining the trace concentration inequalities (5.8) with the variational formulation (1.20) of the operator norm.

**Theorem 5.4.** For any $N \geq 8r\log(31)$ we have

$$\mathbb{E}\left[\|\mathcal{H}\|_{op}\right] \wedge \mathbb{E}\left[\|\mathcal{H}_N\|_{op}\right] \leq 6\sqrt{r}\,\lambda_1(P)$$

This also implies that

$$\mathbb{E}\left[\lambda_1(P_N)/\lambda_1(P)\right] \leq 1 + 6\sqrt{r/N}$$

For any $N$ and $\delta$ such that $N/8 \geq \delta + 7r$, the probability of the event

$$\|\mathcal{H}_N\|_{op} \leq 5\sqrt{3}\,\sqrt{\delta+7r}\,\lambda_1(P) \tag{5.9}$$

is greater than $1 - e^{-\delta}$.

**Corollary 5.5.** For any $N$ and $\delta$ such that $N/8 \geq \delta + 7r$, the probability of the event

$$\sup_{1\leq k\leq r} |\lambda_k(P_N) - \lambda_k(P)| \leq 5\sqrt{3}\sqrt{\frac{\delta+7r}{N}}\,\lambda_1(P)$$

is greater than $1 - e^{-\delta}$.

The proof of theorem 5.4 and corollary 5.5 is given in Section 6.10.

# 6

# Proof of the main results

## 6.1 Proof of theorem 3.1

When $m > n/2$ at least one block of $\pi \in \mathcal{P}_{n,m}$ is of unit size and the matrix moment $M_\pi(P) = 0$ resume to the null matrix. This implies that

$$\forall \lfloor n/2 \rfloor < m \leq n \qquad M_{n,m}(P) = 0 \tag{6.1}$$

Let $\langle p, N\rangle$ be the set of all the $(N)_p := N!/(N-p)!$ one to one mappings $\beta$ from $[p] := \{1, \ldots, p\}$ into $[N]$. In this notation we have the decomposition

$$[N]^{[n]} = \cup_{1 \leq p \leq (n \wedge N)} \cup_{\pi : |\pi| = p} \{\alpha^\pi_\beta \ : \ \beta \in \langle p, N\rangle\} \tag{6.2}$$

with the mapping $\alpha^\pi_\beta$ defined by

$$\alpha^\pi_\beta = \sum_{1 \leq i \leq |\pi|} \beta(i) \ 1_{\pi_i}$$

where $\pi_1, \ldots, \pi_{|\pi|}$ stand for the $|\pi|$ blocks of a partition of $[n]$ ordered w.r.t. their smallest element. When $\beta(i) = Id(i) := i$ is the identity mapping sometimes we write $\alpha^\pi$ instead of $\alpha^\pi_{Id}$. A detailed proof of this decomposition can be found in Section 8.6 in [21]. For a further discussion on partitioning and counting onto mappings we also refer to [25]. This decomposition implies that

$$\begin{aligned} & N^{n/2}\,\mathbb{E}\left(H_N^n\right) \\ & = \sum_{1\le m\le \lfloor n/2\rfloor\wedge N} (N)_m\ M_{n,m}(P) \\ & = \sum_{1\le l\le \lfloor n/2\rfloor\wedge N} N^l \sum_{l\le m\le \lfloor n/2\rfloor\wedge N} s(m,l)\ M_{n,m}(P) \end{aligned} \tag{6.3}$$

The second polynomial formulae in (6.3) comes from the formula

$$(N)_m = \sum_{1\le l\le m} s(m,l)\ N^l$$

This ends the proof of the first assertion.

Using (6.3), for any given $n\ge 1$ and any $1\le k\le n$ we also have

$$\mathbb{E}\left(\left[\sqrt{k}\ H_k\right]^{2n}\right) := f_n(k) = \sum_{1\le l\le k} \left(\begin{array}{c} k \\ l \end{array}\right) g_n(l)$$

with

$$g_n(l) := l!\ M_{2n,l}(P)$$

Using the binomial inversion formula we find that

$$\forall 1\le k\le n\quad g_n(k) = \sum_{1\le l\le k} (-1)^{k+l} \left(\begin{array}{c} k \\ l \end{array}\right) f_n(l)$$

Thus, for any $1\le m\le \lfloor n/2\rfloor$ we also have the matrix moment inversion formula

$$M_{n,m}(P) = \sum_{1\le k\le m} (-1)^{m-k}\,\frac{1}{k!(m-k)!}\,\mathbb{E}\left(\left[\sum_{1\le i\le k} (\mathbb{X}_i - P)\right]^n\right)$$

This ends the proof of the theorem. $\square$

## 6.2 Proof of theorem 3.4

For $n=2$ we have $C_2=2$ non-crossing partitions and a single crossing partition

$$N_{2n,n}(0,n,0,\ldots,0) = \{\pi^1,\pi^2\}$$

and

$$\mathcal{P}_{2n,n}(0,n,0,\ldots,0) - N_{2n,n}(0,n,0,\ldots,0) = \{\pi^3\}$$

given by

$$\pi^1 = \{1,2\}\ \{3,4\} \qquad \pi^2 = \{1,4\}\ \{2,3\} \qquad \text{and} \quad \pi^3 = \{1,3\}\ \{2,4\} \tag{6.4}$$

Non-crossing pair partitions can be interpreted in various ways. Firstly, we can interpret non-crossing pair partitions in terms of connections of the ordered set of integers $[2n]$ lying on a horizontal line in the plane with $n$ nonintersecting arcs/lines. In the above example we have the Murasaki diagram

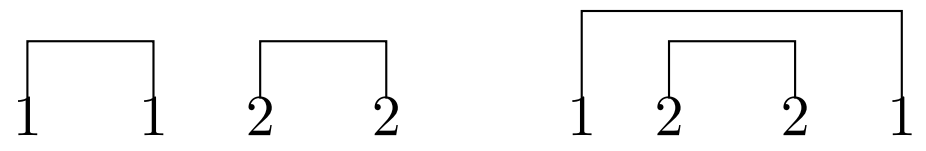

$\pi^1$ $\pi^2$

A given crossing pair partition can be interpreted as balanced parentheses on the ordered set $[2n]$ (with $n$ left and $n$ right). When $n = 2$ we have

$$\pi^1 = \{1,2\}\ \{3,4\} = ()\,() \qquad \text{and} \quad \pi^2 = \{1,4\}\ \{2,3\} = (())$$

They can also be seen as the number of ways of grouping $n+1$ letters with parenthesis

$$\pi^1 = \{1,2\}\ \{3,4\} = ()\,() = ((a_1 \cdot a_2) \cdot a_3)$$

and

$$\pi^2 = \{1,4\}\ \{2,3\} = (()) = (a_1 \cdot (a_2 \cdot a_3)) \tag{6.5}$$

We obtain the balanced parenthesis by removing the letters and the left parenthesis "( " and then replacing each "·" by a left parenthesis "( ". This construction provides a simple way to count the number of non-crossing partitions using Catalan recursion

$$C_2 = C_0\ C_1 + C_1\ C_0$$

For $n = 3$ we have 5 non-crossing pair partitions of 6 elements into pairs. The first $C_2\ C_0 = 2$ partitions

$$\pi_1 = \{1,2\}\{3,4\}\{5,6\} = ()\,()\,() = (((a_1 \cdot a_2) \cdot a_3) \cdot a_4)$$
$$\pi_2 = \{1,3\}\{2,4\}\{5,6\} = (())\,() = ((a_1 \cdot (a_2 \cdot a_3)) \cdot a_4)$$

Then $C_1\ C_1 = 1$ partition given by

$$\pi_3 = \{1,2\}\{2,6\}\{3,4\} = ()\,(()) = ((a_1 \cdot a_2) \cdot (a_3 \cdot a_4))$$

Finally, the last $C_0\ C_2 = 2$ partitions

$$\pi_4 = \{1,6\}\{2,3\}\{4,5\} = (()\,()) = (a_1 \cdot ((a_2 \cdot a_3) \cdot a_4))$$
$$\pi_5 = \{1,6\}\{2,5\}\{3,4\} = ((())) = (a_1 \cdot (a_2 \cdot (a_3 \cdot a_4)))$$

The matrix moments $M_\pi(P)$ associated with pair partition can be sequentially computed using the positive maps defined for any $Q \in \mathcal{S}_r$ by the matrix *linear* functionals

$$\begin{aligned} \Omega(Q) &:= \mathbb{E}\left((\mathbb{X}-P)\mathrm{Tr}\left[Q(\mathbb{X}-P)\right]\right) = 2PQP \\ \Gamma(Q) &:= \mathbb{E}\left((\mathbb{X}-P)Q(\mathbb{X}-P)\right) = \frac{1}{2}\,\Omega(Q) + \overline{\Gamma}(Q) \end{aligned} \tag{6.6}$$

with the mapping

$$\overline{\Gamma}(Q) := \mathrm{Tr}(PQ)\ P$$

When $P = I$ we have

$$\Omega(Q) = 2Q \quad \text{and} \quad \overline{\Gamma}(Q) = \mathrm{Tr}(Q)\ I \tag{6.7}$$

and therefore

$$\Gamma(I) = (1+r)\ I$$

The matrix moments $M_\pi(P)$ associated with non-crossing pair partition are expressed in terms of powers and compositions of the mapping $\Gamma$. For instance, for $n = 2$ the matrix moments associated with the partitions (6.5) are given by adding a $\Gamma$ in front of each left balanced parenthesis and replacing the inside by the identity; that is we have

$$\begin{aligned} M_{\pi^1}(P) &= \mathbb{E}\left((\mathbb{X}-P)^2\right)^2 \\ &= \Gamma\,(I)\ \ \Gamma\,(I) = \Gamma(I)^2 = P^4 + 2\mathrm{Tr}(P)\,P^3 + \mathrm{Tr}(P)^2\,P^2 \end{aligned}$$

and

$$\begin{aligned} M_{\pi^2}(P) &= \mathbb{E}\left((\mathbb{X}-P)\Gamma(I)(\mathbb{X}-P)\right) \\ &= \Gamma\left(\Gamma\left(I\right)\right) \\ &= \Gamma^2\left(I\right) P^4 + \mathrm{Tr}(P)\, P^3 + \left[\mathrm{Tr}(P^3) + \mathrm{Tr}(P)\,\mathrm{Tr}(P^2)\right]\, P \end{aligned}$$

This shows that

$$\begin{aligned} M^+_{2,1}(P) &= M^+_{0,0}(P)\ \Gamma\left(M^+_{0,0}(P)\right) = \Gamma(I) \quad \text{with} \quad M^+_{0,0}(P) = I \\ &= P^2 + \mathrm{Tr}(P)\ P \end{aligned}$$

and

$$\begin{aligned} M^+_{4,2}(P) &= M^+_{0,0}(P)\ \Gamma\left(M^+_{2,1}(P)\right) + M^+_{2,1}(P)\ \Gamma\left(M^+_{0,0}(P)\right) \\ &= I\ \Gamma\left(\Gamma(I)\right) + \Gamma(I)\Gamma(I) \\ &= 2P^4 + 3\mathrm{Tr}(P)\ P^3 + \mathrm{Tr}(P^2)\ P^2 \\ &\qquad\qquad + \left[\mathrm{Tr}(P^3) + \mathrm{Tr}(P^2)\mathrm{Tr}(P)\right]\ P \end{aligned}$$

In the same vein, for $n = 3$ we have the sum of the $C_3 = 5$ terms associated with the balanced parenthesis

$$((()))\ ,\ (()\,())\ ,\ ()\,(())\ ,\ (())\,()\ ,\ ()\,()\,()$$

These balanced strings give the value of the mapping $(\alpha^\pi(1), \ldots, \alpha^\pi(6))$ for the

$$C_3 = C_0\ (C_2) + C_1\ (C_1) + C_2\ (C_0) = 5$$

non-crossing partitions of $[6]$ into 3 blocks:

$$(1,2,3,3,2,1)\quad,\quad(1,2,2,3,3,1)\ ,\ (1,1,2,3,3,2)$$

and

$$(1,2,2,1,3,3)\ ,(1,1,2,2,3,3)$$

This implies that

$$\begin{aligned} M^+_{6,3}(P) &= \Gamma\left(\Gamma^2(I) + \Gamma(I)^2\right) + \Gamma(I)\ \Gamma^2(I) + \left[\Gamma^2(I) + \Gamma(I)^2\right]\ \Gamma(I) \\ &= M^+_{0,0}(P)\ \Gamma\left(M^+_{4,2}(P)\right) + M^+_{2,1}(P)\ \Gamma\left(M^+_{2,1}(P)\right) \\ &\qquad\qquad + M^+_{4,2}(P)\ \Gamma\left(M^+_{0,0}(P)\right) \end{aligned}$$

from which we check that

$$\begin{aligned} M^+_{6,3}(P) &= 5\ P^6 + 10\ \mathrm{Tr}(P)\ P^5 + 6\ \mathrm{Tr}(P)^2\ P^4 \\ &+ \left[3\ \mathrm{Tr}(P^3) + 3\ \mathrm{Tr}(P^2)\ \mathrm{Tr}(P) + \mathrm{Tr}(P)^3\right]\ P^3 \\ &+ 2\left[\mathrm{Tr}(P)\ \mathrm{Tr}(P^3) + \mathrm{Tr}(P^2)\ \mathrm{Tr}(P)^2\right]\ P^2 \\ &+ \left[2\ \mathrm{Tr}(P^5) + 3\ \mathrm{Tr}(P^4)\ \mathrm{Tr}(P) \right. \\ &\quad \left. + \mathrm{Tr}(P^3)\ \mathrm{Tr}(P^2) + \mathrm{Tr}(P^3)\mathrm{Tr}(P)^2 + \mathrm{Tr}(P^2)^2\ \mathrm{Tr}(P)\right]\ P \end{aligned}$$

For $n = 4$ we have the sum of the $C_4 = 14$ terms. The first $C_0\ C_3 = 5$ terms

$$(((()))) \ ,((()\,())) \ , \ (()\,(())) \ , \ ((())\,()) \ , \ (()\,()\,())$$

The $C_1\ C_2 = 2$ and the $C_2\ C_1 = 2$ terms

$$()((())) \ , \ ()(()\,()) \quad \text{and} \quad ((()))() \ , \ (()\,())()$$

and finally the $C_3\ C_0 = 5$ terms

$$((())) \ ,(()\,())\,() \ , \ ()\,(())\,() \ , \ (())\,()\,() \ , \ ()\,()\,()\,()$$

This yields the decomposition

$$\begin{aligned} & M^+_{8,4}(P) \\ &= M^+_{0,0}(P)\ \Gamma(M^+_{6,3}(P)) + M^+_{2,1}(P)\ \Gamma(M^+_{4,2}(P)) \\ &\qquad + M^+_{4,2}(P)\ \Gamma(M^+_{2,1}(P)) + M^+_{6,3}(P)\ \Gamma(M^+_{0,0}(P)) \\ &= \Gamma^2\left(\Gamma^2(I) + \Gamma(I)^2\right) + \Gamma\left(\Gamma(I)\ \Gamma^2(I)\right) + \Gamma\left(\left[\Gamma^2(I) + \Gamma(I)^2\right]\ \Gamma(I)\right) \\ &\quad + \Gamma(I)\ \Gamma\left(\Gamma^2(I) + \Gamma(I)^2\right) + \left(\Gamma^2(I) + \Gamma(I)^2\right)\ \Gamma(\Gamma(I)) \\ &\quad + \Gamma\ (\Gamma^2(I) + \Gamma(I)^2) + \Gamma(I)\ \Gamma^2(I) + ([\Gamma^2(I) + \Gamma(I)^2]\ \Gamma(I))\ \Gamma(I) \end{aligned}$$

More generally, we have the Catalan-type matrix polynomial recursion

$$M^+_{2n,n}(P) := \sum_{k+l=n-1} M^+_{2k,k}(P)\ \Gamma(M^+_{2l,l}(P)) \tag{6.8}$$

The matrix polynomial $M^{+}_{2n,n}(P)$ is a linear combination of monomials of the form

$$\left[\prod_{1\leq i\leq m} \mathrm{Tr}(P^i)^{v_i}\right] P^w \quad \text{with} \quad \left\{\begin{array}{l} 1\leq w\leq 2n \qquad 0\leq m\leq n \\ \\ \sum_{1\leq i\leq m} v_i\leq n \quad (\Longrightarrow m\leq n) \\ \\ \sum_{1\leq i\leq m}\ i\,v_i+w=2n \end{array}\right.$$

In addition, the coefficients don't depend on the matrix $P$. We check this property by induction w.r.t. the parameter $n$ using the recursion (6.8). This ends the proof of (3.14).

We let $M^{[Q_0,\ldots,Q_{2n}],+}_{2n,n}(P)$ the matrix moments associated with a collection of matrices $Q_i$ defined by

$$M^{[Q_0,\ldots,Q_{2n}],+}_{2n,n}(P)$$

$$:= \sum_{\pi\in\mathcal{Q}^{+}_{2n,n}} \mathbb{E}\left(Q_0(\mathbb{X}_{\alpha^{\pi}(1)}-P)Q_1(\mathbb{X}_{\alpha^{\pi}(2)}-P)Q_2\ldots(\mathbb{X}_{\alpha^{\pi}(2n)}-P)Q_{2n}\right) \tag{6.9}$$

For homogeneous sequences $Q_i = Q$ we write $M^{Q,+}_{2n,n}(P)$ instead of $M^{(Q,\ldots,Q),+}_{2n,n}(P)$.

Extending the arguments given above to the matrix polynomials (6.9), we obtain the following technical lemma.

**Lemma 6.1.** We have the functional recursion

$$M^{[Q_0,\ldots,Q_{2n}],+}_{2n,n}(P)$$

$$= \sum_{k+l=n-1} M^{[Q_0,\ldots,Q_{2k}],+}_{2k,k}(P)\ \ \Gamma\left[M^{[Q_{2k+1},\ldots,Q_{2n-1}],+}_{2l,l}(P)\right]\ Q_{2n} \tag{6.10}$$

with the initial condition $M^{Q,+}_{0,0}(P)=Q$.

For instance we have

$$
\begin{aligned}
M_{2,1}^{[Q_0,Q_1,Q_2],+}(P) &= Q_0\Gamma(Q_1)Q_2\\
M_{4,2}^{[Q_0,\ldots,Q_4],+}(P) &= Q_0\Gamma\left[Q_1\Gamma(Q_2)Q_3\right]Q_4 + Q_0\Gamma(Q_1)Q_2\Gamma(Q_3)Q_4\\
M_{6,3}^{[Q_0,\ldots,Q_6],+}(P) &= Q_0\Gamma(Q_1)Q_2\ \Gamma(Q_3)Q_4\ \Gamma(Q_5)Q_6\\
&\quad + Q_0\Gamma\left[Q_1\Gamma(Q_2)Q_3\right]Q_4\ \Gamma(Q_5)Q_6\\
&\quad + Q_0\Gamma(Q_1)Q_2\ \Gamma\left[Q_3\Gamma(Q_4)Q_5\right]Q_6\\
&\quad + Q_0\ \Gamma\left[Q_1\Gamma(Q_2)Q_3\Gamma(Q_4)Q_5\right]Q_6\\
&\quad + Q_0\ \Gamma\left[Q_1\Gamma\left[Q_2\Gamma(Q_3)Q_4\right]Q_5\right]Q_6
\end{aligned}
$$

For any $0 \leq m \leq n$ we let $\Xi_{n,m}$ be the set of partition $\tau = (\tau_1, \ldots, \tau_{m+1})$ of $[2n-1]$ be into $(m+1)$ increasing sequences

$$\tau_j = (i_1^j, \ldots, i_k^j) \quad \text{with} \quad 0 < i_1^j < \ldots < i_k^j < 2n$$

We also set

$$(PQ)_{\tau_j} = (PQ_{i_1^j})(PQ_{i_2^j})\ldots(PQ_{i_k^j})$$

In this notation, using the induction (6.10) and recalling that $\Gamma$ only increases the number of traces by 1, we check that for any $n \geq 1$ the matrix polynomial $M_{2n,n}^{[Q_0,\ldots,Q_{2n}],+}(P)$ is a linear combination of monomials of the form

$$\left[\prod_{1\leq i\leq m} \mathrm{Tr}\left((PQ)_{\tau_i}\right)\right] \quad Q_0\ (PQ)_{\tau_{m+1}}Q_{2n} \tag{6.11}$$

with

$$0 \leq m \leq n \quad \text{and} \quad \tau \in \Xi_{n,m}$$

Arguing as above, we check that the matrix polynomial $M_{2n,n}^{[P],+}(P)$ is a linear combination of monomials of the form

$$\mathrm{Tr}_v(P)\ P^w$$

with the parameters

$$\left\{\begin{array}{lll} 3 \leq w \leq 4n+1 & 0 \leq m \leq n & \sum_{1\leq i\leq m} v_i \leq n \\ \\ \sum_{1\leq i\leq m}\ iv_i + w = 4n+1 & & \end{array}\right.$$

Rewritten in terms of the matrix functional $\overline{\Gamma}$ the Catalan-type recursive formulae (3.15) reduces to the formula

$$\Sigma_n(P) := \sum_{k+l=n-1} \Sigma_k(P)\ \overline{\Gamma}(\Sigma_l(P)) \quad \text{with} \quad \Sigma_0(P) = I$$

A simple induction w.r.t. the parameter $n$ shows that the matrix polynomial $\Sigma_n(P)$ is a linear combination of monomial for $n$ trace terms given by of the form

$$\mathrm{Tr}_v(P)\ P^w$$

with the parameters

$$\begin{cases} 1 \leq w \leq n \qquad 0 \leq m \leq n \qquad \sum_{1\leq i\leq m} v_i = n \\ \\ \sum_{1\leq i\leq m}\ i\, v_i + w = 2n \end{cases}$$

This assertion is also a direct consequence of the formula (1.14) stated in theorem 3.5. Let

$$\widetilde{M}_{2n,n}(P) := M^+_{2n,n}(P) - \Sigma_{2n,n}(P) \Longrightarrow \widetilde{M}_{0,0}(P) = 0$$

Combining the properties discussed above with the easily checked induction

$$\begin{aligned} &\widetilde{M}_{2n,n}(P) \\ &= \sum_{k+l=n-1} \widetilde{M}_{2k,k}(P)\ \overline{\Gamma}\left(\widetilde{M}_{2l,l}(P)\right) + \frac{1}{2} \sum_{k+l=n-1} M_{2k,k}(P)\ \Omega\left(M_{2l,l}(P)\right) \\ &\qquad + \sum_{k+l=n-1} \left[\widetilde{M}_{2k,k}(P)\ \overline{\Gamma}\left(\Sigma_{2l,l}(P)\right) + \Sigma_{2k,k}(P)\ \overline{\Gamma}\left(\widetilde{M}_{2l,l}(P)\right)\right] \end{aligned}$$

we check that the matrix polynomial $\widetilde{M}_{2n,n}(P)$ is a linear combination of monomial for $n$ trace terms given by

$$\mathrm{Tr}_v(P)\ P^w$$

with the parameters

$$\begin{cases} 1 \leq w \leq 2n \qquad 0 \leq m \leq n \qquad \sum_{1\leq i\leq m} v_i < n \\ \\ \sum_{1\leq i\leq m}\ i\, v_i + w = 2n \end{cases}$$

The matrix moments $M_\pi(P)$ associated with crossing pair partitions require more sophisticated calculations. Nevertheless they can be computed in terms of powers and compositions of the mappings $(\Gamma,\Omega)$ using sequentially conditioning formulae. For instance, the $\pi^3$-matrix moment associated with the partition $\pi^3$ discussed in (6.4) is given by

$$\begin{aligned} M_{\pi^3}(P) &= \mathbb{E}\left(\mathbb{E}\left([(\mathbb{X}_1-P)(\mathbb{X}_2-P)]\,[(\mathbb{X}_1-P)(\mathbb{X}_2-P)]\ \mid \mathbb{X}_2\right)\right) \\ &= \mathbb{E}\left([P(\mathbb{X}_2-P)P+P\,\mathrm{Tr}((\mathbb{X}_2-P)P)]\,(\mathbb{X}_2-P)\right) \end{aligned}$$

This yields the formula

$$\begin{aligned} & M_{\pi^3}(P) \\ & = P\ \Gamma(P)+P\ \Omega(P)=P\ \ [\Gamma(P)+\Omega(P)]=3\ P^4+\mathrm{Tr}(P^2)\,P^2 \end{aligned} \tag{6.12}$$

We let $M_{2n,n}^{[Q_0,\ldots,Q_{2n}]}(P)$, and respectively $M_{2n,n}^{[Q_0,\ldots,Q_{2n}],-}(P)$, the matrix moments defined as in (6.9) by replacing the set $\mathcal{Q}_{2n,n}^{+}$ by the set of crossing partitions $\mathcal{Q}_{2n,n}$, and respectively $\mathcal{Q}_{2n,n}^{-}$.

Arguing as above, for any matrices $Q_1,Q_2,Q_3$ we have

$$\begin{aligned} & \mathbb{E}\left((\mathbb{X}_1-P)Q_1(\mathbb{X}_2-P)Q_2(\mathbb{X}_1-P)Q_3(\mathbb{X}_2-P)\right) \\ & = \mathbb{E}\left((\mathbb{X}_1-P)Q_1\mathbb{E}\left[(\mathbb{X}_2-P)Q_2(\mathbb{X}_1-P)Q_3(\mathbb{X}_2-P)\ \mid \mathbb{X}_1\right]\right) \\ & = \mathbb{E}\left((\mathbb{X}_1-P)Q_1\ \Gamma\left(Q_2(\mathbb{X}_1-P)Q_3\right)\right) \\ & = \mathbb{E}\left((\mathbb{X}_1-P)\ Q_1PQ_2\ (\mathbb{X}_1-P)\right)\ Q_3P \\ & \qquad\qquad +\mathbb{E}\left(\mathrm{Tr}\left((\mathbb{X}_1-P)(Q_3PQ_2)\right)(\mathbb{X}_1-P)\right)\ Q_1P \end{aligned}$$

This yields the formula

$$\begin{aligned} & \mathbb{E}\left(Q_0(\mathbb{X}_1-P)Q_1(\mathbb{X}_2-P)Q_2(\mathbb{X}_1-P)Q_3(\mathbb{X}_2-P)Q_4\right) \\ & = \mathbb{E}\left[Q_0PQ_1PQ_2PQ_3PQ_4\right]+2\ \mathbb{E}\left[Q_0PQ_3PQ_2PQ_1PQ_4\right] \\ & \qquad\qquad +\mathbb{E}\left[\mathrm{Tr}(Q_1PQ_2P)\ Q_0PQ_3PQ_4\right] \end{aligned} \tag{6.13}$$

for any random matrices $Q_i$, with $0 \leq i \leq 4$, independent of $\mathbb{X}_1$ and $\mathbb{X}_2$. Any $\pi$-matrix moment $M_\pi(P)$ associated with a crossing partition $\pi \in \mathcal{Q}_{2n,n} - \mathcal{Q}^+_{2n,n}$ has the form (6.13) with some random matrix $Q_0$ associated with a non-crossing partition $\pi^{[m]} \in \mathcal{Q}^+_{2m,m}$ for some $0 \leq m \leq n-2$. When $m = 0$ we use the convention $Q_0 = I$.

Observe that

$$\mathbb{E}\left[Q_0 P Q_1 P Q_2 P Q_3 P Q_4\right] \quad \text{and} \quad \mathbb{E}\left[Q_0 P Q_3 P Q_2 P Q_1 P Q_4\right]$$

can be interpreted as a matrix moment

$$\begin{aligned} & M_\pi^{[R_0,\ldots,R_{2((n-m)-2)}]}(P) \\ & = \mathbb{E}\Big(R_0(\mathbb{X}_{\alpha^\pi(1)} - P)R_1(\mathbb{X}_{\alpha^\pi(2)} - P)R_2 \ldots \\ & \qquad \ldots (\mathbb{X}_{\alpha^\pi(2((n-m)-2))} - P)R_{2((n-m)-2)}\Big) \end{aligned}$$

associated with some deterministic matrices $R_i \in \{P^k \,:\, 0 \leq k \leq 4\}$ and some partition $\pi \in \mathcal{Q}_{2((n-m)-2),((n-m)-2)}$. By (6.11) when $\pi \in \mathcal{Q}^+_{2((n-m)-2),(n-m)-2}$, the matrix moment $M_\pi^{[R_0,\ldots,R_{2((n-m)-2)}]}(P)$ is a linear combination of monomials of the form

$$\mathrm{Tr}_v(P)\ P^w$$

with the parameters

$$\begin{cases} 1 \leq w \leq 2(n-m) \qquad 0 \leq m \leq (n-m)-2 \\ \sum_{1\leq i\leq m} v_i \leq (n-m)-2 \\ \sum_{1\leq i\leq m}\ i v_i + w = 2(n-m) \end{cases}$$

Assume that $(Q_1, Q_2)$ and $(Q_3, Q_4)$ are independent. In this case we have

$$\mathbb{E}\left[\mathrm{Tr}(Q_1 P Q_2 P)\ Q_0 P Q_3 P Q_4\right] = \mathrm{Tr}(\mathbb{E}\left[Q_1 P Q_2 P\right])\ \mathbb{E}\left[Q_0 P Q_3 P Q_4\right]$$

Arguing as above

$$\mathbb{E}\left[Q_1 P Q_2 P\right] \quad \text{and} \quad \mathbb{E}\left[Q_0 P Q_3 P Q_4\right]$$

can be can be interpreted as the matrix moments $M_\alpha^{[S_0,\ldots,S_k]}(P)$ and $M_\beta^{[T_0,\ldots,T_l]}(P)$ associated with a collection of deterministic matrices $S_i \in \{P^k \;:\; 0 \leq k \leq 2\}$ and $T_j \in \{P^k \;:\; 0 \leq k \leq 2\}$ and some partitions $\alpha \in \mathcal{Q}_{2k,k}$ and $\beta \in \mathcal{Q}_{2l,l}$ with $k + l = (n-m) - 2$.

By (6.11) when $\alpha \in \mathcal{Q}^+_{2k,k}$, the matrix moment $M_\alpha^{[S_0,\ldots,S_k]}(P)$ is a linear combination of monomials of the form

$$\mathrm{Tr}_v(P)\; P^w$$

with the parameters

$$\begin{cases} 1 \leq w \leq 2(k+1) \qquad 0 \leq m \leq k \qquad \sum_{1\leq i\leq m} v_i \leq k \\ \\ \sum_{1\leq i\leq m}\; iv_i + w \leq 2(k+1) \end{cases}$$

In addition, when $\beta \in \mathcal{Q}^+_{2k,k}$, the matrix moment product

$$M_\alpha^{[S_0,\ldots,S_k]}(P)\; M_\beta^{[T_0,\ldots,T_l]}(P)$$

is a linear combination of monomials of the form

$$\mathrm{Tr}_v(P)\; P^w \quad \text{with} \quad \begin{cases} 1 \leq w \leq 2(n-m) \qquad 0 \leq m \leq k \\ \\ \sum_{1\leq i\leq m} v_i \leq n-m-2 \\ \\ \sum_{1\leq i\leq m}\; iv_i + w = 2(n-m) \end{cases}$$

Whenever $(Q_1, Q_2)$ and $(Q_3, Q_4)$ are not independent we have

$$\mathrm{Tr}(Q_1PQ_2P)\; Q_0PQ_3PQ_4 \stackrel{law}{=} \mathrm{Tr}(\overline{Q}_1(\mathbb{X}-P)\overline{Q}_2)\; Q_0\overline{Q}_3(\mathbb{X}-P)\overline{Q}_4$$

for some matrices $Q_0$ and $\overline{Q}_i$ independent of $\mathbb{X}$. In this situation, taking the expectation and using (6.6) we check that

$$\mathbb{E}\left[\mathrm{Tr}(Q_1PQ_2P)\; Q_0PQ_3PQ_4\right] = 2\;\mathbb{E}\left[Q_0\overline{Q}_3P\overline{Q}_2\overline{Q}_1P\overline{Q}_4\right]$$

Here again $\mathbb{E}\left[Q_0\overline{Q}_3P\overline{Q}_2\overline{Q}_1P\overline{Q}_4\right]$ can be interpreted as the matrix moments $M_\alpha^{[U_0,\ldots,U_v]}(P)$ associated with a collection of deterministic matrices $U_i$ and some $\alpha \in \mathcal{Q}_{2((n-m)-3),((n-m)-3)}$. By (6.11) when $\alpha \in \mathcal{Q}^+_{2((n-m)-3),((n-m)-3)}$, the matrix moment $M_\alpha^{[U_0,\ldots,U_v]}(P)$ is a linear combination of monomials of the form

$$\mathrm{Tr}_v(P)\ P^w$$

with the parameters

$$\begin{cases} 1 \leq w \leq 2((n-m)-1) \qquad 0 \leq m \leq (n-m)-3 \\ \\ \sum_{1\leq i\leq m} v_i \leq (n-m)-3 \\ \\ \sum_{1\leq i\leq m}\ iv_i + w = 2((n-m)-1) \end{cases}$$

In summary, we have proved that the $\pi$-matrix moment $M_\pi(P)$ associated with a given crossing partition $\pi \in \mathcal{Q}^-_{2n,n}$ can be expressed as a sum of matrix moments products

$$M_\alpha^{[R_0,\ldots,R_{2k}]}(P)\ M_\beta^{[S_0,\ldots,S_{2l}]}(P)$$

associated with some partitions $(\alpha,\beta) \in (Q_{2k,k} \times \mathcal{Q}_{2l,l})$ with $0 \leq k+l \leq n-2$, and some matrices $R_i, S_i \in \{P^k\ :\ 0 \leq k \leq 4\}$ whose values depend on $\pi$. Iterating this procedure whenever $\alpha$ and/or $\beta$ are crossing, we can write $M_\pi(P)$ as a sum of matrix moments products

$$M_{\gamma_1}^{[T_0^{(1)},\ldots,T_{2k_1}^{(1)}]}(P)\ \ldots\ M_{\gamma_q}^{[T_0^{(m)},\ldots,T_{2k_q}^{(1)}]}(P)$$

for some $2 \leq q \leq n-2$ and some *non-crossing* partitions $\gamma_i \in Q^+_{2k_i,k_i}$ with $0 \leq k_1 + \ldots + k_q \leq n-2$. In addition, $M_{2k_i,k_i}^{[T_0^{(1)},\ldots,T_{2k_i}^{(1)}]}(P)$ is a linear combination of monomials of the form

$$\mathrm{Tr}_v(P)\ P^w$$

with the parameters

$$\begin{cases} 1 \leq w \leq l_i \qquad 0 \leq m \leq k_i \qquad \sum_{1\leq j\leq m} v_j \leq k_i \\ \\ \sum_{1\leq j\leq m}\ jv_j + w \leq l_i \end{cases}$$

for some $l_i \geq 2k_i$. We conclude that $M^-_{2n,n}(P)$ is a linear combination of monomials of the form

$$\mathrm{Tr}_v(P)\ P^w$$

with the parameters

$$\begin{cases} 1 \leq w \leq 2n \qquad 0 \leq m \leq n-2 \qquad \sum_{1\leq j\leq m} v_j \leq n-2 \\ \\ \sum_{1\leq j\leq m} \; jv_j + w \leq 2n \end{cases}$$

The quantity $2n$ corresponds to maximal power of $P$ we can obtain using crossing partitions. It arises when we consider partitions $\pi$ with the maximal number of pairwise crossing pairs, non-crossing between them. For instance for an even number $n$ of pairs, pairwise crossing and non-overlapping we have

$$(6.12) \implies M_\pi(P) = \left(3 \; P^4 + \mathrm{Tr}(P^2)\, P^2\right)^{n/2}$$

The proof of theorem 3.4 is now complete. □

## 6.3 Proof of theorem 3.5

We check that the matrices $\Sigma_n(P)$ defined in (1.14) satisfy the recursion (3.15). Observe that

$$\Sigma_n(P)P = \mathrm{Tr}(P) \sum_{\pi\in\mathcal{N}_n} \left[\prod_{i\geq 0} \mathrm{Tr}\left(P^{1+i}\right)^{r_i(\pi)}\right] \frac{P^{|\pi_1|+1}}{\mathrm{Tr}(P^{|\pi_1|+1})}$$

This yields the formula

$$\mathrm{Tr}\left(\Sigma_n(P)P\right) = \mathrm{Tr}(P) \sum_{\pi\in\mathcal{N}_n} \left[\prod_{i\geq 0} \mathrm{Tr}\left(P^{1+i}\right)^{r_i(\pi)}\right]$$

Thus, for any $k + l = n$, we have

$$\mathrm{Tr}\left[\Sigma_k(P)P\right] \; \left[\Sigma_l(P)P\right]$$

$$= \mathrm{Tr}(P)^2 \sum_{(\pi,\overline{\pi})\in(\mathcal{N}_l\times\mathcal{N}_k)} \left[\prod_{i\geq 0} \mathrm{Tr}\left(P^{1+i}\right)^{r_i(\pi)+r_i(\overline{\pi})}\right] \frac{P^{|\pi_1|+1}}{\mathrm{Tr}(P^{|\pi_1|+1})}$$

We associate with a pair of partitions $(\pi, \overline{\pi}) \in (\mathcal{N}_l \times \mathcal{N}_k)$ a partition $\tau \in \mathcal{N}_{n+1}$ with blocks

$$\begin{aligned} \tau_1 := \pi_1 \cup \{n+1\} \leq \tau_i := \pi_i \leq \tau_{|\pi|+j} &:= (k+1) + \overline{\pi}_j \\ &= \{(k+1) + u : u \in \overline{\tau}_j\} \end{aligned}$$

for any $2 \leq i \leq |\pi|$ and $1 \leq j \leq |\overline{\pi}|$. Observe that

$$\begin{aligned}
r_{|\pi_1|+1}(\tau) &= r_{|\pi_1|+1}(\pi) + r_{|\pi_1|+1}(\overline{\pi}) + 1 \\
r_{|\pi_1|}(\tau) &= r_{|\pi_1|}(\pi) + r_{|\pi_1|}(\overline{\pi}) - 1 \\
r_i(\tau) &= r_i(\pi) + r_{|\pi_1|}(i) \qquad \forall i \geq 1 \quad \text{and} \quad i \notin \{|\pi_1|, |\pi_1|+1\}
\end{aligned}$$

$$\Longrightarrow \sum_{i\geq 1} r_i(\tau) = \sum_{i \geq 1} [r_i(\pi) + r_i(\overline{\pi})]$$

On the other hand, we have

$$\mathrm{Tr}\,(P)^{n+1-\sum_{i\geq 1} r_i(\pi)+r_i(\overline{\pi})} \left[\prod_{i\geq 1} \mathrm{Tr}\left(P^{1+i}\right)^{r_i(\pi)+r_i(\overline{\pi})}\right] \frac{P^{|\pi_1|+1}}{\mathrm{Tr}(P^{|\pi_1|+1})}$$

$$= \mathrm{Tr}\,(P)^{r_0(\tau)} \left[\prod_{i\geq 1,\ i\notin\{|\pi_1|,|\pi_1|+1\}} \mathrm{Tr}\left(P^{1+i}\right)^{r_i(\tau)}\right]$$

$$\times\ \mathrm{Tr}\left(P^{|\tau_1|}\right)^{r_{|\tau_1|-1}(\tau)+1} \ \mathrm{Tr}\left(P^{1+|\tau_1|}\right)^{r_{|\tau_1|}(\tau)-1} \frac{P^{|\tau_1|}}{\mathrm{Tr}(P^{|\tau_1|})}$$

This implies that

$$\sum_{k+l=n} \mathrm{Tr}\,[\Sigma_k(P)P] \ \ [\Sigma_l(P)P]$$

$$= \mathrm{Tr}(P) \sum_{\tau\in\mathcal{N}_{n+1}} \mathrm{Tr}\,(P)^{r_0(\tau)} \left[\prod_{i\geq 1} \mathrm{Tr}\left(P^{1+i}\right)^{r_i(\tau)}\right] \frac{P^{|\tau_1|}}{\mathrm{Tr}(P^{|\tau_1|+1})}$$

$$= \Sigma_{n+1}(P)$$

This ends the proof of the matrix polynomial formula (1.14).

Now we come to the proof of the Kreweras trace-type formula (3.17). For any $\tau \in \mathcal{N}_{n+1}$ we set

$$\begin{aligned}
&k_{|\tau_1|-1}(\tau) := r_{|\tau_1|-1}(\tau) + 1 \qquad k_1(\tau) = 1 + r_0(\tau) \\
&k_{|\tau_1|+1}(\tau) := r_{|\tau_1|}(\tau) - 1 \qquad \text{and} \qquad k_{i+1}(\tau) := r_i(\tau)
\end{aligned}$$

for any $i \geq 0$ s.t. $i \notin \{0, |\tau_1| - 1, |\tau_1|\}$.

Observe that

$$\sum_{i\geq 1} k_i(\tau) = \; 1 + r_0(\tau) + \sum_{i\geq 1} r_i(\tau) = 1 + n$$

We also have

$$\begin{aligned}
&\sum_{i\geq 1} \; i \; k_i(\tau) \\
&= 1 + r_0(\tau) + \sum_{j\geq 1,\; j\notin\{|\tau_1|-1,|\tau_1|\}} (j+1)\; r_j(\tau) \\
&\qquad + |\tau_1| \; \left(r_{|\tau_1|-1} + 1\right) + (|\tau_1| + 1) \; \left(r_{|\tau_1|} - 1\right) \\
&= \; r_0(\tau) + \textstyle\sum_{j\geq 1} \; r_j(\tau) + \sum_{j\geq 1} \; j \; r_j(\tau) = 2n
\end{aligned}$$

Taking the trace in the above matrix polynomial formula we find that

$$\begin{aligned}
\mathrm{Tr}(\Sigma_n(P)) &= \sum_{\tau\in\mathcal{N}_n} \; \mathrm{Tr}\,(P)^{k_1(\tau)} \; \left[ \prod_{i\geq 1,\; i\notin\{|\tau_1|-1,|\tau_1|\}} \mathrm{Tr}\left(P^{1+i}\right)^{k_{i+1}(\tau)} \right] \\
&\qquad \times \; \mathrm{Tr}\left(P^{1+|\tau_1|}\right)^{r_{|\tau_1|}(\tau)-1} \; \mathrm{Tr}\left(P^{|\tau_1|}\right)^{r_{|\tau_1|-1}(\tau)+1} \\
&= \sum_{\tau\in\mathcal{N}_n} \left[ \prod_{i\geq 1} \mathrm{Tr}\left(P^i\right)^{k_i(\tau)} \right]
\end{aligned}$$

This ends the proof of (3.17). The proof of the theorem is now completed.

□

## 6.4 Proof of theorem 3.7

Observe that $\mathrm{cl}(\mathcal{N}_{n,m})$ coincides with the subset of partitions $\pi \in \mathcal{N}_{n+1,m}$ s.t. $\iota(\pi) = 1$. Also observe that for any given $\pi \in \mathcal{N}_{n,m}$ the partition $\Xi(\mathrm{cl}(\pi))$ has the same combinatorial structure as $\Xi(\pi)$ but it contains one more singleton than $\Xi(\pi)$. This yields the formula

$$\sum_{\pi\in\mathcal{N}_{n,m}}\left[\prod_{i\geq 1}\mathrm{Tr}(P^i)^{r_i(\Xi(\pi))}\right]P$$

$$=\sum_{\pi\in\mathcal{N}_{n,m}}\left[\prod_{i\geq 1}\mathrm{Tr}(P^i)^{r_i(\Xi(\mathrm{cl}(\pi)))}\right]\frac{P^{\iota(\mathrm{cl}(\pi))}}{\mathrm{Tr}(P^{\iota(\mathrm{cl}(\pi))})}$$

from which we prove that

$$\sum_{\pi_1\in\mathcal{N}_{n_1,m_1}}\sum_{\pi_2\in\mathcal{N}_{n_2,m_2}}\left[\prod_{i\geq 1}\mathrm{Tr}(P^i)^{r_i(\Xi(\pi_1))+r_i(\Xi(\pi_2))}\right]\frac{P^{\iota(\pi_1)}}{\mathrm{Tr}\left(P^{\iota(\pi_1)}\right)}P$$

$$=\sum_{\pi_1\in\mathcal{N}_{n_1,m_1}}\sum_{\pi_2\in\mathcal{N}_{n_2,m_2}}$$

$$\left[\prod_{i\geq 1}\mathrm{Tr}(P^i)^{r_i(\Xi(\pi_1)\oplus\Xi(\mathrm{cl}(\pi_2)))}\right]\frac{P^{\iota(\pi_1)+\iota(\mathrm{cl}(\pi_2))}}{\mathrm{Tr}\left(P^{\iota(\pi_1)}\right)\mathrm{Tr}\left(P^{\iota(\mathrm{cl}(\pi_2))}\right)}$$

Following the arguments and construction in Section 2.5 we have

$$\forall i\notin\{\iota(\pi_1),\iota(\mathrm{cl}(\pi_2))\}\qquad r_i(\Xi(\pi_1)\oplus\Xi(\mathrm{cl}(\pi_2)))=r_i(\Xi(\pi_1\oplus\mathrm{cl}(\pi_2)))$$

as well as

$$\forall i\in\{\iota(\pi_1),\iota(\mathrm{cl}(\pi_2))\}\qquad r_i(\Xi(\pi_1)\oplus\Xi(\mathrm{cl}(\pi_2)))-1=r_i(\Xi(\pi_1\oplus\mathrm{cl}(\pi_2)))$$

and

$$r_{\iota(\pi_1)+\iota(\mathrm{cl}(\pi_2))}(\Xi(\pi_1\oplus\mathrm{cl}(\pi_2)))=r_{\iota(\pi_1)+\iota(\mathrm{cl}(\pi_2))}(\Xi(\pi_1)\oplus\Xi(\mathrm{cl}(\pi_2)))+1$$

This implies that

$$\sum_{\pi_1\in\mathcal{N}_{n_1,m_1}}\sum_{\pi_2\in\mathcal{N}_{n_2,m_2}}\left[\prod_{i\geq 1}\mathrm{Tr}(P^i)^{r_i(\Xi(\pi_1))+r_i(\Xi(\pi_2))}\right]\frac{P^{\iota(\pi_1)}}{\mathrm{Tr}\left(P^{\iota(\pi_1)}\right)}P$$

$$=\sum_{\pi_1\in\mathcal{N}_{n_1,m_1}}\sum_{\pi_2\in\mathcal{N}_{n_2,m_2}}\left[\prod_{i\geq 1}\mathrm{Tr}(P^i)^{r_i(\Xi(\pi_1\oplus\mathrm{cl}(\pi_2)))}\right]\frac{P^{\iota(\pi_1\oplus\mathrm{cl}(\pi_2))}}{\mathrm{Tr}\left(P^{\iota(\pi_1\oplus\mathrm{cl}(\pi_2))}\right)}$$

Using (2.21), this shows that (1.15) satisfies (3.19). □

## 6.5 Proof of theorem 4.1

Combining (6.3) with the estimates (3.2) we check that

$$\|\mathbb{E}\left(H_N^n\right)\|_F \leq n!\ (8\mathrm{Tr}(P))^n\ N^{-n/2} \sum_{1\leq m\leq \lfloor n/2\rfloor\wedge N} \frac{(N)_m}{m!}\ 2^{-(2m+1)}$$

This ends the proof of (4.2). Also observe that for any $n \leq N$ we have

$$\mathbb{E}\left[\mathcal{H}_N^{2n}\right] = \frac{(N)_n}{N^n}\ M_{2n,n}(P) + \sum_{1\leq m<n} \frac{(N)_m}{N^n}\ M_{2n,m}(P) \tag{6.14}$$

We check (4.3) using the estimates (3.2). We also have

$$1 - \frac{(n-1)^2}{N} \leq \frac{(N)_n}{N^n} \leq 1$$

Using the estimates (3.2) we readily check (4.4). The proof is now completed. □

## 6.6 Proof of theorem 4.2

Using the fact that

$$[I-2sP]^{-1}P\ = \partial_s \log\left([I-2sP]^{-1/2}\right) = -\frac{1}{2}\ \partial_s \log\left(I-2sP\right)$$

formula (4.5) can alternatively be rewritten as

$$\begin{aligned}\mathbb{E}\left(\exp\left[t\,\mathbb{X}\right]\right)\ &= I + \int_0^t \exp\left[\mathrm{Tr}\left(\log\left[(I-2sP)^{-1/2}\right]\right)\right] \\ &\qquad\qquad \times \partial_s \log\left[(I-2sP)^{-1/2}\right]\ ds \\ &\leq I - \frac{1}{2}\log\left(I-2tP\right)\ \exp\left[-\frac{1}{2}\ \mathrm{Tr}\left(\log\left(I-2tP\right)\right)\right]\end{aligned}$$

The proof of theorem 4.2 is based on the following technical lemma.

**Lemma 6.2.** For any $0 \leq t \leq 1/(2\mathrm{Tr}(P))$ we have

$$\begin{aligned}&\mathbb{E}\left(\exp\left[\frac{t}{2}\frac{\mathbb{X}}{\mathrm{Tr}(P)}\right]\right) \\ &\leq I + \frac{t}{2}\ \frac{P}{\mathrm{Tr}(P)} + \frac{t^2}{4}\ \left[\left(\frac{P}{\mathrm{Tr}(P)}\right)^2 + \frac{1}{2}\ \frac{P}{\mathrm{Tr}(P)}\right] + t\,L(t)\ \varpi_0\left(\frac{P}{\mathrm{Tr}(P)}\right)\end{aligned} \tag{6.15}$$

with

$$\varpi_0(Q) := Q \left[\frac{1}{6\sqrt{7}}\; I \;+ \frac{1}{12}\; Q + \; Q^2\right]$$

In addition, we have

$$\begin{aligned}
&\mathbb{E}\left(\exp\left[\frac{t}{2}\frac{\mathbb{X}-P}{\mathrm{Tr}(P)}\right]\right) \\
&\leq I + \frac{t^2}{8}\;\left[\left(\frac{P}{\mathrm{Tr}(P)}\right)^2 + \frac{P}{\mathrm{Tr}(P)}\right] + t\,L(t)\;\varpi_0\left(\frac{P}{\mathrm{Tr}(P)}\right)
\end{aligned} \tag{6.16}$$

*Proof.* Using (3.12), for any $t \geq 0$ we check that

$$\exp[tP] \leq \mathbb{E}\left(\exp[t\mathbb{X}]\right) \leq I + Pt + \frac{t^2}{2}\;\left[2P^2 + P\mathrm{Tr}(P)\right] + \sum_{n\geq 3}\frac{t^n}{n!}\;\mathbb{E}\left(\mathbb{X}^n\right)$$

as well as

$$\mathbb{E}\left(\exp[t(\mathbb{X}-P)]\right) \leq I + \frac{t^2}{2}\;\left[P^2 + P\,\mathrm{Tr}(P)\right] + \sum_{n\geq 3}\frac{t^n}{n!}\;\mathbb{E}\left(\mathbb{X}^n\right)$$

On the other hand, using (A.1) we also check that

$$\begin{aligned}
&\sum_{n\geq 3}\frac{t^n}{n!}\;\mathbb{E}\left(\mathbb{X}^n\right) \leq \sum_{n\geq 3}\frac{(2t)^n}{2n} \\
&\times\left[\frac{(2(n-1))!}{2^{2(n-1)}\;(n-1)!^2}\;\mathrm{Tr}(P)^{n-1}\;\; P + \frac{(2(n-2))!}{2^{2(n-2)}\;(n-2)!^2}\;\mathrm{Tr}(P)^{n-2}\;\; P^2\right] \\
&\qquad\qquad + (2tP)^3\sum_{n\geq 0}\frac{(2t)^n}{2(n+3)}\;\sum_{0\leq k\leq n}\frac{(2k)!}{2^{2k}\;k!^2}\;\mathrm{Tr}(P)^k\;\; P^{n-k}
\end{aligned}$$

By (A.2) and (A.3) this yields the estimate

$$\begin{aligned}
\sum_{n\geq 3}\frac{t^n}{n!}\;\mathbb{E}\left(\mathbb{X}^n\right)\; \leq\; &\frac{(2t)^3}{2}\;\sum_{n\geq 0}\frac{(2t\mathrm{Tr}(P))^n}{n+3}\;\frac{1}{\sqrt{3n+7}}\;P\,\mathrm{Tr}(P)^2 \\
&+\frac{(2t)^3}{2}\;\sum_{n\geq 0}\frac{(2t\mathrm{Tr}(P))^n}{n+3}\;\frac{1}{\sqrt{3n+4}}\;\;P^2\mathrm{Tr}(P) \\
&+(2tP)^3\sum_{n\geq 0}\frac{(2t\mathrm{Tr}(P))^n}{2(n+3)}\;\frac{2n+1}{\sqrt{3n+1}}
\end{aligned} \tag{6.17}$$

This implies that

$$\sum_{n\geq 3} \frac{t^n}{n!}\,\mathbb{E}\left(\mathbb{X}^n\right) \leq \frac{(2t)^2}{2} \sum_{n\geq 1} \frac{(2t\mathrm{Tr}(P))^n}{n+2}\,\frac{1}{\sqrt{3n+4}}\,P\,\mathrm{Tr}(P)$$
$$+\frac{(2t)^2}{2} \sum_{n\geq 1} \frac{(2t\mathrm{Tr}(P))^n}{n+2}\,\frac{1}{\sqrt{3n+1}}\,P^2+\frac{(2tP)^3}{1-2t\mathrm{Tr}(P)}$$

from which we check that

$$\begin{aligned}
&\sum_{n\geq 3} \frac{t^n}{n!}\,\mathbb{E}\left(\mathbb{X}^n\right)\\
&\leq\; -\left(2t\mathrm{Tr}(P)\right)^2\,\left(\frac{1}{4}\frac{P}{\mathrm{Tr}(P)}+\frac{1}{2}\left(\frac{P}{\mathrm{Tr}(P)}\right)^2\right)\log\left(1-2t\mathrm{Tr}(P)\right)\\
&\qquad\qquad+\frac{(2tP)^3}{1-2t\mathrm{Tr}(P)}\\
&\leq \frac{(2t\mathrm{Tr}(P))^3}{1-2t\mathrm{Tr}(P)}\,\left[\frac{1}{4}\,\frac{P}{\mathrm{Tr}(P)}+\frac{1}{2}\left(\frac{P}{\mathrm{Tr}(P)}\right)^2+\left(\frac{P}{\mathrm{Tr}(P)}\right)^3\right]
\end{aligned}$$

In the last assertion we have used the estimate $-\log\left(1-x\right)\leq x/(1-x)$. This estimate can be slightly improved. Indeed, using (6.17) we check that

$$\begin{aligned}
&\sum_{n\geq 3} \frac{t^n}{n!}\,\mathbb{E}\left(\mathbb{X}^n\right)\\
&\leq \frac{(2t\mathrm{Tr}(P))^3}{1-2t\mathrm{Tr}(P)}\,\left[\frac{1}{6\sqrt{7}}\frac{P}{\mathrm{Tr}(P)}+\frac{1}{12}\,\left(\frac{P}{\mathrm{Tr}(P)}\right)^2+\left(\frac{P}{\mathrm{Tr}(P)}\right)^3\right]
\end{aligned}$$

This ends the proof of the lemma follows elementary manipulations, thus it is skipped. □

We now return to prove theorem 4.2. The estimate (6.15) can be rewritten as follows

$$\begin{aligned}
&\mathbb{E}\left(\exp\left[\frac{t}{2}\frac{\mathbb{X}}{\mathrm{Tr}(P)}\right]\right)\\
&\leq I+\frac{t}{2}\,\frac{P}{\mathrm{Tr}(P)}+t^2\,\left[\frac{1}{4}\,\left(\frac{P}{\mathrm{Tr}(P)}\right)^2+\frac{1}{8}\,\frac{P}{\mathrm{Tr}(P)}\right]\\
&\quad+\left(\frac{1}{1-t}-1\right)\,t^2\,\left[\frac{1}{6\sqrt{7}}\,\frac{P}{\mathrm{Tr}(P)}\,+\frac{1}{12}\,\left(\frac{P}{\mathrm{Tr}(P)}\right)^2+\left(\frac{P}{\mathrm{Tr}(P)}\right)^3\right]\\
&\leq I+\frac{t}{2}\,\frac{P}{\mathrm{Tr}(P)}+\frac{t^2}{1-t}\,\left[\frac{1}{4}\,\left(\frac{P}{\mathrm{Tr}(P)}\right)^2+\frac{1}{8}\,\frac{P}{\mathrm{Tr}(P)}\right]\\
&\qquad\qquad+\frac{t^2}{1-t}\,\left[\frac{1}{6\sqrt{7}}\,\frac{P}{\mathrm{Tr}(P)}\,+\frac{1}{12}\,\left(\frac{P}{\mathrm{Tr}(P)}\right)^2+\left(\frac{P}{\mathrm{Tr}(P)}\right)^3\right]\\
&\leq I+\frac{t}{2}\,\frac{P}{\mathrm{Tr}(P)}+\frac{t^2}{1-t}\,\left[\frac{1}{3}\,\left(\frac{P}{\mathrm{Tr}(P)}\right)^2+\frac{1}{5}\,\frac{P}{\mathrm{Tr}(P)}+\left(\frac{P}{\mathrm{Tr}(P)}\right)^3\right]
\end{aligned}$$

This yields

$$\mathbb{E}\left(\exp\left[\frac{t}{2}\,\frac{\mathbb{X}}{\mathrm{Tr}(P)}\right]\right)\leq I+\frac{t}{2}\,\frac{P}{\mathrm{Tr}(P)}+\frac{t^2}{1-t}\,\varpi_1\left(\frac{P}{\mathrm{Tr}(P)}\right) \tag{6.18}$$

The exponential matrix estimate (4.7) is now a consequence of (3.12) and the estimate (6.18). In the same vein, we have

$$\begin{aligned}
&\mathbb{E}\left(\exp\left[t\left(\mathbb{X}-I\right)\right]\right)\\
&\leq I+\frac{(2t\mathrm{Tr}(P))^2}{1-2t\mathrm{Tr}(P)}\,\left[\frac{1}{5}\,\frac{P}{\mathrm{Tr}(P)}\,+\frac{5}{24}\,\left(\frac{P}{\mathrm{Tr}(P)}\right)^2+\left(\frac{P}{\mathrm{Tr}(P)}\right)^3\right]
\end{aligned}$$

This ends the proof of (4.8). The proof of the theorem is completed. □

## 6.7 Proof of theorem 4.4

The proof of theorem 4.4 is partly based on the following technical lemma of its own interest.

**Lemma 6.3.** For any $t \in \mathbb{R}$ we have

$$\log \|\mathcal{E}(t)\|_F \leq \frac{t^2}{2}\ \mathbb{E}(\|\mathcal{H}\|_F^2) \tag{6.19}$$

In addition, for any $m \geq 1$ we have the estimate

$$\left\|\mathcal{E}(t) - \sum_{n \leq m} \frac{t^{2n}}{(2n)!}\ \mathbb{E}(\mathcal{H}^{2n})\right\|_F$$

$$\leq \exp\left[\frac{t^2}{2}\ \mathbb{E}(\|\mathcal{H}\|_F^2)\right] \frac{1}{\sqrt{2\pi(m+1)}} \left(\frac{et^2}{2(m+1)}\ \mathbb{E}(\|\mathcal{H}\|_F^2)\right)^{m+1} \tag{6.20}$$

*Proof.* By (3.2) we have

$$\frac{t^{2n}}{(2n)!}\ \|\mathbb{E}(\mathcal{H}^{2n})\|_F \leq \frac{t^{2n}}{(2n)!}\ \upsilon_{2n}\ \mathbb{E}\left[\|\mathcal{H}\|_F^2\right]^n = \frac{1}{n!}\left(\frac{t^2}{2}\ \mathbb{E}\left[\|\mathcal{H}\|_F^2\right]\right)^n$$

This ends the proof of (6.19). On the other hand, for any $x \geq 0$ we have

$$\sum_{n \geq m} \frac{x^n}{n!} \leq e^x\ \frac{x^m}{m!}\ \left(\leq e^x\ \frac{1}{\sqrt{2\pi m}}\ \left(\frac{ex}{m}\right)^m\right)$$

This implies that

$$\sum_{n > m} \frac{t^{2n}}{(2n)!}\ \|\mathbb{E}(\mathcal{H}^{2n})\|_F$$

$$\leq \exp\left(\frac{t^2}{2}\ \mathbb{E}(\|\mathcal{H}\|_F^2)\right) \frac{1}{\sqrt{2\pi(m+1)}} \left(\frac{et^2}{2m}\ \mathbb{E}(\|\mathcal{H}\|_F^2)\right)^{m+1}$$

This ends the proof of (6.20). The proof of the lemma is completed. ☐

We have

$$\mathbb{E}\left[\cosh\left(t\,\mathcal{H}_N\right)\right] - \mathbb{E}\left(\exp\left(t\mathcal{H}_n\right)\right)$$

$$= \sum_{2 \leq n \leq N} \frac{t^{2n}}{(2n)!}\ \left[\mathbb{E}(\mathcal{H}_N^{2n}) - \mathbb{E}(\mathcal{H}^{2n})\right]$$

$$+ \sum_{n > N} \frac{t^{2n}}{(2n)!}\ \left(\mathbb{E}\left[\mathcal{H}_N^{2n}\right] - \mathbb{E}\left[\mathcal{H}^{2n}\right]\right)$$

On the other hand, for any $n \leq N$ we have

$$\epsilon_{2n}(N) = \sum_{0\leq m<n} \frac{4^m}{N^m} \frac{n!}{(n-m)!}$$

Using (4.3) and (4.4) we obtain the estimate

$$N \|\mathbb{E}[\cosh(t\mathcal{H}_N)] - \mathbb{E}(\exp(t\mathcal{H}_n))\|_F$$

$$\leq \left(\frac{t^2}{2}\mathbb{E}\left[\|\mathcal{H}\|_F^2\right]\right)^2 \sum_{2\leq n\leq N} \frac{1}{(n-2)!} \frac{(n-1)^2}{n(n-1)} \left(\frac{t^2}{2}\mathbb{E}\left[\|\mathcal{H}\|_F^2\right]\right)^{n-2}$$

$$+\frac{1}{2} \sum_{2\leq n\leq N} (4\sqrt{2}\,t\,\mathrm{Tr}(P))^{2n} \sum_{0\leq m<n} \frac{(4/N)^m}{(n-m)!}$$

$$+N \,\|\sum_{n>N} \frac{t^{2n}}{(2n)!} \mathbb{E}\left[\mathcal{H}^{2n}\right]\|_F$$

$$+ \frac{1}{N} (2t\mathrm{Tr}(P))^2 \left(\sqrt{2}t\mathrm{Tr}(P)\right)^{2N} \frac{1}{1-(2\sqrt{2}t\mathrm{Tr}(P))^2/N}$$

$$\times \sum_{0\leq m<N} \left(\frac{4}{N}\right)^m \frac{1}{(N-m)!}$$

as soon as $\sqrt{N} > 2\sqrt{2}\,t\,\mathrm{Tr}(P)$. We use the decomposition

$$\sum_{2\leq n\leq N} a_n \sum_{0\leq m<n} b_{n,m} = \sum_{2\leq k\leq N} a_k\, b_{k,0} + \sum_{1\leq l<N} \sum_{k=l+1}^{N} a_k\, b_{k,l}$$

which is valid for any collection of numbers $a_n$ and $b_{n,m}$, to check that

$$\sum_{2\leq n\leq N} (4\sqrt{2}\,t\,\mathrm{Tr}(P))^{2n} \sum_{0\leq m<n} \frac{(4/N)^m}{(n-m)!}$$

$$= \sum_{2\leq k\leq N} \frac{(4\sqrt{2}\,t\,\mathrm{Tr}(P))^{2k}}{k!}$$

$$+ \sum_{1\leq l<N} \left(\frac{1}{N}\,(8\sqrt{2}\,t\,\mathrm{Tr}(P))^2\right)^l \sum_{k=1}^{N-l} (4\sqrt{2}\,t\,\mathrm{Tr}(P))^{2k}\,\frac{1}{k!}$$

$$\leq \sum_{2\leq k\leq N} \frac{(4\sqrt{2}\,t\,\mathrm{Tr}(P))^{2k}}{k!}$$

$$+ \frac{4(e-1)}{N}\,(4\sqrt{2}\,t\,\mathrm{Tr}(P))^4\,\frac{1}{1-\frac{4}{N}\,(4\sqrt{2}\,t\,\mathrm{Tr}(P))^2}$$

as soon as $4\sqrt{2}\,t\,\mathrm{Tr}(P) < 1$. In addition, when $N \geq 8$ we have

$$\sum_{2\leq n\leq N} (4\sqrt{2}\,t\,\mathrm{Tr}(P))^{2n} \sum_{0\leq m<n} \frac{(4/N)^m}{(n-m)!}$$

$$\leq (4\sqrt{2}\,t\,\mathrm{Tr}(P))^4 \sum_{k\geq 0} \frac{1}{(k+2)(k+1)}\,\frac{(4\sqrt{2}\,t\,\mathrm{Tr}(P))^{2k}}{k!}$$

$$+ \frac{8(e-1)}{N}\,(4\sqrt{2}\,t\,\mathrm{Tr}(P))^4$$

$$\leq \left[\frac{e}{2} + \frac{8(e-1)}{N}\right]\,(4\sqrt{2}\,t\,\mathrm{Tr}(P))^4$$

Choosing $\left(4\sqrt{2}\, t\mathrm{Tr}(P) \leq\right) 4\sqrt{2}\, t\, \mathbb{E}\left[\|\mathcal{H}\|_F^2\right]^{1/2} < 1$ and $N \geq 8$ we check that

$$N\ \|\mathbb{E}\left[\cosh\left(t\,\mathcal{H}_N\right)\right] - \mathbb{E}\left(\exp\left(t\mathcal{H}_n\right)\right)\|_F$$

$$\leq \left(\frac{t^2}{2}\mathbb{E}\left[\|\mathcal{H}\|_F^2\right]\right)^2\ \exp\left(\frac{t^2}{2}\ \mathbb{E}\left[\|\mathcal{H}\|_F^2\right]\right)$$

$$+\left[\frac{e}{4}+\frac{4(e-1)}{N}\right]\ (4\sqrt{2}\, t\,\mathrm{Tr}(P))^4 + N\ \|\sum_{n>N}\ \frac{t^{2n}}{(2n)!}\ \mathbb{E}\left[\mathcal{H}^{2n}\right]\|_F$$

$$+\frac{(e-1)}{4^3N}\ \left(4\sqrt{2}t\mathrm{Tr}(P)\right)^4\ \left(\sqrt{2}t\mathrm{Tr}(P)\right)^{2(N-1)}$$

This yields the estimate

$$N\ \|\mathbb{E}\left[\cosh\left(t\,\mathcal{H}_N\right)\right] - \mathbb{E}\left(\exp\left(t\mathcal{H}_n\right)\right)\|_F$$

$$\leq \left[\frac{e}{4}+\frac{5(e-1)}{N}\right]\ (4\sqrt{2}\, t\,\mathrm{Tr}(P))^4$$

$$+\left(\frac{t^2}{2}\ \mathbb{E}(\|\mathcal{H}\|_F^2)\right)^2\exp\left(\frac{t^2}{2}\ \mathbb{E}(\|\mathcal{H}\|_F^2)\right)$$

$$\times\ \left[1+\frac{1}{N^{3/2}}\ \frac{e^2}{\sqrt{2\pi}}\ \left(\frac{et^2}{2(N+1)}\ \mathbb{E}(\|\mathcal{H}\|_F^2)\right)^{N-1}\right]$$

$$\leq t^4\ \mathbb{E}(\|\mathcal{H}\|_F^2)^2\left(4^4\ e+\frac{e}{4}+\frac{5(e-1)}{N}\ 4^5\ +\frac{1}{N^{3/2}}\ \frac{e^3}{4\sqrt{2\pi}}\ \right)$$

We also have

$$\mathbb{E}\left[\sinh\left(t\,\mathcal{H}_N\right)\right] = \sum_{n\geq 1}\ \frac{t^{2n+1}}{(2n+1)!}\ \mathbb{E}(\mathcal{H}_N^{2n+1})$$

Using (4.2) we find that

$$\|\mathbb{E}\left[\sinh\left(t\,\mathcal{H}_N\right)\right]\|_F$$

$$\leq \frac{t}{2}\frac{\mathrm{Tr}(P)}{N^{1/2}}\sum_{n\geq 1}(2t^2\mathrm{Tr}(P)^2)^n\sum_{1\leq m\leq n\wedge N}\left(\frac{4}{N}\right)^{n-m}\frac{1}{m!}$$

This implies that

$$\|\mathbb{E}\left[\sinh\left(t\,\mathcal{H}_N\right)\right]\|_F$$

$$\leq t^3\,(e-1)\,\frac{\mathrm{Tr}(P)^3}{N^{1/2}}\,\frac{1}{1-2t^2\mathrm{Tr}(P)^2}\leq 2t^3\,(e-1)\,\frac{\mathrm{Tr}(P)^3}{N^{1/2}}$$

as soon as $N\geq 4$ and $4t^2\mathrm{Tr}(P)^2\leq 1$. This ends the proof of (4.11). This ends the proof of the theorem. □

## 6.8 Proof of theorem 5.1

Observe that

$$\mathbb{E}\left(\exp\left(\frac{t}{\sqrt{2r}}\,\mathrm{Tr}\left(A\mathcal{H}_N\right)\right)\right)=\exp\left[\frac{t^2}{2}\,r^{-1}\mathrm{Tr}((AP)^2)\right.$$

$$\left.+\left(\frac{2}{rN}\right)^{1/2}\sum_{n\geq 3}\frac{t^n}{n}\left(\frac{2}{rN}\right)^{(n-3)/2}r^{-1}\mathrm{Tr}((AP)^n)\right]$$

Also recall that $\mathrm{Tr}\left[(AP)^{2n}\right]\geq 0$ and $\mathrm{Tr}\left[(AP)^{2n+1}\right]\in\mathbb{R}$ for any $n\geq 0$. We also have the estimates

$$\begin{aligned}r^{-1}|\mathrm{Tr}\left[(AP)^n\right]| &\leq r^{-1}|\lambda^\star(AP)|^n+r^{-1}\sum_{2\leq i\leq r}\left[\frac{|\lambda_i(AP)|}{\lambda^\star(AP)}\right]^n\\ &\leq r^{-1}\left[\lambda^\star(AP)\right]^n+1\end{aligned}$$

We conclude that

$$\sup_{r\geq 1}\lambda^\star(AP)<\infty\Longrightarrow\sum_{n\geq 3}\frac{t^n}{n}\sup_{r\geq 1}|\,r^{-1}\,\mathrm{Tr}((AP)^n)|<\infty$$

The end of the proof of the first assertion of the theorem is now completed. Using (2.3), for any $n \geq 2$ we check that

$$\begin{aligned} \left|\frac{1}{n}\left(\frac{2}{rN}\right)^{(n-2)/2} r^{-1}\mathrm{Tr}((AP)^n)\right| &\leq \frac{1}{2}\left(\frac{2}{rN}\right)^{(n-2)/2} \alpha_A(P)\ \beta_A(P)^n \\ &= \rho\ \beta^n \end{aligned}$$

with the parameters

$$\rho = \frac{rN}{4}\ \alpha_A(P) \quad \text{and} \quad \beta = \sqrt{\frac{2}{rN}}\ \beta_A(P)$$

The end of the proof is a direct consequence of lemma 2.1. For instance we check that

$$\begin{aligned} &\left|\mathbb{E}\left(\left(\frac{1}{\sqrt{2r}}\ \mathrm{Tr}\left(A\mathcal{H}_N\right)\right)^{2n}\right) - \frac{(2n)!}{2^n n!}\ \left[r^{-1}\mathrm{Tr}((AP)^2)\right]^n\right| \\ &\leq \frac{(2n)!}{2}\ \beta_A(P)^{2n}\ \left(\frac{8}{rN}\right)^n \sum_{1\leq k<n} \frac{1}{k!}\ \left(\frac{rN}{8}\right)^k\ \alpha_A(P)^k \\ &\leq \frac{(2n)!}{2}\ \beta_A(P)^{2n}\ \alpha_A(P)^{n-1}\ \left(\frac{8}{rN}\right)^n \sum_{1\leq k<n} \frac{1}{k!}\ \left(\frac{rN}{8}\right)^k \\ &\leq \frac{(e-1)}{2}\ (2n)!\ \alpha_A(P)^{n-1}\ \beta_A(P)^{2n}\ \left(\frac{8}{rN}\right)^n \left(1 \vee \frac{rN}{8}\right)^{n-1} \end{aligned}$$

This end the proof of (3.7). □

## 6.9 Proof of theorem 5.2

Using the sub-Gaussian estimate (5.7) we check that

$$\mathbb{P}\left(\pm\mathrm{Tr}(A\mathcal{H}_N) \geq u\right) \leq \exp\left(-\frac{1}{4}\ \frac{u^2}{\mathrm{Tr}((AP)^2) + 2\|AP\|_F^2}\right)$$

for any

$$u \leq \sqrt{N}\ \frac{2\|AP\|_F^2 + \mathrm{Tr}((AP)^2)}{2\|AP\|_F}$$

We end the proof of the first assertion by choosing

$$\delta = \frac{1}{4}\,\frac{u^2}{\mathrm{Tr}((AP)^2) + 2\|AP\|_F^2} \leq \frac{N}{8}\,\left(1 + \frac{\mathrm{Tr}((AP)^2)}{2\|AP\|_F^2}\right) \Longleftarrow \delta \leq \frac{N}{8}$$

By lemma A.2, we have

$$2t\mathrm{Tr}(AP) < \sqrt{N} \Longrightarrow \log \mathbb{E}\left(\exp\left(-t\ \mathrm{Tr}\left(A\mathcal{H}_N\right)\right)\right) \leq \frac{t^2}{2}\ \left(2\mathrm{Tr}((AP)^2)\right)$$

Thus, for any $u \geq 0$ and any $t \geq 0$ s.t. $2\,t\,\mathrm{Tr}(AP) < \sqrt{N}$ we have

$$\mathbb{P}\left(-\mathrm{Tr}(A\mathcal{H}_N) \geq u\right) \leq \mathbb{E}\left[\exp\left(-u\,t - t\,\mathrm{Tr}(A\mathcal{H}_N)\right)\right]$$

This yields

$$\mathbb{P}\left(-\mathrm{Tr}(A\mathcal{H}_N) \geq u\right) \leq \exp\left(-ut + t^2\ \mathrm{Tr}((AP)^2)\right)$$

Choosing $N$ and $t$ such that

$$t = \frac{u}{2\mathrm{Tr}((AP)^2)} \leq \frac{\sqrt{N}}{2\mathrm{Tr}(AP)} \Longleftrightarrow \sqrt{N} \geq u\ \frac{\mathrm{Tr}(AP)}{\mathrm{Tr}((AP)^2)}$$

we find that

$$\mathbb{P}\left(\mathrm{Tr}(A\mathcal{H}_N) \leq -u\right) \leq \exp\left(-\frac{1}{4}\ \frac{u^2}{\mathrm{Tr}((AP)^2)}\right)$$

This implies that

$$0 \leq \delta \leq \frac{N}{4}\ \frac{\mathrm{Tr}((AP)^2)}{\mathrm{Tr}(AP)^2}$$

$$\Longrightarrow \mathbb{P}\left(\mathrm{Tr}(A\mathcal{H}_N) \leq -2\ \sqrt{\delta\ \mathrm{Tr}((AP)^2)}\right) \leq \exp\left(-\delta\right)$$

In the same vein, using the sub-Gaussian Laplace estimate (5.7) we have

$$\mathbb{P}\left(\mathrm{Tr}(A\mathcal{H}_N) \geq u\right) \leq \mathbb{E}\left[\exp\left(-u\,t + t\,\mathrm{Tr}(A\mathcal{H}_N)\right)\right]$$

This implies that

$$\mathbb{P}\left(\mathrm{Tr}(A\mathcal{H}_N) \geq u\right) \leq \exp\left(-ut + \frac{t^2}{2}\ \left[6\mathrm{Tr}((AP)^2)\right]\right)$$

for any $t < \frac{\sqrt{N}}{4\mathrm{Tr}(AP)}$. Choosing $N$ and $t$ such that

$$t = \frac{u}{6\mathrm{Tr}((AP)^2)} \leq \frac{N^{1/2}}{4\mathrm{Tr}(AP)} \Longleftrightarrow \sqrt{N} \geq \frac{2}{3}\; u\; \frac{\mathrm{Tr}(AP)}{\mathrm{Tr}((AP)^2)}$$

we conclude that

$$\mathbb{P}\left(\mathrm{Tr}(A\mathcal{H}_N) \geq u\right) \leq \exp\left(-\frac{u^2}{12\mathrm{Tr}((AP)^2)}\right)$$

Replacing $u$ by

$$\delta = \frac{u^2}{12\mathrm{Tr}((AP)^2)} \leq \frac{3}{4^2}\; N\; \frac{\mathrm{Tr}((AP)^2)}{\mathrm{Tr}(AP)^2}$$

we prove that

$$\mathbb{P}\left(\mathrm{Tr}(A\mathcal{H}_N) \geq 2\, \sqrt{3\delta\mathrm{Tr}((AP)^2)}\right) \leq \exp\left(-\delta\right)$$

This ends the proof of the theorem. □

## 6.10 Proof of theorem 5.4

We need the following lemma.

**Lemma 6.4.** For any collection of symmetric matrices $(A_i)_{i\in I}$ indexed by some finite set $I$ we have

$$\begin{aligned} & 4\ \log|I| \leq N \\ \Longrightarrow\ & \mathbb{E}\left(\max_{i\in I}\mathrm{Tr}(A_i\mathcal{H}_N)\right) \\ & \leq 2\,\sqrt{\log|I|}\ \max_{i\in I}\ \sqrt{[\mathrm{Tr}((A_iP)^2) + 2\|A_iP\|_F^2]} \end{aligned} \tag{6.21}$$

In addition, we have

$$\mathbb{E}\left(\max_{i\in I}\mathrm{Tr}(A_i\mathcal{H})\right) \leq 2\,\sqrt{\log|I|}\ \max_{i\in I}\ \sqrt{[\mathrm{Tr}((A_iP)^2)]}$$

*Proof.* Applying Jensen's inequality to the logarithm we check that

$$0 < t \leq \frac{\sqrt{N}}{4\max_{i\in I}\|A_iP\|_F}$$

$$\Longrightarrow \quad \mathbb{E}\left(\max_{i\in I}\mathrm{Tr}(A_i\mathcal{H}_N)\right) \leq \frac{1}{t}\ \log\sum_{i\in I}\mathbb{E}\left(\exp\left(t\ \mathrm{Tr}(A_i\mathcal{H}_N)\right)\right)$$

In this situation, using (5.7) we have the estimate

$$\begin{aligned}
&\mathbb{E}\left(\max_{i\in I}\mathrm{Tr}(A_i\mathcal{H}_N)\right)\\
&\leq \frac{1}{t}\ \log\sum_{i\in I}\ \exp\left(t^2\ \left[\mathrm{Tr}((A_iP)^2)+2\|A_iP\|_F^2\right]\right)\\
&\leq \frac{\log|I|}{t} + t\ \max_{i\in I}\ \left[\mathrm{Tr}((A_iP)^2)+2\|A_iP\|_F^2\right]
\end{aligned}$$

as soon as $0 \leq t \leq \frac{\sqrt{N}}{4\max_{i\in I}\|A_iP\|_F}$. We prove (6.21) by choosing

$$t\ =\ \sqrt{\frac{\log|I|}{\max_{i\in I}\ \left[\mathrm{Tr}((A_iP)^2)+2\|A_iP\|_F^2\right]}}\ \leq\ \frac{\sqrt{N}}{4\max_{i\in I}\|A_iP\|_F}$$

as soon as $4\log|I| \leq N$.

This ends the proof of (6.21). Using the Laplace trace formula (5.1) the last assertion follows the same line of arguments, thus it is skipped. This ends the proof of the lemma. □

We now move to prove theorem 5.4. We recall that,

$$\|\mathcal{H}_N\|_{op} = \sup_{x,y\in\mathbb{B}}\langle\mathcal{H}_Nx,y\rangle = \sup_{A\in\mathbb{A}}\mathrm{Tr}(A\mathcal{H}_N)$$

with the sets $(\mathbb{A},\mathbb{B})$ introduced in (1.21). Observe that for any $x,y\in\mathbb{B}$

$$\mathrm{Tr}((xy'P)^2)+2\|xy'P\|_F^2\ =\ \langle x,Py\rangle^2+2\langle x,y\rangle\langle x,P^2y\rangle\ \leq\ 3\ \lambda_1(P)^2$$

$$\Longrightarrow\quad \max_{A\in\mathbb{A}}\ \sqrt{\left[\mathrm{Tr}((AP)^2)+2\|AP\|_F^2\right]}\ \leq\ \sqrt{3}\ \lambda_1(P) \tag{6.22}$$

We recall that $\mathbb{B}$ can be covered by an $\epsilon$-net $\mathbb{B}_\epsilon\subset\mathbb{B}$ of cardinality $|\mathbb{B}_\epsilon| \leq (1+2/\epsilon)^r$, in the sense that for any $x\in\mathbb{B}$ there exists some

$x_\epsilon \in \mathbb{B}_\epsilon$ such that $\|x - x_\epsilon\| \leq \epsilon$. Using the decomposition

$$xy' - x_\epsilon y'_\epsilon = (x - x_\epsilon)(y - y_\epsilon)' + x_\epsilon(y - y_\epsilon)' + (x - x_\epsilon)y'_\epsilon$$

we check that $\mathbb{A}$ can be covered by an $\epsilon$-net $\mathbb{A}_\epsilon \subset \mathbb{A}$ w.r.t. the Frobenius norm of cardinality $|\mathbb{A}_\epsilon| \leq (1 + 6/\epsilon)^{2r}$. In addition, we have

$$\begin{aligned} |\mathrm{Tr}(\mathcal{H}_N(xy' - x_\epsilon y'_\epsilon))| &\leq \frac{2}{3}\,\epsilon \sup_{A\in\mathbb{A}} \mathrm{Tr}(A\mathcal{H}_N) \\ \Longrightarrow \quad \left(1 - \frac{2}{3}\,\epsilon\right) \sup_{A\in\mathbb{A}} \mathrm{Tr}(A\mathcal{H}_N) &\leq \sup_{A\in\mathbb{A}_\epsilon} \mathrm{Tr}(A\mathcal{H}_N) \end{aligned} \tag{6.23}$$

Taking the expectation we conclude that

$$\mathbb{E}\left[\|\mathcal{H}_N\|_{op}\right] \leq\ 2\sqrt{3}\,\left(1 - \frac{2}{3}\,\epsilon\right)^{-1} \sqrt{\log\left(1 + 6/\epsilon\right)}\,\sqrt{r}\,\lambda_1(P)$$

as soon as $N \geq 8r\log\left(1 + 6/\epsilon\right)$. Choosing $\epsilon = 2/10$ we find that

$$N \geq 8r\log(31)$$

$$\Longrightarrow\ \mathbb{E}\left[\|\mathcal{H}_N\|_{op}\right] \leq\ 6\,\sqrt{r}\,\lambda_1(P) \Longrightarrow \mathbb{E}\left[\lambda_1\left(P_N\right)/\lambda_1(P)\right] \leq\ 1 + 6\,\sqrt{\frac{r}{N}}$$

This ends the proof of the first assertion.

Now we come to the proof of (5.9). Combining (5.8) with the estimate (6.22) we check that for any $A \in \mathbb{A}$ the probability of the event

$$\mathrm{Tr}(A\mathcal{H}_N) \geq 2\,\sqrt{3\delta}\,\lambda_1(P) \quad \text{is less than } e^{-\delta}$$

as soon as $0 \leq 8\delta \leq N$. Using the estimate

$$\begin{aligned} &\mathbb{P}\left(\max_{i\in I} \mathrm{Tr}(A_i\mathcal{H}_N) \geq 2\,\sqrt{3\delta}\,\lambda_1(P)\right) \\ &\leq \sum_{i\in I} \mathbb{P}\left(\mathrm{Tr}(A_i\mathcal{H}_N) \geq 2\,\sqrt{3\delta}\,\lambda_1(P)\right) \leq (1 + 6/\epsilon)^{2r}\; e^{-\delta} \end{aligned}$$

which is valid for any $\epsilon$-net $\mathbb{A}_\epsilon = \{A_i\ :\ i \in I\} \subset \mathbb{A}$, we find that the probability of the event

$$\sup_{A\in\mathbb{A}_\epsilon} \mathrm{Tr}(A\mathcal{H}_N) \leq 2\,\sqrt{3\delta}\,\lambda_1(P) \quad \text{is greater than } 1 - (1 + 6/\epsilon)^{2r}\; e^{-\delta}$$

Using (6.23) we also check that the probability of the event

$$\|\mathcal{H}_N\|_{op} \leq 2\sqrt{3}\left(1-\frac{2}{3}\,\epsilon\right)^{-1}\;\sqrt{\delta+2r\;\log\left(1+6/\epsilon\right)}\;\lambda_1(P)$$

is greater than $1-e^{-\delta}$, for any $(N,\delta,\epsilon)$ such that

$$N/8 \geq \delta + 2r\log\left(1+6/\epsilon\right)$$

We end the proof of (5.9) by choosing $\epsilon = 2/10$. This ends the proof of the theorem. □

We come to the proof of corollary 5.5. Following the proof of proposition 5.4 and using (5.4), we check that the probability of the event

$$\|\mathcal{H}\|_{op} \leq \lambda_1(P)\;\left(1-\frac{2}{3}\,\epsilon\right)^{-1}\;2\sqrt{\delta+2r\;\log\left(1+6/\epsilon\right)}$$

is greater than $1-e^{-\delta}$, for any $\delta \geq 0$ and $\epsilon \in ]0,1[$. Choosing $\epsilon = 2/10$ we conclude that for any $\delta \geq 0$ the probability of the event

$$\|\mathcal{H}\|_{op} \leq \;\frac{5}{2}\;\lambda_1(P)\;\sqrt{\delta+7r}\quad\text{is greater than } 1-e^{-\delta} \tag{6.24}$$

By Weyl's inequality (2.5) we have

$$\sup_{1\leq k\leq r}|\lambda_k(P_N)-\lambda_k(P)| \leq \|P_N-P\|_{op} = \frac{1}{\sqrt{N}}\;\|\mathcal{H}_N\|_{op}$$

This completes the proof of the corollary.

# Appendices

# A

# Appendix

This Appendix contains the proofs of a number of results and formulae in the order in which the appear in the main text. In particular, the proof of lemma 2.1, formula (3.2), formula (3.7), formulae (3.5) and (3.8), corollary 3.2, formula (4.5), formula (4.15), formulae (5.1), (5.2), and (5.3), formula (5.6), and formula (5.7) are given here in order.

## A.1 Proof of lemma 2.1

For any $n \geq 0$ and $0 \leq q < p$ we have

$$B_{pn+q}(\overbrace{0,\ldots,0}^{p-1 \text{ terms}}, x_p, x_{p+1}, \ldots, x_{pn+q})$$

$$= \sum_{1\leq k\leq n} \frac{(pn+q)!}{k!} \sum_{m_1+\ldots+m_k=p(n-k)+q,\ m_i\geq 0} \prod_{1\leq i\leq k} \frac{x_{p+m_i}}{(p+m_i)!}$$

We check this claim applying the summation formula

$$\sum_{k\geq 1}\sum_{n\geq kp} a_n\, b_{n,k} = \sum_{1\leq n} \sum_{0\leq q<p} a_{np+q} \sum_{1\leq k\leq n} b_{np+q,k}$$

to the equation

$$\exp\left[\sum_{n\geq p} \frac{t^n}{n!}\, x_n\right]$$

$$= 1 + \sum_{k\geq 1}\sum_{n\geq kp} t^n\, \frac{1}{k!} \sum_{m_1+\ldots+m_k=p(n-k)+q,\ m_i\geq 0} \prod_{1\leq i\leq k} \frac{x_{p+m_i}}{(p+m_i)!}$$

$$= 1 + \sum_{1\leq n}\ \sum_{0\leq q<p} \frac{t^{np+q}}{(np+q)!} \sum_{1\leq k\leq n} \frac{(np+q)!}{k!}$$

$$\times \sum_{m_1+\ldots+m_k=p(n-k)+q,\ m_i\geq 0} \prod_{1\leq i\leq k} \frac{x_{p+m_i}}{(p+m_i)!}$$

Also observe that

$$\exp\left[\sum_{n\geq p} \frac{t^n}{n!}\, x_n\right]$$

$$= 1 + \sum_{1\leq n} \frac{t^{np}}{(np)!} \left[\frac{(np)!}{n!} \left(\frac{x_p}{p!}\right)^n\right.$$

$$\left. + \sum_{1\leq k<n} \frac{(np)!}{k!} \sum_{m_1+\ldots+m_k=p(n-k)} \prod_{1\leq i\leq k} \frac{x_{p+m_i}}{(p+m_i)!}\right]$$

$$+ \sum_{1\leq n}\ \sum_{1\leq q<p} \frac{t^{np+q}}{(np+q)!} \sum_{1\leq k\leq n} \frac{(np+q)!}{k!}$$

$$\times \sum_{m_1+\ldots+m_k=p(n-k)+q,\ m_i\geq 0} \prod_{1\leq i\leq k} \frac{x_{p+m_i}}{(p+m_i)!}$$

Under our assumption, we have

$$\sum_{1\leq k<n} \frac{(np)!}{k!} \sum_{m_1+\ldots+m_k=p(n-k)} \prod_{1\leq i\leq k} \frac{x_{p+m_i}}{(p+m_i)!}$$

$$\leq \beta^{pn} \sum_{1\leq k<n} \frac{(np)!}{k!} \rho^k \sum_{m_1+\ldots+m_k=p(n-k)} 1$$

$$\leq \beta^{pn} \sum_{1\leq k<n} \frac{(np)!}{k!} \rho^k \, 2^{p(n-k)+k-1}$$

$$= \frac{1}{2} (2\beta)^{pn} \sum_{1\leq k<n} \frac{(np)!}{k!} \left(\frac{\rho}{2^{p-1}}\right)^k$$

In the same vein, we have

$$\sum_{1\leq k\leq n} \frac{(np+q)!}{k!} \sum_{m_1+\ldots+m_k=p(n-k)+q,\ m_i\geq 0} \prod_{1\leq i\leq k} \frac{x_{p+m_i}}{(p+m_i)!}$$

$$\leq \beta^{pn+q} \sum_{1\leq k\leq n} \frac{(np+q)!}{k!} \rho^k \sum_{m_1+\ldots+m_k=p(n-k)+q,\ m_i\geq 0} 1$$

$$\leq \frac{1}{2} (2\beta)^{pn+q} \sum_{1\leq k\leq n} \frac{(np+q)!}{k!} \left(\frac{\rho}{2^{p-1}}\right)^k$$

This completes the proof of this lemma. □

## A.2 Proof of (3.2)

The proof is partly based on the estimate

$$\frac{2^{-n}}{n!} \mathbb{E}(\mathbb{X}^n) \leq \frac{1}{2n} \sum_{0\leq k<n} \frac{2^{-k}}{k!} \upsilon_{2k} \operatorname{Tr}(P)^k \quad P^{n-k} \tag{A.1}$$

For $n = 1, 2$ the inequality (A.1) reduces to an equality. The proof of the above formulae can be found in [10]. We also have the identity

$$\sum_{0\leq k\leq n} \frac{(2k)!}{2^{2k}k!^2} = \frac{(2n+1)!}{2^{2n}n!^2} \leq \frac{2n+1}{\sqrt{3n+1}} \tag{A.2}$$

The proof of the first assertion in the above display can be found in [70] (cf. formula (1.38) in Section 1.2). The r.h.s. estimate comes from the central binomial estimates

$$\frac{1}{\sqrt{4n}} \leq \frac{(2n)!}{2^{2n} n!^2} \leq \frac{1}{\sqrt{3n+1}} \tag{A.3}$$

The proof of the above estimate can be found in [39]. Using (A.1) we check the trace estimates

$$\begin{aligned}\mathrm{Tr}\left[\mathbb{E}\left(\mathbb{X}^n\right)\right] = \mathbb{E}\left(\|X\|^{2n}\right) &\leq \frac{(2n)!}{2^n n!}\ \mathrm{Tr}(P)^n \\ &\leq \frac{2^n}{\sqrt{3n+1}}\ n!\ \mathrm{Tr}(P)^n \leq 2^{n-1}\ n!\ \mathrm{Tr}(P)^n\end{aligned}$$

On the other hand, we have the rather crude estimates

$$\|\mathbb{E}\left[M_{2n,n}(P)\right]\| \leq \frac{(2n)!}{2^n n!}\ \mathbb{E}\left(\|\mathbb{X}-P\|^2\right)^n \tag{A.4}$$

for any compatible matrix norm. For instance, for the Frobenius norm we find that

$$\|M_{2n,n}(P)\|_F \leq \frac{(2n)!}{2^n n!}\ \left[\mathrm{Tr}\left(\mathbb{E}(\mathcal{H}^2)\right)\right]^n = \frac{(2n)!}{2^n n!}\ \left[\mathrm{Tr}(P^2)+\mathrm{Tr}(P)^2\right]^n$$

Using (A.1) and (A.3) we have the rather crude estimate

$$\begin{aligned}&\mathbb{E}\left[\|\mathbb{X}-P\|_F^n\right] \\ &\leq 2^{n-1}\ \left[\mathbb{E}\left[\|X\|^{2n}\right]+\|P\|_F^n\right] \\ &= 2^{n-1}\ \left[2^{n-1}\ n!+1\right]\ \mathrm{Tr}(P)^n \leq \frac{n!}{2}\ (4\mathrm{Tr}(P))^n\end{aligned}$$

Therefore, for any partition $\pi$ of $[n]$ with $m$ blocks $\pi_i$ of size $|\pi_i| = n_i \geq 2$ we have

$$\|M_\pi(P)\|_F \leq \prod_{1\leq i\leq m} \mathbb{E}\left[\|\mathbb{X}-P\|_F^{n_i}\right] \leq 2^{-m}\ (4\mathrm{Tr}(P))^n \prod_{1\leq i\leq m} n_i!$$

This implies that

$$\|M_{n,m}(P)\|_F \leq \frac{2^{-m}}{m!} \sum_{n_1+\ldots+n_m=n} \frac{n!}{n_1!\ldots n_m!} \prod_{1\leq i\leq m} \mathbb{E}\left[\|\mathbb{X}-P\|_F^{n_i}\right]$$

$$\leq \frac{n!}{m!}\; 2^{-m}\; (4\mathrm{Tr}(P))^n \; \left(\begin{array}{c} n-2m+m-1 \\ n-2m \end{array}\right)$$

$$\leq \frac{n!}{m!}\; (8\mathrm{Tr}(P))^n \; 2^{-(2m+1)}$$

The last assertion comes from the fact that the number of solutions to the equation

$$j_1+j_2+\ldots+j_p = q \quad \text{with} \quad j_i \geq 0 \quad \text{is given by} \quad \left(\begin{array}{c} p+q-1 \\ q \end{array}\right)$$

and $\left(\begin{array}{c} n \\ p \end{array}\right) \leq 2^n$ for any $0 \leq p \leq n$. This ends the proof of (3.2). □

## A.3 Proof of (3.7)

Recall that

$$P = I \implies \mathbb{X} \overset{law}{=} \chi_r \, \mathbb{U} \implies \mathbb{X} - I = (\chi_r - 1)\, \mathbb{U} + (\mathbb{U} - I) \qquad \text{(A.5)}$$

where $\mathbb{U} := \overline{X}\,\overline{X}'$ stands for the random projection matrix associated with an uniform random vector $\overline{X} := X/\sqrt{\|X\|}$ on the unit sphere, and $\chi_r = \|X\|^2$ is an independent $\chi$-square random variable with $r$ degrees of freedom.

We use the binomial formula

$$\mathbb{E}\left[(\mathbb{X}-I)^n\right] = (-1)^n \left(1-\frac{1}{r}\right) I + \mathbb{E}\left(\left[\|X\|^2-1\right]^n\right) \frac{1}{r}\; I$$

and

$$\mathbb{E}\left[\|X\|^{2k}\right] = \prod_{0\leq l<k} (r+2l)$$

The above formula can be checked using (A.5). Notice that $\mathbb{U}^n = \mathbb{U}$ and $(\mathbb{U}-I)^n = (-1)^{n+1}(\mathbb{U}-I) \implies \mathbb{U}(\mathbb{U}-I)^n = 0 = \mathbb{U}(\mathbb{U}-I)^n$ for any

$n \geq 1$. This yields

$$\mathbb{E}\left[(\mathbb{X}-I)^n\right] = \; r^{n-1} \; \left[(-1)^n \; \left(\frac{1}{r}\right)^{n-1} \right.$$

$$\left. + \sum_{1\leq l\leq n} \left(\begin{array}{c} n \\ l \end{array}\right) \; \left(-\frac{1}{r}\right)^{n-l} \; \prod_{0\leq j<l} \left(1+\frac{2j}{r}\right)\right] \; I$$

Integrating sequentially the non-crossing indices, for any partition of $[n]$ of type $\mu \vdash n$ we find that

$$M_\pi(I)$$

$$= 1_{\mu_1=0} \; \prod_{2\leq k\leq n} \; \mathbb{E}\left((\mathbb{X}-I)^k\right)^{\mu_k} = r^{n-m} \; 1_{\mu_1=0} \; \prod_{2\leq k\leq n}$$

$$\left[(-1)^n \; \left(\frac{1}{r}\right)^{k-1} + \sum_{1\leq l\leq k} \left(\begin{array}{c} k \\ l \end{array}\right) \; \left(-\frac{1}{r}\right)^{k-l} \; \prod_{0\leq j<l} \left(1+\frac{2j}{r}\right)\right]^{\mu_k} \; I$$

Also observe that for any $2t < 1$ we have

$$\log \mathbb{E}\left(e^{t(\|X\|^2-1)}\right) = (r-1)\; t + \frac{r}{2} \sum_{n\geq 2} \; \frac{(2t)^n}{n}$$

This yields the formula

$$\mathbb{E}\left(\left[\|X\|^2-1\right]^n\right)$$

$$= B_n\left(r-1,\; r\; 2^{2-1}\; (2-1)!, \ldots,\; r\; 2^{n-1}\; (n-1)!\right) \geq 1$$

In this situation, the polynomial formula (3.12) reduces to

$$\mathbb{E}\left[(\mathbb{X}-I)^n\right] = \left[\frac{1}{r}\; B_n\left(r-1,\; r\; 2^{2-1}\; (2-1)!, \ldots,\; r\; 2^{n-1}\; (n-1)!\right)\right.$$

$$\left. +(-1)^n \; \left(1-\frac{1}{r}\right)\right] \; I$$

A more thorough discussion on more general rank-one matrix moments can be found in [10]. Finally observe that

$$r^{-(n-1)}\mathbb{E}\left[(\mathbb{X}-I)^n\right] = \left[(-1)^n \left(\frac{1}{r}\right)^{n-1} + \prod_{0\leq j<n}\left(1+\frac{2j}{r}\right) + \sum_{1\leq l<n}\binom{n}{l}\left(-\frac{1}{r}\right)^{n-l}\prod_{0\leq j<l}\left(1+\frac{2j}{r}\right)\right] I$$

This implies that

$$r^{-(n-1)}\mathbb{E}\left[(\mathbb{X}-I)^n\right] \longrightarrow_{r\to\infty} (1_{n>1}+\mathrm{O}(1/r))\ I$$

The end of the proof of (3.9) is now easily completed. This ends the proof of (3.7). □

## A.4 Proof of (3.5) and (3.8)

For any $n\geq 1$ and any collection of symmetric matrices $Q_1,\ldots,Q_n$ with $Q_1=I$ we have the reduction formula

$$\mathbb{E}((\mathbb{X}Q_2)(\mathbb{X}Q_3)\ldots(\mathbb{X}Q_n)\mathbb{X}) = \mathbb{E}\left((\mathbb{X}Q_1)\ \mathrm{Tr}\left[(\mathbb{X}Q_2)(\mathbb{X}Q_3)\ldots(\mathbb{X}Q_n)\right]\right)$$

$$= \sum_{1\leq m\leq n}\sum_{\tau\in\mathcal{P}_{n,m}} 2^{n-m} \times\left[\prod_{1\leq k\leq n}(k-1)!^{r_k(\tau)}\right]\left[\prod_{2\leq i\leq m}\mathrm{Tr}(Q_{\tau_i})\right]\ (Q_{\tau_1})_{sym} \tag{A.6}$$

In the above display, $\tau_1\leq\ldots\leq\tau_m$ stands for the $m$ ordered blocks of the partition $\tau$ ($\Longrightarrow \tau_1\ni 1$), $r_k(\tau)$ stands for the number of blocks of size $k$, and

$$\tau_i=\{j_1^i,\ldots,j_{|\tau_i|}^i\}\quad\text{with}\quad j_1^i\leq\ldots\leq j_{|\tau_i|}^i \Longrightarrow Q_{\tau_i}=Q_{j_1^i}Q_{j_2^i}\ldots Q_{j_{|\tau_i|}^i}$$

The proof of the above formulae can be found in [10]. We further assume that $\pi\in\mathcal{P}_{n,m}-\mathcal{N}_{n,m}$ is an $\mathbb{X}$-connected crossing partition of $[n]$ in the sense that $\mathbb{X}_\pi$ cannot be written in terms of the product of two independent random matrices. In this situation, we have

$$\mathbb{X}_\pi = [(\mathbb{X}_1 Q_2)(\mathbb{X}_1 Q_3)\ldots(\mathbb{X}_1 Q_{m_1})\mathbb{X}_1]\, Q_{m_1+1} \quad \text{with} \quad m_1 = |\pi_1|$$

and some random matrices $Q_j$ independent of $\mathbb{X}_1$ and such that

$$h(\pi) := |\{2 \leq i \leq m_1 \ : \ Q_i = I\}| \leq m_1 - 2$$

For any rank one matrix $\mathbb{X} = XX'$, for some column vector $X \in \mathbb{R}^r$, and for any collection of matrices $(R_i)_{1\leq i\leq p}$ and $(S_i)_{0\leq i\leq q}$ we have

$$\begin{aligned}&\mathrm{Tr}\,[(\mathbb{X}R_1)\ldots(\mathbb{X}R_p)]\ \ S_0\ \ [(\mathbb{X}S_1)\ldots(\mathbb{X}S_p)]\\ &= S_0\ \ [(\mathbb{X}R_1)\ldots(\mathbb{X}R_p)]\ \ [(\mathbb{X}S_1)\ldots(\mathbb{X}S_p)]\end{aligned}$$

Combining the above decomposition with the reduction formula (A.6) we check that

$$\mathbb{E}\left(\mathbb{X}_\pi\right) = \sum_{l=0}^{h(\pi)} \sum_{k\in J_\pi} r^l \ \alpha_{k,l}(\pi)\ \mathbb{E}\left(\mathbb{X}_{\varphi_{k,l}(\pi)}\right)$$

for any $\mathbb{X}$-connected crossing partition $\pi$ of $[n]$ with $b$ blocks.

In the above display, $\alpha_{k,l}(\pi) \geq 0$ stands for some non-negative parameters whose values doesn't depend on $r$, $\varphi_{k,l}(\pi)$ a collection of partitions on $[n-m_1]$ with $(b-1)$ blocks, and $J_\pi$ is a set with cardinality less than $2B(m_1)$, with the Bell number $B(m_1) := \sum_{1\leq k\leq m_1} S(m_1, k)$.

We further assume that for any $\overline{\pi} \in \mathcal{P}_{\overline{n},\overline{m}} - \mathcal{N}_{\overline{n},\overline{m}}$ we have

$$\mathbb{E}(\mathbb{X}_{\overline{\pi}}) = c_{\overline{\pi}}\ r^{\overline{n}-\overline{m}}\ \ (1 + \mathrm{O}(1/r))\ \ r^{-1}\ I$$

for any $1 \leq \overline{m} \leq \overline{n} < n$ and some finite constant $c_\pi$ whose values doesn't depend on $r$. This result is clearly met for $n = 5$. In this case $\mathcal{P}_{4,2} - \mathcal{N}_{4,2}$ reduces to the partition with $b = 2$ blocks

$$\pi_1 = \{1, 3\} \leq \pi_2 = \{2, 4\} \Longrightarrow \mathbb{E}(\mathbb{X}_1\mathbb{X}_2\mathbb{X}_1\mathbb{X}_2) = 3\ (r + 2)\ I$$

Under the above induction hypothesis and applying (3.4) when $\varphi_{k,l}(\pi) \in \mathcal{N}_{n,m}$ we have

$$\begin{aligned}\mathbb{E}\left(\mathbb{X}_\pi\right) &= c_\pi\ r^{m_1-2}\ r^{(n-m_1)-(b-1)}\ \ (1 + \mathrm{O}(1/r))\ \ I\\ &= c_\pi\ r^{n-b}\ \ (1 + \mathrm{O}(1/r))\ \ r^{-1}\ I\end{aligned}$$

When $\pi$ is crossing but non-necessarily $\mathbb{X}$-connected we have

$$\mathbb{E}\left(\mathbb{X}_\pi\right) = \mathbb{E}\left(C_1\right) \ldots \mathbb{E}\left(C_k\right)$$

with some $\mathbb{X}$-connected components $C_i$, with $1 \leq i \leq k$. Each of them, say $C_i$ is associated with a partition of $[n_i]$ with $b_i$ blocks for some parameters $n = n_1 + \ldots + n_k$ and $b = b_1 + \ldots + b_k$. Since $\pi$ is crossing at least one of the components, say $C_j$ is crossing. In this case we have

$$\begin{aligned}\mathbb{E}\left(\mathbb{X}_\pi\right) &= c_\pi \ r^{n_1-b_1} \quad r^{n_k-b_k} \ \left(1+\mathrm{O}(1/r)\right) \ r^{-1} \ I \\ &= c_\pi \ r^{n-b} \ \left(1+\mathrm{O}(1/r)\right) \ r^{-1} \ I\end{aligned}$$

This ends the proof of (3.5). □

The proof of (3.8) relies on the decomposition

$$(\mathbb{X}_{\alpha^\pi(1)} - I) \ldots (\mathbb{X}_{\alpha^\pi(n)} - I) = \sum_{\epsilon=(\epsilon_1,\ldots,\epsilon_n)\in\{0,1\}^n} (-1)^{n-\sum_{1\leq i\leq n}\epsilon_i} \ \mathbb{X}_\pi^\epsilon$$

with

$$\mathbb{X}_\pi^\epsilon := \ \mathbb{X}_{\alpha^\pi(1)}^{\epsilon_1} \ldots \mathbb{X}_{\alpha^\pi(n)}^{\epsilon_n}$$

Using (2.14) we check that

$$\pi \in \mathcal{Q}_{n,b} \quad \text{and} \quad b > \lfloor n/2 \rfloor \implies \mathbb{E}(M_\pi(P)) = 0$$

We further assume that $\pi \in \mathcal{Q}_{n,b}^-$ is a crossing partition with a number of blocks $b \leq \lfloor n/2 \rfloor$.

For any $\epsilon = (\epsilon_1, \ldots, \epsilon_n) \in \{0,1\}^n$ with $\sum_{1\leq i\leq n} \epsilon_i = m$ we have

$$\mathbb{X}_\pi^\epsilon \stackrel{law}{=} \mathbb{X}_{\varphi_\epsilon(\pi)}$$

for some partition $\varphi_\epsilon(\pi)$ of $[m]$ with a number of blocks

$$b - \lfloor (n-m)/2 \rfloor \leq b_m := |\varphi_\epsilon(\pi)| \leq b \wedge m$$

For any $m > n - 2b$ the partition $\varphi_\epsilon(\pi)$ has at least $b_m \geq b - \lfloor (n - m)/2 \rfloor \geq 1$ blocks. In this situation, since $\varphi_\epsilon(\pi)$ may be a crossing partition by (3.4) and (3.5) we have

$$\begin{aligned}\mathbb{E}\left(\mathbb{X}_{\varphi_\epsilon(\pi)}\right) &= r^{m-b_m} \ \left(1+\mathrm{O}(1/r)\right) \ I \\ &\leq r^{m+\lfloor (n-m)/2\rfloor - b} \ \left(1+\mathrm{O}(1/r)\right) \ I \\ &\leq r^{(n+m)/2-b} \ \left(1+\mathrm{O}(1/r)\right) \ r^{-1} \ I \\ &\leq r^{n-b} \ \left(1+\mathrm{O}(1/r)\right) \ r^{-1} \ I\end{aligned}$$

In the last assertion we have used the fact that $m \leq n$ and $\lfloor k/2 \rfloor \leq k-1$, which is valid for any $k \geq 1$. In the reverse angle, when $m \leq n-2b$ the partition $\varphi_\epsilon(\pi)$ may have a single block. In this situation, arguing as above we have

$$\mathbb{E}\left(\mathbb{X}_{\varphi_\epsilon(\pi)}\right) = r^{m-1}\ (1+\mathrm{O}(1/r))\ I \leq r^{n-2b-1}\ (1+\mathrm{O}(1/r))\ I \\ \leq r^{n-b}\ (1+\mathrm{O}(1/r))\ r^{-1}\ I$$

This ends the proof of (3.8). □

## A.5 Proof of corollary 3.2

Combining (1.4) with theorem 3.1, for any $2N \geq n \geq 1$ we have

$$\begin{aligned}\mathbb{E}\left[P_N^n\right] &= P^n + \sum_{1\leq k\leq n}\ \sum_{1\leq l\leq \lfloor k/2\rfloor} \frac{1}{N^{k-l}}\ (n)_k\ \partial_{k,l}(P) \\ &= P^n + \sum_{1\leq l\leq \lfloor n/2\rfloor}\ \sum_{2l\leq k\leq n} \frac{1}{N^{k-l}}\ (n)_k\ \partial_{k,l}(P)\end{aligned}$$

The last assertion comes from the iterative series formula

$$\sum_{1\leq k\leq n}\ \sum_{1\leq l\leq \lfloor k/2\rfloor} a_{k,l} = \sum_{1\leq l\leq \lfloor n/2\rfloor}\ \sum_{2l\leq k\leq n} a_{k,l}$$

applied to the array

$$a_{k,l} := \frac{1}{N^{k-l}}\ (n)_k\ \partial_{k,l}(P)$$

Assume that $n = 2m$. In this case, we have

$$\sum_{1\leq l\leq \lfloor n/2\rfloor}\ \sum_{2l\leq k\leq n} a_{k,l} = \sum_{1\leq l\leq m}\ \sum_{2l\leq k\leq n} a_{l+(k-l),l} = \sum_{1\leq l\leq m}\ \sum_{k=l}^{2m-l} b_{k,l}$$

with the array

$$b_{k,l} := a_{l+k,l} = \frac{1}{N^k}\ (n)_{l+k}\ \partial_{l+k,l}(P)$$

Observe that

$$\begin{aligned}\sum_{1\leq l\leq m}\sum_{k=l}^{2m-l} b_{k,l} &= \sum_{1\leq p\leq m}\sum_{1\leq k\leq p} b_{p,k}+\sum_{1\leq p<m}\sum_{1\leq k\leq m-p} b_{m+p,k}\\ &=\sum_{1\leq p\leq m}\frac{1}{N^p}\sum_{1\leq k\leq p}(n)_{p+k}\ \partial_{p+k,k}(P)\\ &\quad+\sum_{1\leq p<m}\frac{1}{N^{m+p}}\sum_{1\leq k\leq m-p}(n)_{k+m+p}\ \partial_{m+p+k,k}(P)\end{aligned}$$

This yields the formula

$$\begin{aligned}\sum_{1\leq l\leq m}\sum_{k=l}^{2m-l} b_{k,l} &= \sum_{1\leq p\leq m}\frac{1}{N^p}\sum_{p<k\leq 2p}(n)_k\ \partial_{k,k-p}(P)\\ &\quad+\sum_{m<p<2m}\frac{1}{N^p}\sum_{p<k\leq 2m}(n)_k\ \partial_{k,k-p}(P)\\ &=\sum_{1\leq p<n}\frac{1}{N^p}\sum_{p<k\leq(2p)\wedge n}(n)_k\ \partial_{k,k-p}(P)\end{aligned}$$

For odd parameters $n=(2m+1)$ we have

$$\begin{aligned}&\sum_{1\leq k\leq 2m+1}\sum_{1\leq l\leq\lfloor k/2\rfloor} a_{k,l}\\ &=\sum_{1\leq l\leq m} a_{2m+1,l}+\sum_{1\leq l\leq m}\sum_{2l\leq k\leq 2m} a_{k,l}\\ &=\sum_{1\leq p\leq m}\frac{1}{N^p}\sum_{p<k\leq 2p}(n)_k\ \partial_{k,k-p}(P)+\frac{1}{N^{n-1}}\ \delta_{n,1}(P)\\ &\quad+\sum_{m<p<2m}\frac{1}{N^p}\left[\partial_{n,n-p}(P)+\sum_{p<k\leq 2m}(n)_k\ \partial_{k,k-p}(P)\right]\\ &=\sum_{1\leq p\leq m}\frac{1}{N^p}\sum_{p<k\leq 2p}(n)_k\ \partial_{k,k-p}(P)\\ &\quad+\frac{1}{N^{n-1}}\ \delta_{n,1}(P)+\sum_{m<p<2m}\frac{1}{N^p}\sum_{p<k\leq n}(n)_k\ \partial_{k,k-p}(P)\end{aligned}$$

This yields the formula

$$\sum_{1\leq k\leq n}\ \sum_{1\leq l\leq \lfloor k/2\rfloor} a_{k,l} = \sum_{1\leq p\leq m} \frac{1}{N^p} \sum_{p<k\leq (2p)\wedge n} (n)_k\ \partial_{k,k-p}(P) + \sum_{m<p<n} \frac{1}{N^p} \sum_{p<k\leq (2p)\wedge n} (n)_k\ \partial_{k,k-p}(P)$$

This ends the proof of the first assertion. The last differential-type matrix moment formula comes from the fact that

$$\begin{aligned}
\mathbb{E}(P_N^n) &= \frac{1}{N^n} \sum_{1\leq m\leq n} (N)_m\ M^{\circ}_{n,m}(P) \\
&= \frac{1}{N^n} \sum_{1\leq l\leq m\leq n} s(m,l)\ N^l\ M^{\circ}_{n,m}(P) \\
&= \sum_{1\leq l\leq n} \frac{1}{N^{n-l}} \sum_{l\leq m\leq n} s(m,l)\ M^{\circ}_{n,m}(P) \\
&= \sum_{0\leq k<n} \frac{1}{N^{k}} \sum_{n-k\leq m\leq n} s(m,n-k)\ M^{\circ}_{n,m}(P)
\end{aligned}$$

This ends the proof of the corollary. □

## A.6 Proof of (4.5)

Observe that for any $n\geq 1$ and $m\geq 0$ we have

$$\mathbb{X}^{n+m} = \|X\|^{2m}\ \mathbb{X}^n = (\mathrm{Tr}\,[\mathbb{X}])^m\ \ \mathbb{X}^n \Longrightarrow \mathbb{X}^n \exp\,[t\,\mathbb{X}] = \mathbb{X}^n \exp\,[t\,\mathrm{Tr}\,[\mathbb{X}]]$$

In addition, we have the exponential formulae

$$\begin{aligned}
\mathbb{X}^n &= \|X\|^{2(n-1)}\ \mathbb{X} = \frac{\mathbb{X}}{\mathrm{Tr}(\mathbb{X})}\ \|X\|^{2n} \\
&\Longrightarrow \exp\,[t\,\mathbb{X}] = I + \frac{\mathbb{X}}{\mathrm{Tr}(\mathbb{X})}\ \exp\,[t\,\mathrm{Tr}\,[\mathbb{X}]]
\end{aligned} \tag{A.7}$$

Assume that $P=I$. In this case the projection matrix $\mathbb{X}/\mathrm{Tr}(\mathbb{X})$ and the $\chi$-square variable $\mathrm{Tr}(\mathbb{X})$ are independent. This yields for any $t<1/2$ the formula

$$\begin{aligned}
\mathbb{E}\,(\exp\,[t\,\mathbb{X}]) &= I + \mathbb{E}\,(\mathbb{X}/\mathrm{Tr}(\mathbb{X}))\ \ \mathbb{E}\,(\exp\,[t\,\mathrm{Tr}\,[\mathbb{X}]]) \\
&= \left(1 + \frac{1}{r\ (1-2t)^{r/2}}\right)\ I
\end{aligned}$$

For any $n \geq 1$ and any parameter $t$ s.t. $I - 2tP$ is invertible we have the change of variable formula

$$\mathbb{E}\left(\mathbb{X}^n \exp\left[t\,\mathbb{X}\right]\right) = \mathbb{E}\left(\mathbb{X}^n\ \exp\left[t\,\mathrm{Tr}\left[\mathbb{X}\right]\right]\right) = \mathbb{E}\left(\exp\left[t\,\mathrm{Tr}\left[\mathbb{X}\right]\right]\right)\ \mathbb{E}\left[\mathbb{Y}_t^n\right] \tag{A.8}$$

where $\mathbb{Y}_t$ stands for Wishart distribution with a one degree of freedom and covariance matrix $P_t := [I - 2tP]^{-1}P$. We check (A.8) recalling that $\mathbb{X}$ has a Wishart distribution with a one degree of freedom and covariance matrix $P$.

Using formulae (A.1), the matrix moments (A.8) can be computed explicitly in terms of the covariance matrix $P$. For instance, we have

$$\partial_t \mathbb{E}\left(\exp\left[t\,\mathbb{X}\right]\right) = \mathbb{E}\left(\exp\left(t\,\mathrm{Tr}\left[\mathbb{X}\right]\right)\right)\left[I - 2tP\right]^{-1}P \iff (4.5)$$

□

## A.7 Proof of (4.15)

Applying (4.14) we prove that

$$\mathbb{E}\left(\mathrm{Tr}\left[\exp\left(\frac{t}{2}\ \frac{\mathcal{H}_N}{\mathrm{Tr}(P)}\right)\right]\right) \leq r\ \exp\left(NL\left(\frac{t}{\sqrt{N}}\right)\ \varpi_+\left(\frac{P}{\mathrm{Tr}(P)}\right)\right)$$

Combining the above estimate with the Laplace transform method for random matrices presented in proposition 3.1 in [93] we check that

$$\begin{aligned}
&\mathbb{P}\left(\frac{1}{2}\frac{\lambda_1(\mathcal{H}_N)}{\mathrm{Tr}(P)} \geq v\right) \\
&\leq \inf_{t\in[0,1]}\left\{e^{-vt}\ \mathbb{E}\left(\mathrm{Tr}\left[\exp\left(\frac{t}{2}\ \frac{\mathcal{H}_N}{\mathrm{Tr}(P)}\right)\right]\right)\right\} \\
&\leq r\ \exp\left[-N\varpi_+\left(\frac{P}{\mathrm{Tr}(P)}\right)\ L^\star\left(\frac{v}{\sqrt{N}\,\varpi_+\left(\frac{P}{\mathrm{Tr}(P)}\right)}\right)\right]
\end{aligned}$$

We end the proof of (4.15) using (4.6). This ends the proof of the estimate (4.15).

Combining (4.15) with the almost sure estimate

$$\lambda_1(P_N) \leq \lambda_1(P) + \frac{1}{\sqrt{N}}\ \lambda_1(\mathcal{H}_N)$$

we also check that the probability of the event

$$\lambda_1(P_N) \leq \lambda_1(P)$$
$$+2\,\mathrm{Tr}(P)\left[\frac{\delta+\log(r)}{N}+2\sqrt{\frac{\delta+\log(r)}{N}\,\varpi_+\left(\frac{P}{\mathrm{Tr}(P)}\right)}\right] \tag{A.9}$$

is greater than $1-e^{-\delta}$. □

## A.8 Proof of (5.1), (5.2), and (5.3)

We use the fact that

$$\mathrm{Tr}([A\mathcal{H}]) = \sum_{u\in([r]\times[r])} A(u)\ \mathcal{H}(u)$$

where $\mathcal{H}(u)$ are Gaussian random variables with covariance matrix

$$C(u,v) = \mathbb{E}(\mathcal{H}(u)\mathcal{H}(v)) = 2(P\,\widehat{\otimes}\,P)_{u,v}$$

$$\Longrightarrow \log\mathbb{E}\left(\exp\left[t\,\mathrm{Tr}(A\mathcal{H})\right]\right) = \frac{t^2}{2}\sum_{u,v\in([r]\times[r])} C(u,v)\ A(u)\ A(v)$$

This ends the proof of (5.1). □

To check (5.2) we observe that

$$\|\mathcal{H}\|_F^2 = \sum_{u\in([r]\times[r])} \mathcal{H}(u)^2$$

This implies that

$$\log\mathbb{E}\left(\exp\left(t\|\mathcal{H}\|_F^2\right)\right) = \frac{1}{2}\sum_{n\geq 1}\frac{(2t)^n}{n}\,\mathrm{Tr}(C^n)$$

This assertion comes from the fact that $\sum_{u\in([r]\times[r])}\ \mathcal{H}(u)^2$ is a $\chi$-square random variable associated with a centered Gaussian $r^2$-dimensional vector with covariance matrix $C$. On the other hand, using (1.3) we have

$$C = 2\ (P\,\widehat{\otimes}\,P) \quad \text{and} \quad (P\,\widehat{\otimes}\,P)^n = (P^n\,\widehat{\otimes}\,P^n)$$

$$\Longrightarrow \mathrm{Tr}(C^n) = 2^n\ \mathrm{Tr}(P^n\,\widehat{\otimes}\,P^n) = 2^{n-1}\ \left[\mathrm{Tr}(P^n)^2+\mathrm{Tr}(P^{2n})\right]$$

This ends the proof of (5.2). □

Also observe that

$$\mathbb{E}\left(\|\mathcal{H}\|_F^{2n}\right) = B_n\left(\mathrm{Tr}(P)^2 + \mathrm{Tr}(P^2), 4^{2-1}\ (2-1)!\left[\mathrm{Tr}(P^2)^2 + \mathrm{Tr}(P^4)\right],\right.$$
$$\left.\ldots, 4^{n-1}\ (n-1)!\left[\mathrm{Tr}(P^n)^2 + \mathrm{Tr}(P^{2n})\right]\right)$$

with the complete Bell polynomials $B_n$ defined in (2.16).

Now we come to the proof of (5.3). Let $e_k$ be the $r$-column vectors with null entries but the $k$-th unit one. Choosing $A = e_{k,l}$ in (5.5), we find that

$$\begin{aligned}\mathrm{Tr}\left(A\mathcal{H}_N\right) &= \mathcal{H}_N(k,l)\\ \mathrm{Tr}((AP)^2) &= \frac{P(k,l)^2 + P(k,l)P(l,l)}{2}\\ &= (P \,\widehat{\otimes}\, P)_{(k,l),(k,l)} = \mathbb{E}(\mathcal{H}(k,l)^2)/2\end{aligned}$$

and

$$\|AP\|_F^2 = \sigma_{k,l}^2 := \frac{P^2(k,k) + P^2(l,l)}{4} + 2^{-1}P^2(k,k)\ 1_{k=l}$$

This implies that

$$\mathbb{E}\left(\cosh\left(t\ \mathcal{H}_N(k,l) - \frac{t^2}{2}\ \mathbb{E}(\mathcal{H}(k,l)^2)\right)\right)$$
$$\leq \exp\left[-2t^2\ \sigma_{k,l}^2 \log\left(1 - \frac{2|t|}{N^{1/2}}\ \sigma_{k,l}\right)\right]$$

On the other hand, for any $0 < t \leq (2\sigma_{k,l})^{-1}$ we have

$$\mathbb{E}\left(\exp\left[t\ |\mathcal{H}_N(k,l)|\right]\right) \leq 2\ \mathbb{E}\left(\cosh\left(t\ \mathcal{H}_N(k,l)\right)\right)$$

Using Hölder's inequality we find the rather crude estimate

$$\|\mathcal{H}_N\|_F \leq \sum_{1\leq k,l\leq r} |\mathcal{H}_N(k,l)|$$
$$\Longrightarrow\quad \mathbb{E}\left(\exp\left[t\ \|\mathcal{H}_N\|_F - \frac{(rt)^2}{2}\ \mathbb{E}(\|\mathcal{H}\|_F^2)\right]\right)$$
$$\leq\ 2\ \exp\left[-2(rt)^2\ (r\sigma)^2 \log\left(1 - \frac{2r^2t}{N^{1/2}}\ \sigma\right)\right]$$

as soon as $0 < t \leq (2r\sigma)^{-1}$ with $\sigma := \max_{k,l} \sigma_{k,l}$. We conclude that

$$\exists t > 0 \quad \text{s.t.} \quad \sup_{N\geq 1} \mathbb{E}\left[\exp\left(t\,\|\mathcal{H}_N\|_F\right)\right] < \infty$$

$$\Longleftrightarrow \exists c < \infty \quad \text{s.t.} \quad \sup_{N\geq 1} \mathbb{E}\left[\|\mathcal{H}_N\|_F^n\right] \leq c^n\; n!$$

for some constant $c$ whose value depends on $P$. The proof of the equivalence between the growth of the moments and the Hardy's exponential condition can be found in [53, 79] (see for instance Section 5, theorem 3 in [79]). This ends the proof (5.3). □

## A.9 Proof of (5.6)

The proof of (5.6) is a direct consequence of the following lemma.

**Lemma A.1.** For any $n \geq 1$, $A = (xy' + yx')/2$ with $x, y \in \mathbb{R}^r$ we have

$$4\,(AP)^n A = a_n(x,y)\;\; (xy' + yx') + b_n(x,y)\,\left[\langle x, Px\rangle\;\, yy' + \langle y, Py\rangle\;\, xx'\right]$$

with the symmetric functions

$$\begin{aligned} a_n(x,y) \;&=\; 2^{-(n-1)} \sum_{0\leq l\leq \lfloor n/2\rfloor} \binom{n}{2l} \;\langle x, Py\rangle^{n-2l}\;\; \left(\langle x, Px\rangle\,\langle y, Py\rangle\right)^l \\ b_n(x,y) &:= 2^{-(n-1)} \sum_{0\leq l\leq \lfloor (n-1)/2\rfloor} \binom{n}{2l+1} \\ &\qquad\qquad \times \langle x, Py\rangle^{n-1-2l}\;\; \left(\langle x, Px\rangle\,\langle y, Py\rangle\right)^l \end{aligned}$$

*Proof.* For any $n \geq 1$, $A = (xy' + yx')/2$ with $x, y \in \mathbb{R}^r$ we have

$$(AP)^n = p_n\; AP + q_n BP \quad \text{with} \quad B = \langle x, Px\rangle\;\, yy' + \langle y, Py\rangle\;\, xx'$$

and some polynomial functions $p_n, q_n$ of maximal degree $(n-1)$ in the couple of variables $u := \langle x, Py\rangle$ and $v := \langle x, Px\rangle\,\langle y, Py\rangle$. By Cauchy Schwartz inequality we have $u^2 \leq v$.

We check this assertion by induction w.r.t. the parameter $n \geq 1$. For $n = 1$ the result is clear with $p_n = 1$ and $q_n = 0$. Assume the assertion has been checked up to rank $(n-1)$. In this case we have

$$(AP)^n = p_{n-1}\;(AP)^2 + q_{n-1} BPAP$$

Observe that

$$4\ (AP)^2 = (xy' + yx')P(xy' + yx')P = 2u\ AP + BP$$

and

$$\begin{aligned} 2\ BPAP = &\ (\langle x, Px\rangle\ yy'P + \langle y, Py\rangle\ xx'P)\,(xy'P + yx'P) \\ &= \langle x, Px\rangle\langle x, Py\rangle\ yy'P \\ &\qquad + \langle x, Px\rangle\ \langle y, Py\rangle\ (yx' + xy')P \\ &\qquad\qquad + \langle y, Py\rangle\langle x, Py\rangle\ xx'P \\ &= u\ BP + 2v\ AP \end{aligned}$$

This implies that

$$(AP)^n = p_{n-1}\ \left(\frac{u}{2}\ AP + \frac{BP}{4}\right)) + q_{n-1}\left[\frac{u}{2}\ BP + v\ AP\right]$$

from which we find the recursion

$$\begin{bmatrix} p_n \\ q_n \end{bmatrix} = \begin{bmatrix} u/2 & v \\ 1/4 & u/2 \end{bmatrix} \begin{bmatrix} p_{n-1} \\ q_{n-1} \end{bmatrix}$$

The eigenvalues of the above matrix are given by

$$\alpha_- = \frac{1}{2}\left[u - \sqrt{v}\right] \leq 0 \leq \alpha_+ = \frac{1}{2}\left[u + \sqrt{v}\right]$$

The eigenvectors $U_-$ and $U_+$ are given by

$$U_- = \begin{bmatrix} 2\ \sqrt{v} \\ -1 \end{bmatrix} \quad \text{and} \quad U_+ = \begin{bmatrix} 2\ \sqrt{v} \\ 1 \end{bmatrix}$$

This implies that

$$\begin{bmatrix} 1 \\ 0 \end{bmatrix} = \begin{bmatrix} p_1 \\ q_1 \end{bmatrix} = a_-\ U_- + a_+\ U_+ = \frac{1}{4\sqrt{v}}\ (U_- + U_+)$$

$$\begin{aligned} \Longrightarrow \begin{bmatrix} p_{n+1} \\ q_{n+1} \end{bmatrix} &= a_-\ \alpha_-^n\ U_- + a_+\ \alpha_+^n\ U_+ \\ &= \frac{1}{4\sqrt{v}} \begin{bmatrix} 2\sqrt{v}\,(\alpha_+^n + \alpha_-^n) \\ \alpha_+^n - \alpha_+^n \end{bmatrix} \end{aligned}$$

On the other hand we have

$$2^n\ (\alpha_+^n + \alpha_-^n) = \sum_{0\le k\le n} \begin{pmatrix} n \\ k \end{pmatrix} u^{n-k}\ \left(1+(-1)^k\right)\ v^{k/2}$$
$$= 2 \sum_{0\le l\le \lfloor n/2\rfloor} \begin{pmatrix} n \\ 2l \end{pmatrix} u^{n-2l}\ v^l$$

In the same vein, we have

$$2^n\ (\alpha_+^n - \alpha_-^n) = \sum_{0\le k\le n} \begin{pmatrix} n \\ k \end{pmatrix} u^{n-k}\ \left(1-(-1)^k\right)\ v^{k/2}$$

$$= 2\sqrt{v} \sum_{0\le l\le \lfloor (n-1)/2\rfloor} \begin{pmatrix} n \\ 2l+1 \end{pmatrix} u^{n-1-2l}\ v^l$$

This implies that

$$p_{n+1} = \frac{1}{2^n} \sum_{0\le l\le \lfloor n/2\rfloor} \begin{pmatrix} n \\ 2l \end{pmatrix} u^{n-2l}\ v^l$$

and

$$q_{n+1} = \frac{1}{2^{n+1}} \sum_{0\le l\le \lfloor (n-1)/2\rfloor} \begin{pmatrix} n \\ 2l+1 \end{pmatrix} u^{n-1-2l}\ v^l$$

This ends the proof of the lemma. □

We end the proof of the spectral estimate stated in (5.6) using the formula

$$\sum_{0\le l\le \lfloor n/2\rfloor} \begin{pmatrix} n \\ 2l \end{pmatrix} = 2^{n-1} = \sum_{0\le l\le \lfloor (n-1)/2\rfloor} \begin{pmatrix} n \\ 2l+1 \end{pmatrix}$$

which can be found in e.g. [70, 71]. This ends the proof of (5.6). □

## A.10 Proof of (5.7)

The proof of the sub-Gaussian estimates (5.7) is based on the following technical lemma.

**Lemma A.2.** For any symmetric matrix $A$ and any $2|t|\|AP\|_F \leq \sqrt{N}$ we have

$$\log \mathbb{E}\left(\exp\left(t \ \mathrm{Tr}\left(A\mathcal{H}_N\right)\right)\right)$$

$$\leq t^2 \ \mathrm{Tr}((AP)^2) - 2t^2 \ \|AP\|_F^2 \ \log\left(1 - 2|t| \ N^{-1/2} \ \|AP\|_F\right)$$

In addition, for any $A \geq 0$ and any $2t\mathrm{Tr}(AP) < \sqrt{N}$ we have

$$\log \mathbb{E}\left(\exp\left(t \ \mathrm{Tr}\left(A\mathcal{H}_N\right)\right)\right))$$

$$\leq t^2 \ \mathrm{Tr}((AP)^2) \left[1 - 2 \ \log\left[1 - \frac{2t}{N^{1/2}} \ \mathrm{Tr}(AP)\right]\right]$$

and

$$\log \mathbb{E}\left(\exp\left(-t \, \mathrm{Tr}\left(A\mathcal{H}_N\right)\right)\right)$$

$$\leq t^2 \, \mathrm{Tr}((AP)^2) +$$

$$\tfrac{t}{2} \sqrt{N} \, \mathrm{Tr}(AP) \ \left(1 - \tfrac{2t}{\sqrt{N}} \, \mathrm{Tr}(AP)\right) \log\left[1 - \tfrac{1}{\sqrt{N}} \left(2t\lambda_r\left(A^{1/2}PA^{1/2}\right)\right)^2\right]$$

*Proof.* Using (5.5) we find the estimate

$$\mathbb{E}\left(\exp\left(t \ \mathrm{Tr}\left(A\mathcal{H}_N\right) - t^2 \ \mathrm{Tr}((AP)^2)\right)\right)$$

$$\leq \exp\left[2t^2 \ \|AP\|_F^2 \ \sum_{n\geq 1} \frac{1}{n} \ \left(\frac{2|t|}{N^{1/2}} \ \|AP\|_F\right)^n\right]$$

$$= \exp\left[-2t^2 \ \|AP\|_F^2 \ \log\left(1 - \frac{2|t|}{N^{1/2}} \ \|AP\|_F\right)\right]$$

Whenever $A \geq 0$ and $2t\mathrm{Tr}(AP) < \sqrt{N}$, using (2.3) we have

$$\mathbb{E}\left(\exp\left(t \, \mathrm{Tr}\left(A\mathcal{H}_N\right) - t^2 \, \mathrm{Tr}((AP)^2)\right)\right) \geq 1$$

as well as

$$\exp\left[\sum_{n\geq 3} \frac{t^n}{n} \frac{2^{n-1}}{N^{n/2-1}} \ \mathrm{Tr}(AP)^n\right]$$

$$\leq \ \exp\left[-2t^2 \ \mathrm{Tr}((AP)^2) \log\left(1 - \frac{2t}{N^{1/2}} \ \mathrm{Tr}(AP)\right)\right]$$

On the other hand, we have the decomposition

$$\sum_{m\geq 1} \frac{(-t)^{m+2}}{m+2} \frac{1}{N^{m/2}} 2^{m+1} \operatorname{Tr}((AP)^{m+2})$$

$$= \frac{2t}{\sqrt{N}} \sum_{m\geq 1} \frac{t^{2m+1}}{2m+2} \frac{1}{N^{m-1/2}} 2^{2m} \operatorname{Tr}((AP)^{2m+2})$$

$$- \sum_{m\geq 1} \frac{t^{2m+1}}{2m+1} \frac{1}{N^{m-1/2}} 2^{2m} \operatorname{Tr}((AP)^{2m+1})$$

$$\leq -\left(1 - \frac{2t}{\sqrt{N}} \operatorname{Tr}(AP)\right) \sum_{m\geq 1} \frac{t^{2m+1}}{2m+1} \frac{1}{N^{m-1/2}} 2^{2m} \operatorname{Tr}((AP)^{2m+1})$$

This yields the estimate

$$\mathbb{E}\left(\exp\left(-t \operatorname{Tr}\left(A\mathcal{H}_N\right) - t^2 \operatorname{Tr}((AP)^2)\right)\right)$$

$$\leq \exp\left[-\frac{1}{\sqrt{N}}\left(1 - \frac{2t}{\sqrt{N}} \operatorname{Tr}(AP)\right)\right.$$

$$\left.\times \sum_{m\geq 1} \frac{t^{2m+1}}{2m+1} \frac{1}{N^{m-1}} 2^{2m} \operatorname{Tr}((AP)^{2m+1})\right]$$

$$\leq \exp\left[-\frac{t}{2} \sqrt{N} \operatorname{Tr}(AP) \left(1 - \frac{2t}{\sqrt{N}} \operatorname{Tr}(AP)\right)\right.$$

$$\left.\times \sum_{m\geq 1} \frac{(2t\lambda_r\left(A^{1/2}PA^{1/2}\right))^{2m}}{m} \frac{1}{N^m}\right]$$

as soon as $2t\operatorname{Tr}(AP) < \sqrt{N}$. This ends the proof of the lemma. □

Combining lemma A.2 with the estimate $-\log(1-x) \leq \frac{x}{1-x} \leq 2x \leq 1$, which is valid for any $0 \leq x \leq 1/2$, we obtain the sub-Gaussian property (5.7). □

# References

[1] N.I. Akhiezer. The Classical Moment Problem and Some Related Questions in Analysis. Oliver & Boyd (1965).

[2] F.P. Anaraki and S. Hughes. Memory and computation efficient PCA via very sparse random projections. In International Conference on Machine Learning, pp. 1341–1349 (2014).

[3] T.W. Anderson. An Introduction to Multivariate Statistical Analysis. 3rd Edition. John Wiley & Sons, New York (2003).

[4] G.W. Anderson, A. Guionnet, and O. Zeitouni. An Introduction to Random Matrices. Cambridge University Press (2010).

[5] G. Ben Arous and A. Guionnet. Large deviations for Wigner's law and Voiculescu's non-commutative entropy. Probability Theory and Related Fields. vol. 108, no. 4. pp. 517–542 (1997).

[6] F.R. Bernhart. Catalan, Motzkin, and Riordan numbers. Discrete Mathematics. vol. 204, no. 1–3. pp. 73–112 (1999).

[7] E. Bingham and H. Mannila. Random projection in dimensionality reduction: applications to image and text data. In Proceedings of the seventh ACM SIGKDD international conference on Knowledge discovery and data mining, pp. 245–250. ACM (2001).

[8] C.M. Bishop. Pattern Recognition and Machine Learning. Springer (2006).

[9] A.N. Bishop and P. Del Moral. Matrix product moments in normal variables. arXiv e-print, arXiv:1808.00235 (2018).

[10] A.N. Bishop and P. Del Moral. On the Stability of Matrix-Valued Riccati Diffusions. arXiv e-print, arXiv:1703.00353 (2017).

[11] A.N. Bishop, P. Del Moral, and A. Niclas. A perturbation analysis of stochastic matrix Riccati diffusions. arXiv e-print, arXiv:1709.05071 (2017).

[12] M. Bóna and I. Mezö. Real zeros and partitions without singleton blocks. European Journal of Combinatorics. vol. 51. pp. 500–510 (2016).

[13] J.P. Bouchaud and M. Potters. Theory of Financial Risk and Derivative Pricing From Statistical Physics to Risk Management. Cambridge University Press (2009).

[14] M.W. Browne. Asymptotically distribution-free methods for the analysis of covariance structures. British Journal of Mathematical and Statistical Psychology. vol. 37, no. 1. pp. 62–83 (1984).

[15] M.F. Bru. Wishart Processes. Journal of Theoretical Probability. vol. 4, no. 4. pp. 725–751 (1991).

[16] T. Carleman. Les fonctions quasi analytiques. Collection Borel, Gauthier–Villars, Paris (1926).

[17] X. Chen and W. Chu. Moments on Catalan numbers. Journal of Mathematical Analysis and Applications. vol. 349, no. 2. pp. 311–316 (2009).

[18] S.I. Choi. The inversion formula of the Stieltjes transform of spectral distribution. Journal of the Chungcheong Mathematical Society. vol. 22, no. 3. pp. 519–524 (2009).

[19] J.C. Cox, J.E. Ingersoll, and S.A. Ross. A Theory of the Term Structure of Interest Rate. Econometrica. vol. 53, no. 2. pp. 385–407 (1985).

[20] C. Cuchiero, D. Filipovich, E. Mayerhofer, and J. Teichman. Affine processes on positive semidefinite matrices. The Annals of Applied Probability. vol. 21, no. 2. pp. 397–463 (2011)

[21] P. Del Moral. Feynman-Kac formulae: Genealogical and interacting particle systems with applications. Springer, New York (2004).

[22] P. Del Moral. Mean field simulation for Monte Carlo integration. Chapman & Hall/CRC Monographs on Statistics & Applied Probability (2013).

[23] P. Del Moral and J. Tugaut. On the stability and the uniform propagation of chaos properties of ensemble Kalman-Bucy filters. Annals of Applied Probability. vol. 28, no. 2. pp. 790–850 (2018).

[24] N. Dershowitz and S. Zaks. Ordered trees and non-crossing partitions. Discrete Mathematics. vol. 62, no. 1. pp. 215–218 (1986).

[25] T.A. Dowling. Stirling numbers. In Applications of Discrete Mathematics (eds: J.G. Michaels and K.H. Rosen). Chapter 6. McGraw-Hill (1991).

[26] A. Edleman. Eigenvalues and condition numbers of random matrices. Ph.D Thesis, MIT (1989).

[27] E.J. Elton and M.J. Gruber. Modern Portfolio Theory and Investment Analysis. John Wiley and Sons, New York (1995).

[28] F. Götze and A. Tikhomirov. Rate of convergence in probability to the Marchenko–Pastur law. Bernoulli. vol. 10, no. 3, pp. 503–548 (2004).

[29] P. Graczyk, G. Letac, and H. Massam. The hyperoctahedral group, symmetric group representations and the moments of the real Wishart distribution. Journal of Theoretical Probability. vol. 18, no. 1. pp. 1–42 (2005).

[30] A. Guionnet. Large Random Matrices: Lectures on Macroscopic Asymptotics. Lecture Notes in Mathematics. Springer–Verlag (2009).

[31] Harish-Chandra. Differential operators on a semi-simple Lie algebra. American Journal of Mathematics. vol. 79, no.1. pp. 87—120 (1957).

[32] H. Epstein. Remarks on Two Theorems of E. Lieb, Comm. Math. Phys., vol. 31, pp. 317–325 (1973).

[33] L. Isserlis. On a formula for the product-moment coefficient of any order of a normal frequency distribution in any number of variables. Biometrika, vol. 12, no. 1/2. pp. 134–139 (1918).

[34] O. James, H-N. Lee. Concise Probability Distributions of Eigenvalues of Real-Valued Wishart Matrices. arXiv e-print, arXiv:1402.6757 (2014).

[35] R.D. Jensen. The joint distribution of traces of Wishart matrices and some applications. The Annals of Mathematical Statistics. vol. 41, no. 1. pp. 133–145 (1970).

[36] M. Junge and Q. Xu. On the best constants in some non-commutative martingale inequalities. Bulletin of the London Mathematical Society. vol. 37, no. 2. pp. 243–253 (2005).

[37] M. Junge and Q. Xu. Noncommutative Burkholder/Rosenthal inequalities II: Applications. Israel Journal of Mathematics. vol. 167, no 1. pp. 227–282 (2008).

[38] R. Kan. From Moments of Sums to Moments of Product. Journal of Multivariate Analysis, vol. 99, no. 3. pp. 542–554 (2008).

[39] N.D. Kazarinoff. Geometric inequalities. Random House, New York (1961).

[40] M. Keller-Ressel, W. Schachermayer, and J. Teichmann. Affine processes are regular. Probability Theory and Related Fields. vol. 151. no. 3-4. pp. 591–611 (2011).

[41] M.G. Kendall and A. Stuart. The Advanced Theory of Statistics. C. Griffin & Company, London (1943).

[42] J.S. Kim, K.H. Lee, and S.J. Ho. Weight multiplicities and Young tableaux through affine crystals. arXiv e-print, arXiv:1703.10321 (2017).

[43] I. Kortchemski and C. Marzouk. Simply generated non-crossing partitions. Combinatorics, Probability and Computing. vol. 26, no. 4. pp. 560–592 (2017).

[44] G. Kreweras. Sur les partitions non croisées d'un cycle. Discrete Mathematics. vol. 1, no. 4. pp. 333–350 (1972).

[45] P.R. Krishnaiah and T.C. Chang. On the exact distribution of the smallest root of the Wishart matrix using zonal polynomials. Annals of the Institute of Statistical Mathematics. vol. 23, no. 11. pp. 293–295 (1971).

[46] M.J. Kusner, N. I. Kolkin, S. Tyree, and K.Q. Weinberger. Image Data Compression for Covariance and Histogram Descriptors. arXiv preprint arXiv:1412.1740 (2014).

[47] G. Letac and H. Massam. All invariant moments of the Wishart distribution. Scandinavian Journal of Statistics. vol. 31, no. 2. pp. 295–318 (2004).

[48] G. Letac and H. Massam. The noncentral Wishart as an exponential family, and its moments. Journal of Multivariate Analysis. vol. 99. no. 7. pp. 1393–1417 (2008).

[49] G. Letac and H. Massam. A tutorial on non central Wishart distributions. Technical Report: Lab. Stat. Prob., Toulouse, France (2004).

[50] P.L. Leung. An identity for the noncentral wishart distribution with application. Journal of Multivariate Analysis. vol. 48, no. 1. pp. 107–114 (1994)

[51] S.C. Liaw, H.G. Yeh, F.K. Hwang, and G.J. Chang. A simple and direct derivation for the number of noncrossing partitions. Proceedings of the American Mathematical Society. vol. 126, no. 6. pp. 1579–1581 (1998).

[52] E.H. Lieb. Convex Trace Functions and the Wigner–Yanase–Dyson Conjecture. Advances in Math. vol. 11. pp. 267–288 (1973).

[53] G. Lin and J. Stoyanov. On the Moment Determinacy of Products of Non-identically Distributed Random Variables. arXiv e-print, arXiv:1406.1654 (2014).

[54] F. Lust-Piquard. Inégalités de Khintchine dans $C_p$ $(1 < p < \infty)$. C.R. Acad. Sci. Paris. Sér. I Math. vol. 303, no. 7. pp. 289–292 (1986).

[55] F. Lust-Piquard and G. Pisier. Noncommutative Khintchine and Paley inequalities. Arkiv för Matematik. vol. 29, no. 1–2. pp. 241–260 (1991).

[56] V.A. Marchenko and L. A. Pastur. Distribution of eigenvalues for some sets of random matrices. Matematicheskii Sbornik. vol. 114, no. 4. pp. 507–536 (1967).

[57] H. Markowitz. Portfolio Selection: Efficient Diversification of Investments. John Wiley and Sons, New York (1959).

[58] E. Mayerhofer, O. Pfaffel, and R. Stelzer. On strong solutions for positive definite jump-diffusions. Stochastic Processes and Their Applications. vol. 121, no. 9. pp. 2072–2086 (2011).

[59] M.L. Mehta. Random Matrices. 3rd Edition. Elsevier/Academic Press, Amsterdam (2004).

[60] J.A. Mingo and R. Speicher. Free Probability and Random Matrices. Springer (2016).

[61] L. Mirsky. Symmetric gauge functions and unitarily invariant norms. The Quarterly Journal of Mathematics. vol. 11, no. 1. pp. 50–59 (1960).

[62] R.J. Muirhead. Aspects of Multivariate Statistical Theory. Wiley, New York (1982).

[63] K.P. Murphy. Machine Learning, A Probabilistic Perspective. Adaptive Computation and Machine Learning, The MIT Press (2012).

[64] K.P. Murphy. Conjugate Bayesian analysis of the Gaussian distribution. Technical report (2007).

[65] A. Nica and R. Speicher. Lectures on the combinatorics of free probability. Cambridge University Press (2006).

[66] A. Nica, D. Shlyakhtenko, and R. Speicher. Operator-valued distributions. I. Characterizations of freeness, International Mathematics Research Notices, pp. 1509–1538 (2002).

[67] H.J. Park and E. Ayanoglu. An upper bound to the marginal PDF of the ordered eigenvalues of Wishart matrices and its application to MIMO diversity analysis. In Proc. of the 2010 IEEE International Conference on Communications, Cape Town, South Africa (July, 2010).

[68] K. Phillips. R-Functions to Symbolically Compute the Central Moments of the Multivariate Normal Distribution. Journal of Statistical Software, Code Snippets. vol. 33, no. 1. pp. 1–14 (2010).

[69] F. Pourkamali-Anaraki. Estimation of the sample covariance matrix from compressive measurements. IET Signal Processing, vol. 10, no. 9, 1089–1095 (2016).

[70] J. Quaintance and H.W. Gould. Combinatorial Identities for Stirling Numbers. The unpublished notes of H.W. Gould. World Scientific Publishing, Singapore (2016).

[71] J. Quaintance and H.W. Gould. Table for Fundamentals of Series: Table III: Basic Algebraic Techniques. From the seven unpublished manuscripts of H.W. Gould, edited and compiled by J. Quaintance (2010).

[72] C.E. Rasmussen and C.K.I. Williams. Gaussian Processes for Machine Learning. The MIT Press (2006).

[73] C. Redelmeier and I. Emily. Genus expansion for real Wishart matrices. J. Theoret. Probab., vol. 24, pp. 1044–1062 (2011).

[74] J. Riordan. Enumeration of plane trees by branches and endpoints. Journal of Combinatorial Theory, Series A. vol. 19, no. 2. pp. 214–222 (1975).

[75] M. Rudelson. Random vectors in the isotropic position. Journal of Functional Analysis. vol. 164, no. 1. pp. 60–72 (1999).

[76] L.W. Shapiro. A Catalan triangle. Discrete Math. vol. 14, no. 1. pp. 83–90 (1976).

[77] H. Shin, M.Z. Win, J.H. Lee, and M. Chiani. On the capacity of doubly correlated MIMO channels. IEEE Transactions on Wireless Communications. vol. 5, no. 8. pp. 2253–2265 (2006).

[78] A.N. Shiryaev. Probability. Springer-Verlag, New York (1996).

[79] J. Stoyanov and G. Lin. Hardy's condition in the moment problem for probability distributions. Theory of Probability & Its Applications. vol. 57, no. 4. pp. 699–708 (2013).

[80] J.W. Silverstein and P.L. Combettes. Signal detection via spectral theory of large dimensional random matrices. IEEE Transactions on Signal Processing. vol. 40, no. 8. pp. 2100–2105 (1992).

[81] R. Simion. Noncrossing partitions. Discrete Mathematics. vol. 217, no. 1–3. pp. 367–409 (2000).

[82] R. Simion. Combinatorial statistics on non-crossing partitions. Journal of Combinatorial Theory, Series A. vol. 66, no. 2. pp. 270–301 (1994).

[83] R. Simion and D. Ullman. On the structure of the lattice of noncrossing partitions. Discrete Mathematics. vol. 98, no. 3. pp. 193–206 (1991).

[84] N.J.A. Sloane. The on-line encyclopaedia of integer sequences: https://oeis.org/.

[85] R.P. Stanley. Enumerative Combinatorics, vol. 1 and 2. Cambridge University Press (1999).

[86] R. Speicher. Combinatorial theory of the free product with amalgamation and operator-valued free probability theory, Mem. Amer. Math. Soc., vol. 132, no. 627, pp. x+88 (1998).

[87] R. Speicher and C. Vargas. Free deterministic equivalents, rectangular random matrix models, and operator- valued free probability theory Random Matrices: Theory and Applications (2012).

[88] T. Sugiyama. Percentile points of the largest latent root of a matrix and power calculation for testing $\Sigma = I$. Journal of the Japan Statistical Society, Japanese Issue. vol. 3, no. 1. pp. 1–8 (1973).

[89] T. Tao. Topics in random matrix theory. American Mathematical Society, Providence, Rhode Island (2012).

[90] S. Thorbjornsen. Mixed moments of Voiculescu's Gaussian random matrices, Journal of Functional Analysis, vol. 176, pp. 213–246 (2000).

[91] J.Y. Tourneret, A. Ferrari, and G. Letac. The noncentral Wishart distribution: Properties and application to speckle imaging. In Proceedings of the 2005 IEEE/SP 13th Workshop on Statistical Signal Processing, Bordeaux, France (July, 2005).

[92] D.S. Tracy and S.A Sultan. Higher Order Moments of Multivariate Normal Distribution Using Matrix Derivatives. Stochastic Analysis and Applications. vol. 11, no. 3. pp. 337–348 (1993).

[93] J.A. Tropp. User-Friendly Tail Bounds for Sums of Random Matrices. Foundations of Computational Mathematics. vol. 12, no. 4. pp. 389–434 (2012).

[94] J.A. Tropp. An Introduction to Matrix Concentration Inequalities. Foundations and Trends in Machine Learning. vol. 8, no. 1-2. pp. 1–230 (2015).

[95] D.V. Voiculescu. Operations on certain non-commutative operator-valued random variables, Asterisque, no. 232, pp. 243–275 (1995).

[96] H. Weyl. The laws of asymptotic distribution of the eigenvalues of linear partial differential equation. Math. Ann. vol. 71. pp. 441–479 (1912).

[97] G.C. Wick. The evaluation of the collision matrix. Physical Review. vol. 80, no. 2, pp. 268–272 (1950).

[98] E. Wigner. Characteristic Vectors of Bordered Matrices With Infinite Dimensions. The Annals of Mathematics. vol. 62, no. 3. pp. 548–564 (1955).

[99] J.H. Winters. On the capacity of radio communication systems with diversity in Rayleigh fading environment. IEEE Journal on Selected Areas in Communications. vol. 5, no. 5. pp. 871–878 (1987).

[100] J. Wishart. The generalized product moment distribution in samples from a normal multivariate population. Biometrika, vol. 20A, no. 1/2. pp. 32–52 (1928).

[101] F. Yano and H. Yoshida. Some set partition statistics in non-crossing partitions and generating functions. Discrete Mathematics. vol. 307, no. 24. pp. 3147–3160 (2007).

[102] Y.Q. Yin. Limiting spectral distribution for a class of random matrices. Journal of Multivariate Analysis. vol. 20, no. 1. pp. 50–68 (1986).

[103] Y.Q. Yin and P.R. Krishnaiah. A limit theorem for the eigenvalues of product of two random variables. Journal of Multivariate Analysis. vol. 13, no. 4. pp. 489–507 (1983).

[104] A. Zanella, M. Chiani, and M.Z. Win. On the marginal distribution of the eigenvalues of Wishart matrices. IEEE Transactions on Communications. vol. 57, no. 4. pp. 1050–1060 (2009).